\documentclass[12equipped]{article}

\usepackage{}
\usepackage{amssymb}
\usepackage{amsfonts}

\usepackage{amssymb,amsmath,tikz}
\usepackage{graphicx}
\usepackage[plainpages,backref]{hyperref}
\usepackage{enumerate}
\usepackage{verbatim}
\usepackage{mathrsfs}
\usepackage{latexsym}
\newcommand {\Omit}[1]{}

\usepackage [autostyle, english = american]{csquotes}

\input epsf.tex

\usepackage[all]{xy}\usepackage{}

\usepackage{tikz}
\usetikzlibrary{automata,positioning}
\usepackage[top=1.5in, bottom=1.5in, left=1.5in, right=1.5in]{geometry}

\usetikzlibrary{shapes,shadows}

\usetikzlibrary{arrows,calc}
\tikzset{
>=stealth',
help lines/.style={dashed, thick},
axis/.style={<->},
important line/.style={thick},
connection/.style={thick, dotted},
}
\usetikzlibrary{patterns}
\newlength{\hatchspread}
\newlength{\hatchthickness}
\tikzset{hatchspread/.code={\setlength{\hatchspread}{#1}},
         hatchthickness/.code={\setlength{\hatchthickness}{#1}}}
\tikzset{hatchspread=3pt,
         hatchthickness=0.4pt}
\pgfdeclarepatternformonly[\hatchspread,\hatchthickness]
   {custom north west lines}
   {\pgfqpoint{-2\hatchthickness}{-2\hatchthickness}}
   {\pgfqpoint{\dimexpr\hatchspread+2\hatchthickness}{\dimexpr\hatchspread+2\hatchthickness}}
   {\pgfqpoint{\hatchspread}{\hatchspread}}
   {
    \pgfsetlinewidth{\hatchthickness}
    \pgfpathmoveto{\pgfqpoint{0pt}{\hatchspread}}
    \pgfpathlineto{\pgfqpoint{\dimexpr\hatchspread+0.15pt}{-0.15pt}}
    \pgfusepath{stroke}
   }

\newcommand{\nc}{\newcommand}
\nc{\rnc}{\renewcommand}
\nc{\bb}[1]{{\mathbb #1}}
\nc{\bbA}{\bb{A}}\nc{\bbB}{\bb{B}}
\nc{\bbD}{\bb{D}}
\nc{\bbE}{\bb{E}}\nc{\bbF}{\bb{F}}\nc{\bbG}{\bb{G}}\nc{\bbH}{\bb{H}}
\nc{\bbI}{\bb{I}}\nc{\bbJ}{\bb{J}}\nc{\bbK}{\bb{K}}\nc{\bbL}{\bb{L}}
\nc{\bbM}{\bb{M}}\nc{\bbN}{{\bf N}}\nc{\bbO}{\bb{O}}
\nc{\bbQ}{{\bf Q}}\nc{\bbR}{\bb{R}}\nc{\bbS}{\bb{S}}\nc{\bbT}{\bb{T}}
\nc{\bbU}{\bb{U}}\nc{\bbV}{\bb{V}}\nc{\bbW}{\bb{W}}\nc{\bbX}{\bb{X}}
\nc{\bbY}{\bb{Y}}\nc{\bbZ}{\bb{Z}}
\nc{\bbP}{\mathbf{P}}
\nc{\mbf}[1]{{\mathbf #1}}
\nc{\bfA}{\mbf{A}}\nc{\bfB}{\mbf{B}}\nc{\bfC}{\mbf{C}}\nc{\bfD}{\mbf{D}}
\nc{\bfE}{\mbf{E}}\nc{\bfF}{\mbf{F}}\nc{\bfG}{\mbf{G}}\nc{\bfH}{\mbf{H}}
\nc{\bfI}{\mbf{I}}\nc{\bfJ}{\mbf{J}}\nc{\bfK}{\mbf{K}}\nc{\bfL}{\mbf{L}}
\nc{\bfM}{\mbf{M}}\nc{\bfN}{\mbf{N}}\nc{\bfO}{\mbf{O}}\nc{\bfP}{\mbf{P}}
\nc{\bfQ}{\mbf{Q}}\nc{\bfR}{\mbf{R}}\nc{\bfS}{\mbf{S}}\nc{\bfT}{\mbf{T}}
\nc{\bfU}{\mbf{U}}\nc{\bfV}{\mbf{V}}\nc{\bfW}{\mbf{W}}\nc{\bfX}{\mbf{X}}
\nc{\bfY}{\mbf{Y}}\nc{\bfZ}{\mbf{Z}}
\nc{\bfa}{\mbf{a}}\nc{\bfb}{\mbf{b}}\nc{\bfc}{\mbf{c}}\nc{\bfd}{\mbf{d}}
\nc{\bfe}{\mbf{e}}\nc{\bff}{\mbf{f}}\nc{\bfg}{\mbf{g}}\nc{\bfh}{\mbf{h}}
\nc{\bfi}{\mbf{i}}\nc{\bfj}{\mbf{j}}\nc{\bfk}{\mbf{k}}\nc{\bfl}{\mbf{l}}
\nc{\bfm}{\mbf{m}}\nc{\bfn}{\mbf{n}}\nc{\bfo}{\mbf{o}}\nc{\bfp}{\mbf{p}}
\nc{\bfq}{\mbf{q}}\nc{\bfr}{\mbf{r}}\nc{\bfs}{\mbf{s}}\nc{\bft}{\mbf{t}}
\nc{\bfu}{\mbf{u}}\nc{\bfv}{\mbf{v}}\nc{\bfw}{\mbf{w}}\nc{\bfx}{\mbf{x}}
\nc{\bfy}{\mbf{y}}\nc{\bfz}{\mbf{z}}

\newcommand{\op}{\text{op}}

\nc{\mcal}[1]{{\mathcal #1}}
\nc{\calA}{\mcal{A}}\nc{\calB}{\mcal{B}}\nc{\calC}{\mcal{C}}\nc{\calD}{\mcal{D}}
\nc{\calE}{\mcal{E}} \nc{\calF}{\mcal{F}}\nc{\calG}{\mcal{G}}\nc{\calH}{\mcal{H}}
\nc{\calI}{\mcal{I}}\nc{\calJ}{\mcal{J}}\nc{\calK}{\mcal{K}}\nc{\calL}{\mcal{L}}
\nc{\calM}{\mcal{M}}\nc{\calN}{\mcal{N}}\nc{\calO}{\mcal{O}}\nc{\calP}{\mcal{P}}
\nc{\calQ}{\mcal{Q}}\nc{\calR}{\mcal{R}}\nc{\calS}{\mcal{S}}\nc{\calT}{\mcal{T}}
\nc{\calU}{\mcal{U}}\nc{\calV}{\mcal{V}}\nc{\calW}{\mcal{W}}\nc{\calX}{\mcal{X}}
\nc{\calY}{\mcal{Y}}\nc{\calZ}{\mcal{Z}}
\nc{\fA}{\frak{A}}\nc{\fB}{\frak{B}}\nc{\fC}{\frak{C}} \nc{\fD}{\frak{D}}
\nc{\fE}{\frak{E}}\nc{\fF}{\frak{F}}\nc{\fG}{\frak{G}}\nc{\fH}{\frak{H}}
\nc{\fI}{\frak{I}}\nc{\fJ}{\frak{J}}\nc{\fK}{\frak{K}}\nc{\fL}{\frak{L}}\nc{\fM}{\frak{M}}
\nc{\fN}{\frak{N}}\nc{\fO}{\frak{O}}\nc{\fP}{\frak{P}}
\nc{\fQ}{\frak{Q}}\nc{\fR}{\frak{R}}\nc{\fS}{\frak{S}}\nc{\fT}{\frak{T}}
\nc{\fU}{\frak{U}}\nc{\fV}{\frak{V}}\nc{\fW}{\frak{W}}\nc{\fX}{\frak{X}}
\nc{\fY}{\frak{Y}}\nc{\fZ}{\frak{Z}}
\nc{\fa}{\frak{a}}\nc{\fb}{\frak{b}}\nc{\fc}{\frak{c}} \nc{\fd}{\frak{d}}
\nc{\fe}{\frak{e}}\nc{\fFf}{\frak{f}}\nc{\fg}{\frak{g}}\nc{\fh}{\frak{h}}
\nc{\fri}{\frak{i}}\nc{\fj}{\frak{j}}\nc{\fk}{\frak{k}}\nc{\fl}{\frak{l}}
\nc{\fn}{\frak{n}}\nc{\fo}{\frak{o}}\nc{\fp}{\frak{p}}
\nc{\fq}{\frak{q}}\nc{\fr}{\frak{r}}\nc{\fs}{\frak{s}}\nc{\ft}{\frak{t}}
\nc{\fu}{\frak{u}}\nc{\fv}{\frak{v}}\nc{\fw}{\frak{w}}\nc{\fx}{\frak{x}}
\nc{\fy}{\frak{y}}\nc{\fz}{\frak{z}}

\DeclareMathOperator{\Image}{Im} \DeclareMathOperator{\Sym}{Sym}

 \DeclareMathOperator{\GL}{GL}

\DeclareMathOperator{\Hom}{{Hom}} 
\DeclareMathOperator{\Ext}{{Ext}}
\DeclareMathOperator{\sHom}{{\mathscr{H}om}}

\DeclareMathOperator{\Hilb}{{Hilb}}
\DeclareMathOperator{\Tot}{{Tot}}

 \DeclareMathOperator{\Lie}{Lie}
 \DeclareMathOperator{\tr}{tr}

\DeclareMathOperator{\Grass}{Grass} \DeclareMathOperator{\End}{End}

\DeclareMathOperator{\Coh}{Coh}
\DeclareMathOperator{\Mod}{Mod\hbox{-}}
\DeclareMathOperator{\Modf}{Mod_f\hbox{-}}

 \DeclareMathOperator{\SL}{SL}
  \DeclareMathOperator{\Fl}{Fl}

      \DeclareMathOperator{\crit}{crit}
      \DeclareMathOperator{\Surj}{Surj}

\DeclareMathOperator{\SH}{\textbf{SH}}

\DeclareMathOperator{\Crit}{Crit}

\newcommand{\g}{\mathfrak{g}}
\newcommand{\h}{\mathfrak{h}}

\DeclareMathOperator{\Rep}{Rep}
\DeclareMathOperator{\Res}{Res}

\DeclareMathOperator{\Sh}{Sh}

\newcommand{\surj}{\twoheadrightarrow}
\newcommand{\inj}{\hookrightarrow}

\newcommand{\pt}{\text{pt}}

\newcommand{\PP}{\bbP}

\newcommand{\N}{\bbN}

\newcommand{\affY}{Y_{\hbar_1, \hbar_2, \hbar_3}(\widehat{\mathfrak{gl}(1)})}

\DeclareMathOperator{\fac}{fac}

 \newcommand{\ext}{\fe}

\newcommand{\loc}{loc}
\DeclareMathOperator{\BM}{BM}

\DeclareMathOperator{\coop}{coop}

\DeclareMathOperator{\Perv}{Perv}

\newcount\cols
{\catcode`,=\active\catcode`|=\active
 \gdef\Young(#1){\hbox{$\vcenter
 {\mathcode`,="8000\mathcode`|="8000
  \def,{\global\advance\cols by 1 &}%
  \def|{\cr
        \multispan{\the\cols}\hrulefill\cr
        &\global\cols=2 }%
  \offinterlineskip\everycr{}\tabskip=0pt
  \dimen0=\ht\strutbox \advance\dimen0 by \dp\strutbox
  \halign
   {\vrule height \ht\strutbox depth \dp\strutbox##
    &&\hbox to \dimen0{\hss$##$\hss}\vrule\cr
    \noalign{\hrule}&\global\cols=2 #1\crcr
    \multispan{\the\cols}\hrulefill\cr%
   }
 }$}}
}

\newtheorem{thm}{Theorem}[subsection]
\newtheorem{defn}[thm]{Definition}
\newtheorem{lmm}[thm]{Lemma}
\newtheorem{rmk}[thm]{Remark}
\newtheorem{prp}[thm]{Proposition}
\newtheorem{conj}[thm]{Conjecture}
\newtheorem{exa}[thm]{Example}

\newcommand{\C}{{\bf C}}
\newcommand{\R}{{\bf R}}
\newcommand{\Z}{{\bf Z}}

\renewcommand{\H}{{\mathcal{H}}}

\title{Cohomological Hall algebras and perverse coherent sheaves on toric Calabi-Yau 3-folds}

\author {Miroslav Rap\v{c}\'{a}k, Yan Soibelman, Yaping Yang, Gufang Zhao}

\begin{document}
\maketitle

\begin{abstract}
We study the Drinfeld double of the (equivariant spherical) Cohomological Hall algebra   in the sense of Kontsevich and Soibelman, associated 
to a smooth toric Calabi-Yau 3-fold $X$. 
By general reasons, the COHA acts on the cohomology of the moduli spaces of certain perverse coherent systems on $X$ via ``raising operators". 
Conjecturally the COHA action extends to an action of the Drinfeld double by adding the ``lowering operators". 

In this paper, we show that the Drinfeld double is  a generalization of the notion of the Cartan doubled Yangian defined earlier by Finkelberg and others. 
We extend this ``$3d$ Calabi-Yau perspective" on the Lie theory furthermore by associating a root system to certain families of $X$.
We formulate a conjecture that the above-mentioned action of the Drinfeld double factors through a shifted Yangian of the root system.
 The shift is explicitly determined by the moduli problem and the choice of stability conditions, and is expressed explicitly in terms of an intersection number in $X$. We check the conjectures in several examples, including a special case of an earlier  conjecture of Costello.
\end{abstract}
\tableofcontents

\section{Introduction}
\label{intro}

The paper is devoted to  Cohomological Hall algebras (COHA) of toric Calabi-Yau 3-folds and their representations. The notion of Cohomological Hall algebras was introduced by Kontsevich and Soibelman in \cite{KS} as an alternative way to define motivic Donaldson-Thomas invariants  of $3$-dimensional Calabi-Yau categories introduced in their earlier paper \cite{KS2}.  The framework for COHA in the loc.cit. was algebraic, namely, Quillen-smooth algebras endowed with potential. In this approach the relation to the geometry of $3$-dimensional Calabi-Yau manifolds is not direct. In order to relate  geometry to algebra one has to choose a generator of the derived category of coherent sheaves, and then realize the latter as the category of dg-modules over the endomorphism algebra of the generator.  Roughly, the  representations appear as the action of the algebra $Hom(P, P)$ on $Hom(P, V)$, where $P$ is an object of our category (the latter can be of algebraic or geometric origin) such that the COHA of the category is $Hom(P,P)$ and $V$ is another object of the category. Details depend on a particular problem. In applications to physics the corresponding classes of representations can appear e.g. in the form of the action of the algebra of closed BPS states on the space of open BPS states, while in mathematics they can appear e.g. in the form of  the action of the cohomology of the stack of $0$-dimensional torsion sheaves on a Calabi-Yau $3$-fold on the cohomology of the moduli space of stable framed sheaves supported on curves and/or surfaces (see \cite{So} for discussion).  

It was already pointed out in \cite{KS} that COHA resembles the quantized enveloping algebra of a (graded) Lie algebra. The existence of such ``BPS Lie algebra" conjectured in the loc.cit, was constructed later in \cite{DM}, \cite{D2}.  Thus one can think of COHA as a sort of quantized enveloping algebra of a $3$-dimensional analog of a generalized Kac-Moody Lie algebra and treat its representations in a similar fashion. In particular, one can think of every chamber in the space of stability structures on the category as a choice of Borel subalgebra in COHA  and hence hope for a sort of highest weight representation theory. In a sense our paper provides a realization of the above ideas in the case of a  class of toric Calabi-Yau $3$-folds.

More precisely, let  $X$ be a smooth toric Calabi-Yau 3-fold which is a resolution of singularities $f: X\to Y$, where $Y$ is affine. Assume that the fibers have dimension at most 1 and each component of the exceptional fiber is  $\PP^1$. Let $D^b\Coh(X)$ be the bounded derived category of coherent sheaves on $X$. The map $f: X\to Y$, by general theory of Bridgeland and Van den Bergh \cite{Br,vdB}, provides $D^b\Coh(X)$ with a $t$-structure whose heart is the module category $\Mod\calA_0$.
When $\calA_0$ is the Jacobian algebra of a quiver with potential $(Q,W)$, the general construction of Kontsevich and Soibelman \cite{KS} yields a COHA of the pair $(Q,W)$. We  denote it by $\H_X$ and call it the COHA associated to this $t$-structure. By general reasons $\H_X$ acts on the 
cohomology of the moduli spaces of certain sheaves (e.g. torsion free sheaves, or sheaves supported on a toric divisor) on $X$ via  ``raising operators". In a technical sense, the goal of this paper is to construct and study various ``doubles" of $\H_X$ as well as their representations on the cohomology of these moduli spaces of sheaves by adding the ``lowering operators". We also consider other versions of COHA, including equivariant and spherical ones.

Some basic properties of the Drinfeld double of  the (equivariant spherical) COHA are summarized below \S~\ref{subsec:intr_double}.
From the perspective of the above-mentioned analogy with quantized Kac-Moody algebras, we expect that there exists the notion of ``generalized root system" associated to the above data of toric manifold. Then   the Drinfeld double is characterized algebraically as the so-called Cartan doubled Yangian of this root system. The root system is a subset in the topological $K$-group
of the category $\Mod\calA_0$ (lattice $\Gamma$ in the notation of \cite{KS2}) endowed with a bilinear form. 
There is a notion of simple root (in the case of quivers with potential they correspond to the vertices of the quiver),
the notion of real and imaginary roots.  {In particular, when $X$ is the toric Calabi-Yau resolution of the singularity given by equation $xy = z^mw^n$, one obtains the root system of the affine Lie super algebra $\widehat{\mathfrak{gl}(m|n)}$}. 
We will discuss generalized root systems in more detail
in the future publications, but we will use the above terminology without further comments in this paper.

Geometrically, the Drinfeld double of COHA is expected to act on the cohomology of the moduli space of perverse coherent systems of $X$ in the sense of \cite{Nagao,NN}. The action  factors through the one of the so-called shifted Yangian{, which is a generalization to our version of root system of the notion studied earlier by Finkelberg and others \cite{Finetal}. In this paper, we give a prediction for the shift  in the shifted Yangian in terms of a certain intersection number (see \S~\ref{subsec:shiftedYangian_expectation})}. In particular, the shift depends on the moduli problem and the choice of  stability condition. 
Representations similar to some of those discussed in this paper were considered recently in the case of ``quiver Yangian"   in \cite{LY}.
In the rest of the Introduction we are going to give more details to some of our results.

\subsection{The Drinfeld double}\label{subsec:intr_double}
Let us  summarize some properties of the Drinfeld double in question. 
Assuming that $\H_X$ has a shuffle description we construct the Drinfeld double $D(\mathcal{H})$ (see \cite{STS, L} for the terminology)  of $\H_X$, 
as well as the Drinfeld double $D(\mathcal{SH})$ of the equivariant spherical subalgebra $\mathcal{SH}:=\mathcal{SH}_X$ of $\H_X$ in \S~\ref{Double COHA}. The Drinfeld double $D(\mathcal{SH})$ generalizes the notion of the Cartan doubled Yangian in \cite{Finetal}. {By construction, $D(\mathcal{SH})$ has the following properties} :
\begin{enumerate}
\item There is a triangular decomposition of $D(\mathcal{SH}) \cong \mathcal{SH}^-\otimes \mathcal{SH}^0\otimes \mathcal{SH}^+$, where $\mathcal{SH}^-, \mathcal{SH}^+$ are isomorphic to $\mathcal{SH}$, and $\mathcal{SH}^0=\C[\psi^{(s)}_{i} \mid i\in I, s\in \Z]$ is the algebra of polynomials of infinitely many variables. 
\item The generators of $\mathcal{SH}^+$ are given by $\{e_{i}^{(r)}\mid  i\in I, s\in \N]$. The relations among the generators can be computed by the shuffle description of $\mathcal{SH}$. 
\item
Similarly, the generators of $\mathcal{SH}^-$ are given by $\{f_{i}^{(r)}\mid  i\in I, s\in \N]$. The relations among the generators can be computed by the shuffle description of $\mathcal{SH}$. 
\item
The action of $\mathcal{SH}^0$ on $\mathcal{SH}^+$ is determined by formula \eqref{eq:action}. This gives the relations among $\{\psi^{(s)}_{i}\}$ and $\{e_{j}^{(r)}\}$. 
In a similar way we obtain the relations between $\{\psi^{(s)}_{i}\}$ and $\{f_{j}^{(r)}\}$. 
\item
The commutation relations between $\{e_{i}^{(r_1)}\}$ and $\{f_{j}^{(r_2)}\}$ are given in  Proposition \ref{prop:computeE_iF_j}, which are determined by the definition of Drinfeld double \ref{eq:Drinfeld double mul} (by taking $x=b=1$ and $a=e_{i}^{(r_1)}, y= f_{j}^{(r_2)}$ in  \ref{eq:Drinfeld double mul}). 
\end{enumerate} {
For a coweight $\mu$, we define
 the shifted affine Yangian $Y_{\mu}$  as the quotient of the Drinfeld double $D(\mathcal{SH})$ by certain relations determined by $\mu$ among the generators of $\mathcal{SH}^0$ (see Definition~\ref{def:shiftedYangian}).} 

\subsection{Moduli space of perverse coherent systems}
\label{subsec:modulipervcoh}

{
In  this paper, we consider the moduli space of perverse coherent systems depending on an algebraic cycle $\chi$ of the form \[
N_0 X+\sum_{i=1}^r N_iD_i,\] where $N_i\in\mathbf N$  and each $D_i$ is a toric divisor, and a choice of the stability parameter $\zeta$, as to be explained in \S\ref{sec:D4 branes} and \S\ref{sec:D4 branes2}. A special case is the {\it rank-1 perverse coherent system},which is  a pair consisting  of a  perverse coherent sheaf supported on the exceptional fiber of $X\to Y$ and a homomorphism to it from the structure sheaf of $X$, studied earlier in \cite{Nagao,NN}.  In particular, in \S~\ref{sec:D4 branes2} we will discuss  how the  general $\zeta$ and $\chi$ determine a co-weight $\mu$  using the intersection pairing of $X$ \eqref{eq:pairing}.}
{In this framework the class of representations which we propose to study is described in the following conjecture, see more in \S~\ref{sec:cartanDoubled}:}

\begin{conj}\label{conj:introduction}
There is an action of $D(\mathcal{SH}_X)$ on the cohomology of the moduli space of  perverse coherent system in agreement with the general philosophy of \cite{So}. Such action is obtained from an action of the shifted Yangian $Y_{\mu}$ on the cohomology of the above moduli space, which gives a highest weight module of $Y_{\mu}$. \end{conj}


{
In the special case of rank-1 perverse coherent systems \cite{Nagao,NN}, the algebraic cycle $\chi$ is $X$ itself, and the shifted Yangian in Conjecture \ref{conj:introduction} involves shifting of the imaginary root. 
For the moduli spaces associated of perverse coherent systems supported on a toric divisor in $X$, the algebraic cycle $\chi$ is this toric divisor, and the shifted Yangian from Conjecture \ref{conj:introduction} involves shifting of the real roots. }

{In full generality the Conjecture is unaccessible at this time}.
A reasonable but still difficult class of examples should be related to the case  $X=X_{m, n}\to Y=Y_{m,n}$ (see \cite[\S~1]{Nagao}). By the loc. cit. there is a tilting vector bundle $\calP$ on $X_{m, n}$ such that $\calA_0=\End_X(\calP)$ is isomorphic to the Jacobian algebra of the explicitly written quiver with potential $(Q,W)=(Q_{m,n}, W_{m,n})$.
The generating function of (non-equivariant) cohomology of the moduli space of rank-1 perverse coherent systems on $X_{m, n}$
 was studied in \cite{Nagao} (see also \cite{NN} for $X=X_{1, 1}$).   {A toric diagram of one possible resolution of $Y_{m, n}$ is shown in figure \eqref{fig2} below.} 
 
In this paper we study in detail special cases with $Y=Y_{1,0}$, $Y=Y_{1, 1}$, $Y=Y_{2, 0}$. As we will see later they indeed  deserve special attention. {In particular  for the corresponding resolved conifold $X_{1, 1}$ we confirm the conjecture of Kevin Costello that the deformed algebra of $U(\text{Diff}(\C))$  acts on the cohomology of the moduli space of certain quasi-maps from $\bbP^1$ to $\Hilb(\C^2)$, see \S \ref{Costello Conj} for details. Furthermore, we expect that  the case of general $X_{m, n}$ can be reduced  to the one of these three cases.  }  In the end of this subsection we describe  each of our three cases  in terms of the quiver with potential. As we explained above this description gives rise to the corresponding COHA.

If  $Y=Y_{1, 0}=\C^3$ then, by definition, the resolution $X$ is $\C^3$.  In the two other cases $X$ is the total space of a rank-2 vector bundle on $\bbP^1$. Note that the condition that $X$ is a resolution of singularities of its affinization  implies that $X$ is one of the following types:

a)  the resolved conifold  $X_{1, 1}=\Tot(\calO_{\PP^1}(-1)\oplus\calO_{\PP^1}(-1)) \to Y_{1,1}$;

b) resolution of singularities  of the orbifold $Y_{2, 0}$,  $X_{2, 0}=\Tot(\calO_{\PP^1}(-2)\oplus\calO_{\PP^1}) \to Y_{2,0}=\C^2/(\Z/2)\times \C$. 

The quiver with potential for $X=\C^3$ is 
\begin{equation}
\begin{tikzpicture}\label{eq:Q_3}
  \draw ($(0,0)$) circle (.08);
   \draw[->,shorten <=7pt, shorten >=7pt] ($(0,0)$)
   .. controls +(40:1.5) and +(-40:1.5) .. ($(0,0)$);
     \draw[->,shorten <=7pt, shorten >=7pt] ($(0,0)$)
   .. controls +(90+40:1.5) and +(90-40:1.5) .. ($(0,0)$);
 \draw[->,shorten <=7pt, shorten >=7pt] ($(0,0)$)
   .. controls +(180+40:1.5) and +(180-40:1.5) .. ($(0,0)$);
   \node at (0, 1.2) {$B_2$};
     \node at (-1.2, 0) {$B_1$};
     \node at (1.2, 0) {$B_3$};
\node at (0, -1) {$X=\C^3$};   
\node at (2.5, -1) {$W=B_3[B_1,B_2]$};
\end{tikzpicture} 
\end{equation}
The quivers with potential for the other two cases are:

\begin{equation}
\begin{tikzpicture}[scale=0.8]\label{eq:X11}
    	\node at (-2.2, 0) {$\bullet$};   	\node at (2.2, 0) {$\bullet$};  
		\draw[-latex,  bend left=30, thick] (-2, 0.1) to node[above]{$a_2$} (2, 0.1);
		\draw[-latex,  bend left=50, thick] (-2, 0.3) to node[above]{$a_1$} (2, 0.3);
	\draw[-latex,  bend left=30, thick] (2, -0.1) to node[above]{$b_1$} (-2, -0.1);
		\draw[-latex,  bend left=50, thick] (2, -0.3) to node[above]{$b_2$} (-2, -0.3);
\node at (0, -2) {$X=\Tot(\calO(-1)\oplus\calO(-1))$};   
\node at (-0.3, -2.5) {$W=a_1b_1a_2b_2-a_1b_2a_2b_1$};		
    \end{tikzpicture}\,\ \,\ \text{} \,\ \,\ 
\begin{tikzpicture}[scale=0.8]
            	\node at (-2.2, 0) {$\bullet$};   	\node at (-3.5, 0) {$B_3$};   	
	\node at (2.2, 0) {$\bullet$};  \node at (3.5, 0) {$\tilde{B}_3$};  
		\draw[-latex,  bend left=30, thick] (-2, 0.1) to node[above]{$B_2$} (2, 0.1);
		\draw[-latex,  bend left=50, thick] (-2, 0.3) to node[above]{$B_1$} (2, 0.3);
	\draw[-latex,  bend left=30, thick] (2, -0.1) to node[above]{$\tilde{B}_1$} (-2, -0.1);
		\draw[-latex,  bend left=50, thick] (2, -0.3) to node[above]{$\tilde{B}_2$} (-2, -0.3);
				\path (2.2, 0) edge [loop right, min distance=2cm, thick, bend right=40 ] node {} (2.2, 0);
				\path (-2.2, 0) edge [loop left, min distance=2cm, thick, bend left=40 ] node {} (-2.2, 0);
\node at (0, -2) {$X=\Tot(\calO(-2)\oplus\calO)$};   
\node at (0, -2.5) {$W=B_2 \tilde{B}_1\tilde{B}_3-\tilde{B}_2\tilde{B}_3B_1+\tilde{B}_2B_1 B_3-B_2B_3\tilde{B}_1$};					
    \end{tikzpicture}
\end{equation}

\subsection{Imaginary roots and shifted Yangians} 
 
 In this subsection we are going to describe the representation spaces for the doubles of our COHAs.
 
The moduli space of rank-1 perverse coherent systems depends on a choice of stability condition.  We refer the readers to \cite{Nagao} for the definition of stability condition, which we briefly recall in \S~\ref{sec:PTconifold} for  convenience of the reader. 
We understand the  space of  stability conditions in question as a real vector space which contains a hyperplane arrangement. The latter  divides the vector space into the disjoint union of open connected components called chambers.  A stability condition $\zeta$ is said to be {\it generic} if it belongs to the complement of the union of hyperplanes, i.e.  to  the interior of a chamber.  
The moduli space $\fM_\zeta$ of $\zeta$-stable perverse coherent systems  depends only on the chamber that contains $\zeta$ but not on the particular $\zeta$. For $X=X_{m, n}$ the hyperplane arrangement coincides with the one given by root hyperplanes of the affine Lie algebra of type $A_{m+n-1}$. The chamber structure of the space of stability conditions was  described in \cite{Nagao}.  Depending on the chamber the moduli space of stable perverse coherent systems can be identified with the moduli spaces which appear 
in Donaldson-Thomas (DT), Pandharipande-Thomas (PT) or non-commutative Donaldson-Thomas (NCDT) theories. We will sometimes refer to them as DT, PT or NCDT moduli spaces depending on the chamber under consideration.

If $X=X_{1, 1}$ or $X=X_{2, 0}$, a stability condition is given by a pair of real numbers $\zeta=(\zeta_0,\zeta_1)$. 
Walls of the chambers are labeled by non-negative integers, see \cite[Figure 1]{NN}:
\begin{align}
&m\zeta_0+(m+1) \zeta_1=0, m\in \Z_{\geq 0}, && \text{the DT side,}\\
&\zeta_0+\zeta_1=0,  && \text{the imaginary root hyperplane,}\\
&(m+1)\zeta_0+m \zeta_1=0, m\in \Z_{\geq 0}, &&  \text{the PT side}. 
\end{align}
The DT moduli space  corresponds to the chamber which is positioned immediately below the imaginary root hyperplane and above $m\zeta_0+(m+1) \zeta_1=0$, for 
$m \gg 0$.
For the PT moduli space we consider the chamber which is immediately above the imaginary root hyperplane and below $(m+1)\zeta_0+m \zeta_1=0$, for $m \gg 0$. 

One can describe $\fM_\zeta$ in terms of a quiver with potential. For $X=\C^3$ this is straightforward (see  \S~\ref{sec:Hilb} for the details). The case $X=X_{1,1}$ is discussed in  \S~\ref{subsec:COHAaction}. 

Assume in general that  $(Q,W)$ is the quiver with potential such that $\calA_0=\End_X(\calP)$ is isomorphic to the Jacobian algebra of $(Q,W)$ \cite[\S~1]{Nagao}. 
From $(Q,W)$ one constructs a {\it framed quiver} $\tilde{Q}$ with the {\it framed potential} $\tilde{W}$. The framed quiver  is obtained from  $Q$ by adding an additional {\it framing vertex} as well as an arrow from the framing vertex to one of the vertices of $Q$. We define $\tilde{W}=W$.
The case of $X_{1,1}$ is depicted below.
\begin{equation}\label{eq: ext Q}
\begin{tikzpicture}[scale=0.8]
         \node at (-2.2, 0) {$\bullet$};   	
         \node at (-2.2 -0.5, 0) {$V_0$}; 
          \node at (2.2+0.5, 0) {$V_1$};   	
	\node at (2.2, 0) {$\bullet$};  
		\draw[-latex,  bend left=30, thick] (-2, 0.1) to node[above]{$a_2$} (2, 0.1);
		\draw[-latex,  bend left=50, thick] (-2, 0.3) to node[above]{$a_1$} (2, 0.3);
	\draw[-latex,  bend left=30, thick] (2, -0.1) to node[above]{$b_1$} (-2, -0.1);
		\draw[-latex,  bend left=50, thick] (2, -0.3) to node[above]{$b_2$} (-2, -0.3);
\node at (-2.2, -2) {$\square$};  
 \draw[->, thick] (-2.2, -1.8) -- (-2.2, -0.2) ;
 \node at (-2.5, -1) {$\iota$}; 
  	\node at (-2.7, -2) {$1$}; 
	\node at (0, -2.5){Extended quiver $\tilde{Q}$ for $X_{1, 1}=\Tot(\calO(-1)\oplus\calO(-1))$};
    \end{tikzpicture}
\end{equation}  
The stability condition $\zeta$ defines the corresponding stability  condition $\tilde{\zeta}$  for the abelian category of representations of $(\tilde{Q}, \tilde{W})$ (i.e. the abelian category of representations of the corresponding Jacobi algebra). It is given in terms of the slope function $\theta_{\tilde{{\zeta}}}$ (see \S~\ref{sec:PTconifold}). The moduli space $\fM_\zeta$ is isomorphic to moduli space of $\theta_{\tilde{\zeta}}$-stable representations of $(\tilde{Q}, \tilde{W})$ with 1-dimensional framing. 

Notice that $X$ admits an action of the torus $T\simeq (\bf{C}^\ast)^2$ which preserves the canonical class. 
Let  $\calP$ be the  tilting vector bundle on $X$. The tilting object $\calP$ is equivariant with respect to this $T$-action, and hence induces a $T$-action on  $\fM_\zeta$.  Let $\C^*_{fr}$ be the 1-dimensional torus action on $\fM_\zeta$ by scaling the framing. Then $\fM_\zeta$ carries an action of $T\times \C^*_{fr}$.


We state the following result for  $X=X_{1,1}$, which is expected to hold for an arbitrary $X$ (see \cite[Proposition 4.1.1]{So}). 
\begin{prp}\label{prp:increase}
Let $\zeta$ be a generic stability condition. Then there is an action of $\H_X$ on $H^*_{c,T\times \C^*_{fr}}(\fM_\zeta,\varphi_{\tr W}\C)^\vee$. 
\end{prp}
The $\H_X$-action in Proposition \ref{prp:increase} is given by the ``raising operators". We expect the action of the spherical subalgebra $\mathcal{SH}_X\subset \H_X$ in Proposition~\ref{prp:increase} extends to its Drinfeld double $D(\mathcal{SH}_X)$  by adding the action of ``lowering operators".
Furthermore, the action of the Drinfeld double factors through the shifted Yangian  with shifts of the imaginary roots.

Let  $\vec{z}=(z_1, z_2, \cdots, z_{|l|}) (l\in \Z)$ be complex parameters. 
Definition of the  shifted Yangian  $Y_{l}(\vec{z})$ of $\widehat{\mathfrak{gl}(1)}$ as an associative algebra depending  on the parameters $\vec{z}$ is recalled in \S~\ref{sec:shiftedYangian}. 
We refer to it as  the Yangian which is $l$-positively ($-l$-negatively) shifted in parameters $\vec{z}$ if $l>0$ (if $l<0$). When $l=0$ and there are no parameters, and we obtain the usual affine Yangian $\affY$. Here the parameters $\hbar_i, 1\le i
\le 3$ correspond to the weights of the Calabi-Yau action of the torus $(\C^{\ast})^3$, hence $\hbar_1+\hbar_2+\hbar_3=0$. We hope the reader will not confuse the subscripts in the cases of shifted and usual Yangians. 
Let $R_T$ be the cohomology ring $H_{c, T}^\ast(\pt)^\vee$, and $K_T$ its field of fractions. Let $\Hilb^*(\C^3):= \sqcup_{n\in\bbN}\Hilb^n(\C^3)$ be the Hilbert scheme of points in $\C^3$ (see \S\ref{sec:Hilb} for details).

\begin{thm}\label{thm:intro}
\begin{enumerate}
\item Let $X=\C^3$. The action of $\mathcal{SH}_{\C^3}$ on $H^*_{c, T\times \C^*_{fr}}(\Hilb^*(\C^3), \varphi_{\tr W}\C)^\vee$ extends to an action of $Y_{-1}(z_1)$ on $H^*_{c,T\times \C^*_{fr}}(\Hilb^*(\C^3), \varphi_{\tr W}\C)^\vee \otimes_{R_T} K_T$.
\item Let $X=X_{1,1}$ be the resolved conifold. For a generic  stability condition on the PT-side of the imaginary root hyperplane, 
the action of $\H_{X_{1, 1}}$ on $H^*_{c,T\times \C^*_{fr}}(\fM_\zeta, \varphi_{\tr W}\C)^\vee$ induces an action of $Y_{1}(z_1)$ on $H^*_{c,T\times \C^*_{fr}}(\fM_\zeta, \varphi_{\tr W}\C)^\vee \otimes_{R_T} K_T$. 
\end{enumerate}
 \end{thm} 
Theorem \ref{thm:intro} (1) is proved in  \S \ref{sub:Action of the shifted Yangian} and  Theorem \ref{thm:intro} (2) is proved in \S\ref{sec:PTconifold} (see the proof of Theorem \ref{thm:conifoldYangian} there). 

 \begin{rmk}
a) The algebra $Y_{1}(z_1)$ can also be realized as the deformed double current algebra for $\mathfrak{gl}(1)$. The fact that this algebra acts on the PT-moduli space of the resolved conifold was originally conjectured by Costello from holography considerations in M-theory \cite[\S~14 Conjecture]{Costello}.

 b) We expect the same calculation goes through for an arbitrary  $X_{m, n}$. We expect that the action of $\mathcal{H}_{X_{m,n}}$  extends to the action of the shifted versions of affine Yangian\footnote{See also closely related discussion in \cite{FG,GLP,GR,KS,K,LL} and references therein.} of $\mathfrak{gl}(m|n)$ \cite{PR,R}. The above discussion suggests that the choice of stability condition leads to the shift of imaginary roots in the case when we consider the action on the cohomology of the moduli spaces of torsion free sheaves on $X_{m,n}$.  
 \end{rmk} 
 
\subsection{Relations of various algebras}
There are various algebras in the paper. For the convenience of the readers, we summerize the relations as follows. 

Let $X$ be a CY 3 fold. Denote by $\mathfrak{M}_{\zeta, \chi}(X)$ the moduli space of perverse coherent systems depending on an algebraic cycle $\chi$ and the stability condition $\zeta$. Assume  $\mathfrak{M}_{\zeta, \chi}(X)$ carries an action of $T\times G_{fr}$. Then we have the following diagram of morphisms (see Conjecture \ref{conj:introduction}): 
\[
\xymatrix@R=1em @C=1em{
\mathcal{SH}_X \ar@{^{(}->}[d]\, \ar@{}[r]|-*[@]{\subset} &
\mathcal{H}_X \ar@{^{(}->}[d]   \ar[rdd]  &&\\
D(\mathcal{SH}_X) \, \ar@{}[r]|-*[@]{\subset} \ar@{->>}[d]& D(\mathcal{H}_X) &
& 
\\
Y_{\mu} \ar@{-->}[rr]&& \End( H^*_{c, T\times G_{fr}}(\mathfrak{M}_{\zeta}(X),  \varphi)^{\vee}\otimes_{R_T} K_T)
}
\]

In the special case when $X=\C^3$, and $\mathfrak{M}_{\zeta}(X)=\Hilb^*(\C^3)$, the algebra $\mathcal{SH}_X $ is identified with $ \affY^+$, and the quotient of
$D(\mathcal{SH}_{\C^3})$, $Y_\mu$, with $Y_{-1}(z_1)$.
For a general $X$, and general $\mathfrak{M}_{\zeta}(X)$, the shifted algebra $Y_{l}(\vec{z})$ of $\widehat{\mathfrak{gl}_1}$ is a subalgebra of $Y_{\mu}$. 
For example, when $X=X_{1, 1}=\calO(-1)\oplus \calO(-1)$, and $\mathfrak{M}_{\zeta}(X)$ as in Theorem~\ref{thm:intro}(2), we have 
 \[
\xymatrix@R=1em @C=1em{
 \affY^+\cong  \mathcal{SH}_{\C^3}  \ar@{^{(}->}[d] & \mathcal{SH}_{X_{1, 1}} \ar[l] \ar@{^{(}->}[d]
 \\
D(\mathcal{SH}_{\C^3}) \ar@{->>}[dd]&D(\mathcal{SH}_{X_{1, 1}}) \ar[l] \ar@{-->>}[dd]
\\
& &&
\\
Y_{1}(z_1)\ar@{^{(}-->}[r]&Y_{\mu}
\\
}
\]
Here $Y_{\mu}$ is expected to be related to $Y_{1}(\mathfrak{gl}(1|1))$.
%
In general, we describe subalgebras of $Y_{\mu}$ corresponding to each "root", and we expect all such subalgebras determine $Y_{\mu}$ uniquely.

\subsection{Real roots, toric divisors  and shifted Yangians}

We are going to illustrate the idea in a simple example.
Consider the Calabi Yau 3-fold $X_{2, 0}=T^*\PP^1\times \C$ and the effective divisor $D$ to be fiber of 
$T^*\PP^1\times \C \to \PP^1$ over the north pole or over the south pole of $\PP^1$. Depending on the choice of the toric divisor the corresponding 
geometry can be depicted  as in the following toric diagram
\begin{equation}\label{Intro:N or S pole}
\begin{tikzpicture}[scale=0.9]
    \draw[-, thick] (0,1) -- (2, 1) ;
    \draw[-, thick] (0,0) -- (2, 0);
    \draw[-, thick] (0,0) -- (0, 1);
    \draw[-, thick] (0, 1) -- (-1, 2);
    \draw[-, thick] (0, 0) -- (-1, -1);
         \filldraw (0,1) node[anchor=north, yshift=0.8cm, xshift=0.3cm] {$1$};
    \end{tikzpicture} \,\ \,\  \,\ \,\   \,\ \,\  \,\ \,\   \,\ \,\  \,\ \,\   \,\ \,\  \,\ \,\  
    \begin{tikzpicture}[scale=0.9]
    \draw[-, thick] (0,1) -- (2, 1) ;
    \draw[-, thick] (0,0) -- (2, 0);
    \draw[-, thick] (0,0) -- (0, 1);
    \draw[-, thick] (0, 1) -- (-1, 2);
    \draw[-, thick] (0, 0) -- (-1, -1);
         \filldraw (0,0) node[anchor=south, yshift=-0.8cm, xshift=0.3cm] {$1$};
    \end{tikzpicture}
    \end{equation}
In Section \ref{Chainsaw} we obtain a quiver $Q$ with potential $W$ associated with the pair consisting of the Calabi Yau 3-fold $T^*\PP^1\times \C$ and one of
the above toric divisors $D$. The construction  uses a $\Z/2\Z$ symmetry of the quiver with potential 
associated with the pair consisting of the Calabi Yau 3-fold $\C^3$ and the effective divisor which is the coordinate plane.

By the dimension reduction, $(Q, W)$ can be described as the following ``chainsaw quiver" 
\begin{equation}\label{eqn:chainsaw}
\begin{tikzpicture}[scale=0.8]
        	\node at (-2.2, 0) {$\bullet$};   	
	\node at (-2.2, 0.5){$V_0$}; 
	\node at (-3.5, 0){$B_3$};
	\node at (2.2, 0) {$\bullet$};  
	\node at (2.2, 0.5){$V_1$};
	\node at (3.5, 0){$\tilde{B}_3$};
	\draw[-latex, thick] (-2, 0.3) to node[above]{$B_1$} (2, 0.3);
	\draw[-latex, thick] (2, -0.1) to node[below]{$\tilde{B}_1$} (-2, -0.1);
				\path (2.2, 0) edge [loop right, min distance=2cm, thick, bend right=40 ] node {} (2.2, 0);
				\path (-2.2, 0) edge [loop left, min distance=2cm, thick, bend left=40 ] node {} (-2.2, 0);
\node at (0, -3) {$\triangle$};  \node at (0.5, -3) {$1$};  
 \node at (-1.5, -2) {$I_{13}$};   \node at (1.5, -2) {$J_{13}$};  
 \draw[->, thick] (0, -2.8) -- (-2.2, -0.2) ;
 \draw[->, thick]  (2.2, -0.2) --(0, -2.8) ;
    \end{tikzpicture}
\end{equation}   
with relations
\begin{equation*}
\tilde{B}_1\tilde{B}_3-B_3\tilde{B}_1=I_{13} J_{13}, \,\ 
\tilde{B}_3B_1=B_1 B_3
\end{equation*}
Analogously to the K-theory case in \cite{N}, the $(1, 0)$-shifted Yangian of  $\widehat{\mathfrak{gl}(2)}$ (see Definition in \S\ref{Chainsaw}) acts on the cohomology of the chainsaw quiver variety. 
Hence  it also acts on 
\[
\bigoplus _{V_0, V_1} H^*_{c, \GL(V_0)\times  \GL(V_1)\times T}(\Rep(Q, V_0, V_1)^{st}, \varphi_{tr W})^\vee. 
\]
Here the superscript $st$ refers to the locus of stable representations.

In  physics language the above example describes a  stack of $D4$-branes wrapping a union of smooth toric divisors \cite{AJS,GO,NSY,NP}. 
The corresponding moduli space is the moduli space of pure sheaves supported on the divisors, with the generic rank  at each component  specified by the number of the $D4$-branes. We expect that the Yangian shifted by real simple roots  acts on the cohomology of the corresponding moduli space. The shift is determined by the intersection number of the $\bbP^1$ corresponding to the real simple root and the divisor, counted with multiplicity given by the number of $D4$-branes.

\subsection*{Contents of the paper}
In \S~\ref{Double COHA} following the original approach of Kontsevich and Soibelman we recall the definition of COHA  associated to a smooth toric Calabi-Yau 3-fold. 
When the COHA has a shuffle description, we construct the Drinfeld double of the equivariant spherical COHA. 
In the Drinfeld double of the spherical COHA, one has a commutative subalgebra which is isomorphic to the polynomial algebra with generators labelled by simple roots and $\Z$.  We refer this subalgebra as the Cartan doubled subalgebra. 
We define the shifted Yangian 
as the quotient of the Drinfeld double by certain relations  among the  generators of the Cartan doubled subalgebra, which are determined by a coweight. For a specific choice of the quiver with potential, we show in  \S~\ref{Examples COHA} that the Cartan doubled Yangian of a Kac-Moody Lie algebra defined in \cite{Finetal} surjects to the Drinfeld double of the equivariant spherical COHA of this quiver with potential. We also recall the description of COHA of $\C^3$ from our previous paper. Furthermore, we prove that the Drinfeld double corresponding to the resolved conifold has a quotient algebra which is isomorphic to the affine Yangian of  $\widehat{\mathfrak{gl}}(1)$. In  \S~\ref{sec:shiftedYangian} we discuss the shifted affine Yangians of $\widehat{\mathfrak{gl}}(1)$. We also discuss the relationship of COHA to quantized Coulomb branch algebras of the Jordan quiver gauge theory. We connect this discussion with a special case of the conjecture of Costello. 

In \S~\ref{sec:Hilb} we construct  geometrically an action of the 1-negatively shifted affine Yangian on the cohomology of the Hilbert scheme of $\C^3$. In \S~\ref{sec:PTconifold} we construct  geometrically an action of the one-positively shifted affine Yangian on cohomology of the PT moduli space of the resolved conifold. In \S~\ref{sec:D4 branes} we construct geometrically an action of shifted affine Yangian on the cohomology of the moduli of perverse coherent systems supported a toric divisor in the resolved $\Z_2$-orbifold.

In \S~\ref{sec:D4 branes2}, we make a proposal in the general case. Namely, we expect  the action of the Drinfeld double of the spherical COHA on more general moduli spaces of perverse coherent sheaves on the general toric Calabi-Yau 3-fold $X_{m,n}$. The action of the  Drinfeld double is expected to factor through the shifted Yangian, with the shift determined by certain intersection numbers. 
We verify this proposal in  the examples above. 

In \S~\ref{sec:appendix}, we prove a Jeffery-Kirwan-type residue formula for pushforward in critical cohomology. 

\subsection*{Acknowledgments}

We thank Kevin Costello, Ben Davison, Hiraku Nakajima, Nikita  Nekrasov, Masahito Yamazaki for very helpful discussions. 

Part of the work was done when the last two named authors were visiting at the IPMU. 
M.R. was supported by NSF grant 1521446 and 1820912, the Berkeley
Center for Theoretical Physics and the Simons Foundation. Y.S. was partially supported by the Munson-Simu Star Faculty award of KSU.  Part of the work was done when he was supported by the Perimeter Institute for Theoretical physics as a Distinguished Visiting Research Chair.
Y.Y. was partially supported by 
the Australian Research Council (ARC) via the award DE190101231. G.Z. was partially supported by ARC via the award
DE190101222. 

The work was started when all  authors were visiting the Perimeter Institute for Theoretical Physics (PI). 
They are grateful to PI for excellent research conditions.
The Perimeter Institute for Theoretical Physics is supported by
the Government of Canada through the Department of Innovation, Science and Economic
Development and by the Province of Ontario through the Ministry of Research, \& Innovation
and Science.

We are  grateful to the referees for the careful reading the manuscript and for the detailed comments and suggestions which helped us to improve the manuscript. 

\section{The Drinfeld double of COHA associated to a quiver with potential}
\label{Double COHA}

\subsection{The Hopf algebra structure} 
\label{H=Hopf}
Let $R$ be a commutative integral domain of characteristic zero. Let $H$ be a Hopf algebra over $R$. 
By definition, $H$ posses the following $R$-linear operations
\begin{enumerate}
\item multiplication $\mu: H\otimes H\to H$,
\item unit $\eta: R\to H$,
\item comultiplication $\Delta: H \to H\otimes H$, 
\item co-unit $\epsilon: H\to R$
\item antipode $S: H\to H$
\end{enumerate}
which satisfy the following conditions
\begin{enumerate}
\item[(a)] $(H, \mu, \eta)$ is an associative $R$-algebra. 
\item[(b)] $(H, \Delta, \epsilon)$ is an associative $R$-coalgebra. 
\item [(c)] $\Delta, \epsilon$ are algebra homomorphisms. 
\item [(d)] $\mu(S\otimes 1)\Delta=\mu( 1\otimes S)\Delta=\eta\epsilon$. 
\end{enumerate}

We denote by $\H$ the cohomological Hall algebra (COHA for short) associated with  a quiver with potential $(Q,W)$ defined by  Kontsevich and Soibelman \cite{KS}. 
We briefly recall the definition of $\H$. We write $Q=(I,H)$ with $I$ being the set of vertices and $H$ the set of arrows. 
For a dimension vector $v=(v^i)_{i\in I}\in\bbN^I={\bf Z}_{\ge 0}^I$, let $\Rep(Q,v)$ be the algebraic variety (in fact with non-canonical structure of a vector space) parameterizing representations of $Q$ on the $I$-graded complex vector space whose degree $i\in I$-piece is $\C^{v_i}$. On $\Rep(Q,v)$ there is an action of $\GL_v$ by changing  basis. We assume there is an action of $T=(\bf{C}^\ast)^2$ on $\Rep(Q,v)$ that commutes with the action of $\GL_v$ and preserves the potential $W$. Below the $T$-action  is  induced from its Calabi-Yau action on $X$ and hence on $\End_{Coh(X)}(\calP)$). The trace of $W$ gives rise to the regular function $\tr W$ on $\Rep(Q,v)$ which is invariant under the $(\GL_v\times T)$-action. 
In this notation  the equivariant COHA of $(Q,W)$, is a $\bbN^I$-graded vector space 
\[
\H=\H^{(Q,W)}:=\oplus_{v\in\bbN^I} \H_v\] with $\H_v=H^*_{c,\GL_v\times T}(\Rep(Q,v),\varphi_{\tr W})^\vee,$ endowed with the Hall multiplication $m^{\crit}$, 

Set $v=v_1+v_2$. Define $\Rep(Q)_{v_1, v_2}:=\{x\in \Rep(Q,v)\mid x(V_1)\subset V_1\}$.  We write $G$ instead of $G_{v}$ for short when $v$ is understood from the context. Let $P\subset G_v$ be the parabolic subgroup preserving the subspace $V_1$ and $L:=G_{v_1}\times G_{v_2}$ be the Levi subgroup of  $P$. We have the following correspondence of $L$-varieties. 
\begin{equation}\label{basic corresp}
\xymatrix
{\Rep(Q, v_1)\times \Rep(Q, v_2)&\Rep(Q)_{v_1, v_2} \ar[l]_(0.4){p}\ar[r]^(0.4){\eta} &\Rep(Q, v),
}\end{equation}
where $p$ is the natural projection and $\eta$ is the embedding. 
The trace function $\tr W_{v_i}$ on $\Rep(Q, v_i)$ induces a function $\tr W_{v_1}\boxplus \tr W_{v_2}$
on the product $\Rep(Q, v_1)\times \Rep(Q, v_2)$.
We define $\tr(W)_{v_1, v_2}$ on $\Rep(Q)_{v_1, v_2}$ to be
\[
\tr(W)_{v_1, v_2}:=p^*(\tr W_{v_1}\boxplus \tr W_{v_2})=\eta^*(\tr W_{v_1+v_2}).
\]
Note that we have
$
p^{-1} (\Crit(\tr W_{v_1})\times \Crit(\tr W_{v_2}))
\supsetneqq
\eta^{-1}(\Crit(\tr W_{v_1+v_2})).
$

For simplicity, we assume that $\Crit(\tr W_v)\subseteq (\tr W_v)^{-1}(0)$ so that we have the Thom-Sebastiani isomorphism as below. 

The Hall multiplication $m^{\crit}$,  defined in \cite[\S~7.6]{KS}, is the composition of the following morphisms:
\begin{enumerate}
\item The Thom-Sebastiani isomorphism
\begin{align*}
&H_{c, G_{v_1}\times T}^*(\Rep(Q,v_1), \varphi_{\tr W_{v_1}})^\vee \otimes 
H_{c, G_{v_2}\times T}^*(\Rep(Q,v_2), \varphi_{\tr W_{v_2}})^\vee
\\& \cong 
H_{c, L\times T}^*(\Rep(Q,v_1)\times\Rep(Q,v_2\times T), \varphi_{\tr W_{v_1} \boxplus \tr W_{v_2}} )^\vee.
\end{align*}
\item 
Using the fact that $\Rep(Q)_{v_1, v_2}$
is an affine bundle over $\Rep(Q, v_1)\times \Rep(Q, v_2)$, and
$\tr W_{v_1, v_2}$ is the pullback 
of $\tr W_{v_1} \boxplus \tr W_{v_2}$, we have
\begin{align*}
p^*:H_{c, L\times T}^*(\Rep(Q, v_2)\times \Rep(Q, v_2),
 \varphi_{\tr W_{v_1} \boxplus \tr W_{v_2}})^\vee
\cong
&H_{c, L\times T}^*(\Rep(Q)_{v_1, v_2}, 
\varphi_{\tr W_{v_1, v_2} } )^{\vee}
\end{align*}
\item
Using the fact $\tr W_{v_1, v_2}$ is the restriction of $\tr W_{v}$ to
$\Rep(Q)_{v_1, v_2}$. We have
\begin{align*}
\eta_*:H_{c, L\times T}^*(\Rep(Q)_{v_1, v_2}, 
\varphi_{\tr W_{v_1, v_2} } )^{\vee}\to &
H_{c, L\times T}^*(\Rep(Q, v), 
\varphi_{\tr W_{v} } )^{\vee}.
\end{align*}
\item 
Pushforward along 
$G\times_{P} \Rep(Q, v)\to \Rep(Q, v), (g, m)\mapsto gmg^{-1}$,
 we get
\begin{align*}
H_{c,L\times T}&(\Rep(Q, v), 
\varphi_{\tr W_{v} } )^{\vee}\cong H_{c,P\times T}(\Rep(Q, v), 
\varphi_{\tr W_{v} } )^{\vee}\\&\cong H_{c, G\times T}^*(G\times_{P} \Rep(Q, v), 
\varphi_{\tr W_{v} } )^{\vee} \to 
H_{c, G\times T}^*(\Rep(Q, v), 
\varphi_{\tr W_{v} } )^{\vee}.
\end{align*}
\end{enumerate}

We also denote $\sqcup_{v\in\bbN^I}\Rep(Q,v)$ by $\Rep(Q)$.
Let $R:= \C[\Lie(T)]$ be the ring of functions on $\Lie(T)$. 
Then $\H$ is an $R$-algebra.

We now construct the Drinfeld double of $\H$. To do this, we assume $\H$ has a shuffle description. Explicitly, for $v \in \N^I$, let $\g_v:=\prod_{i\in I} \mathfrak{gl}_{v^i}$ be the Lie algebra, 
and $\h_v\subset \g_v$ be its Cartan subalgebra. 
The Weyl group $S_v$ naturally acts on the ring of functions $\C[\h_{v}]$.
For each $v\in \N^I$, $\C[\h_{v}]$ is a polynomial ring with variables $\{x_{1}^{(i)},  x_{2}^{(i)}, \cdots, x_{v^i}^{(i)} \mid i\in I\}$. 
For simplicity, we will denote this set of variables by $x_{[1, v]}$. 

Let $\textbf{H}$ be the shuffle algebra defined as follows. As a vector space, we have 
$$\textbf{H}=\oplus_{v\in \N^I}\textbf{H}_v=  \oplus_{v\in \N^I} R[\h_{v}]^{S_v}.$$ 
Consider the embeddings
\begin{align*}
R[\h_{v_1}]^{S_{v_1}}\subset R[\h_{v_1+v_2}]^{S_{v_1}\times S_{v_2}}, \,\ & x_{[1, v_1]}\mapsto  x_{[1, v_1]}\\
R[\h_{v_2}]^{S_{v_2}}\subset R[\h_{v_1+v_2}]^{S_{v_1}\times S_{v_2}}, \,\ & x_{[1, v_2]}\mapsto  x_{[v_1+1, v_1+v_2]}. 
\end{align*}
For any pair $(p, q)$ of positive integers, let $\Sh(p, q)$ be the subset of $S_{p+q}$ consisting of $(p, q)$-shuffles
(permutations of $\{1, \cdots , p+q \}$ that preserve the relative order of $\{1, \cdots , p\}$ and $\{p + 1, \cdots , p + q\}$).
Set $\Sh(v_1, v_2):=\prod_i \Sh(v_1^{(i)}, v_2^{(i)})$. 
The multiplication of  $\textbf{H}$ given as
\begin{align*}
\star: R[\h_{v_1}]^{S_{v_1}}\otimes R[\h_{v_2}]^{S_{v_2}} \to R[\h_{v_1+v_2}]^{S_{v_1+v_2}}, \,\ 
 f\star g&=\sum_{\sigma\in \Sh(v_1, v_2)} \sigma(f \cdot g \cdot \fac(x_{[1, v_1]}| x_{[v_1+1, v_1+v_2]})), 
\end{align*}
where $\fac(x_{[1, v_1]}| x_{[v_1+1, v_1+v_2]}) \in (R[\h_{v_1}]^{S_{v_1}}\otimes R[\h_{v_2}]^{S_{v_2}})_{\loc}$ is an explicit rational function with denominator
$\prod_{i\in I} \prod_{a=1}^{ v_1^{(i)}} \prod_{b=v_1^{(i)}+1}^{ v_1^{(i)}+v_2^{(i)}}(x^{(i)}_a - x^{(i)}_b)$. Here the subindex $\loc$ stands for localization with respect to the divisor defined by $\fac(x_{A}| x_{B})$.

Moreover, it has the property that
\[
\fac(x_{A_1 \sqcup A_2 }| x_{B})=\fac(x_{A_1}| x_{B})\fac(x_{A_2}| x_{B}), \fac(x_{A }| x_{B_1\sqcup B_2})=\fac(x_{A}| x_{B_1})\fac(x_{A}| x_{B_2}), 
\]
for any subsets $x_{A_1}, x_{A_2}, x_A, x_{B_1}, x_{B_2}, x_B$ of the collection of variables $x_{[1, v]}$.  
In the case that $(Q,W)$ has a cut, the formula of $\fac(x_{A}| x_{B})$ can be found in \cite[Appendix~A]{YZ4}. We do not need the explicit expression in this paper. 

Let $ \mathcal{SH} $ be the equivariant spherical COHA (\cite{SV, RSYZ}). By definition, $\mathcal{SH}$ is a subalgebra of $\H$ generated by $\H_{e_i}$, as  $i$ varies in $I$. We will skip the word ``equivariant" if it does not lead to a confusion. 
Similarly, let $\SH$ be the spherical shuffle algebra, which is generated by $\textbf{H}_{e_i}$, as  $i$ varies in $I$. 
We assume that there is an algebra homomorphism
\[
\H\to \textbf{H}
\]
which is an isomorphism after passing to the localization $-\otimes_{R[\h_{v}]^{S_v}} R(\h_{v})^{S_v}$, for each $v$. 
We have an induced algebra epimorphism 
\[
\mathcal{SH}\to \SH, 
\]
which is an isomorphism after passing to the same localization.

Let $\H^0:=\C[\psi_{i, r}|i\in I, r\in \N]$ be the polynomial ring with infinitely many formal variables. 
Let $
\psi_{i}(z)=1+\sum_{r \geq 0} \psi_{i, r} z^{-r-1}\in\null \H^0[\![ z^{-1}]\!]$
 be the generating series of generators $\psi_{i, r} \in  \H^0$.  
 Similar to \cite{YZ2}, we define the extended shuffle algebra  $\textbf{H}^{\ext}:=\H^0\ltimes \textbf{H}_{\loc}$ using the 
$\H^0$--action on $\textbf{H}_{\loc}$ by 
 \begin{equation}\label{eq:action}
 \psi_{i}(z) g  \psi_{i}(z)^{-1}:=g \frac{\fac(z|x_{[1, v]})}{ \fac(x_{[1, v]}| z)} , \text{for any $g \in \textbf{H}_{v}$, }
  \end{equation}
  where $\loc$ means we invert elements of the form $\fac(x_{[1, v]}| z)$. 
Note that restricting to the spherical subalgebra $\SH$, we obtain that $\H^0$ acts on $\SH$. The  localization is not needed for $\SH$. To say it differently, $ \H^0\ltimes \SH$ is  a subalgebra of  $\H^0\ltimes \textbf{H}_{\loc}$.

Following the  same construction as in \cite{YZ2}, we define  a localized coproduct on $\textbf{H}^{\ext}$ \[
\Delta: \textbf{H}^{\ext}\to \sum_{v_1+v_2=v} (\textbf{H}_{v_1}^{\ext} \otimes \textbf{H}_{v_2}^{\ext})_{\loc}. 
\]
Here  $\loc$ means the localization away from   the union of null-divisors of $\fac(x_{[1,v_1]}|x_{[v_1+1, v_1+v_2]})$ over all $v_1+v_2=v$ and $\psi_i(z)$. In particular, the functions $\psi_i(z), \fac(x_{[1,v_1]}|x_{[v_1+1, v_1+v_2]})$ have inverse in $\textbf{H}^{\ext}_{\loc}$. 

On the symmetric algebra $\H^0$, the map $\Delta$ is given  by
\begin{align}\label{coprodH}
&
\Delta(\psi_i(w))=\psi_i(w)\otimes \psi_i(w), i\in I.
\end{align} 
For an element $f(x^{(i)})\in \H_{e_i}$, we define $\Delta$ on $\H_{e_i}$ by
\[
\Delta(f(x^{(i)}))=\psi_i(x^{(i)})\otimes f(x^{(i)})+f(x^{(i)})\otimes 1 \in \H_{e_i}^{\ext}\otimes \H_{e_i}^{\ext}. 
\]
We extend $\Delta$ on the entire $\textbf{H}$ by requiring $\Delta(a\star b)=\Delta(a)\star \Delta(b)$. 
For a homogeneous element $f(x_{[1, v]})\in \textbf{H}_v$, $\Delta$ is given as 
\begin{equation}\label{eq:coprod}
\Delta(f(x_{[1, v]}))=
\sum_{\{v_1+v_2=v\}}
\frac{\psi_{[1,v_1]}(x_{[v_1+1,v]}) f(x_{[1,v_1]}\otimes x_{[v_1+1,v]})}{\fac( x_{[v_1+1,v]}| x_{[1,v_1]})}, 
\end{equation}
where $\psi_{[1,v_1]}(x_{[v_1+1,v]}) :=\prod_{k\in I} \prod_{ j= v_1^{k}+1, \cdots, v^{k}} 
\psi_k (x_j^{(k)})$. The subindex $[1, v_1]$ indicates that the factor $\psi_k (x_j^{(k)})=1+\sum_r \psi_{k, r}\otimes (x_j^{(k)})^{-r-1}$ lies in $ \textbf{H}_{v_1, \loc}^{\ext}\otimes  \textbf{H}_{e_k, \loc}^{\ext}$. 

Define the co-unit 
\[
\epsilon: \textbf{H}^{\ext}\to R, \psi_i(z) \mapsto 1, f(x_{[1, v]})\mapsto 0. 
\]
for $v\neq 0$, $  \psi_i(z)\in \H^0[[z^{-1}]],  f(x_{[1, v]})\in \textbf{H}_{v}$. 

Define the antipode $S: \textbf{H}^{\ext}_{\loc}\to \textbf{H}^{\ext}_{\loc}$ as follows. 
\begin{align}
&  \psi_i(z) \mapsto \psi_i^{-1}(z) \label{psi}\\
&  f(x_{[1, v]}) \mapsto (-1)^{|v|}  \psi_{[1, v]}^{-1}(x_{[1, v]}) f(x_{[1, v]})  \label{psix}
\end{align}
for $  \psi_i(z)\in \H^0[[z^{-1}]],  f(x_{[1, v]})\in \textbf{H}_{v}$. Here  
$\psi_{[1, v]}(x_{[1, v]})=\prod_{k\in I}\prod_{ j= v_1^{(k)}+1, \cdots, v^{(k)}} 
\psi_k (x_j^{(k)})$. The subindex $[1, v_1]$ indicates that $\psi_{[1, v]}^{-1}(x_{[1, v]})$ lies in $\textbf{H}_{v}^{\ext}$. 
We extend $S$ to $\textbf{H}^{\ext}_{\loc}$ by requiring $S$ to be an anti-homomorphism. 
(i.e. $S(a\star b)=S(b)\star S(a)$). 
In particular, when $v=e_i$, we have
\[
S(E_i(u)) \mapsto -\psi_{i}^{-1}(x_i) E_i(u). 
\]
Clearly, the above assignment \eqref{psi} defines an anti-homomorphism on $\mathcal{H}_{\loc}^{0}$. 
\begin{lmm}
\begin{enumerate}
\item
Choose 
$a=f(x_{[1, v_1]})\in \textbf{H}_{v_1}$, 
$b=g(x_{[1, v_2]})\in \textbf{H}_{v_2}$. We have $S(a\star b)=S(b)\star S(a)$. 
\item
The map $S$ respects the action \eqref{eq:action}. 
\end{enumerate}
\end{lmm}
{\it Proof.}
For (1): We have
\begin{align*}
&S(g(x_{[1, v_2]})\star S(f(x_{[1, v_1]})\\
=&(-1)^{|v_2|} (-1)^{|v_1|} 
\psi_{[1, v_2]}^{-1}(x_{[1, v_2]}) g(x_{[1, v_2]})\star 
\psi_{[v_2+1, v]}^{-1}(x_{[v_2+1, v]}) f(x_{[v_2+1, v]})\\
=&(-1)^{|v|}  \psi_{[1, v_2]}^{-1}(x_{[1, v_2]}) \psi_{[v_2+1, v]}^{-1} \Big(g(x_{[1, v_2]})\cdot \frac{\fac(x_{[v_2+1, v]} | x_{[1, v]})}{\fac(x_{[1, v]}|x_{[v_2+1, v]} )}\Big)\star 
 f(x_{[v_2+1, v]} ) \\
 =&(-1)^{|v|}  \psi_{[1, v]}^{-1}(x_{[1, v]}) 
\sum_{\sigma\in \Sh(v_2, v_1)}\sigma 
 \Big(g(x_{[1, v_2]})\cdot 
   f(x_{[v_2+1, v]}
 \cdot \fac(x_{[v_2+1, v]} | x_{[1, v_2]})\Big)\\
 =&(-1)^{|v|}  \psi_{[1, v]}^{-1}(x_{[1, v]})  \Big(f(x_{[1, v_1]})\star g(x_{[1, v_2]})\Big)\\
 =& S\Big(f(x_{[1, v_1]})\star g(x_{[1, v_2]})\Big). 
 \end{align*}
 
 We now prove (2). 
 We have 
 \begin{align*}
 S( \psi_{i}(z) f(x_{[1, v]})  \psi_{i}(z)^{-1})
 =& \psi_{i}(z) S(f(x_{[1, v]}))  \psi_{i}(z)^{-1}
 =(-1)^{|v|}\psi_{i}(z)   \psi_{[1, v]}^{-1}(x_{[1, v]})S(f(x_{[1, v]}))  \psi_{i}(z)^{-1}\\
 =&(-1)^{|v|}\psi_{[1, v]}^{-1}(x_{[1, v]}) \psi_{i}(z)   f(x_{[1, v]})  \psi_{i}(z)^{-1}\\
 =&(-1)^{|v|}\psi_{[1, v]}^{-1}(x_{[1, v]})   f(x_{[1, v]}) \frac{\fac(z|x_{[1, v]})}{ \fac(x_{[1, v]}| z)})
 \end{align*}

On the other hand, we compute
 \begin{align*}
 S(f(x_{[1, v]}) \frac{\fac(z|x_{[1, v]})}{ \fac(x_{[1, v]}| z)})
 =(-1)^{|v|}\psi_{[1, v]}^{-1}(x_{[1, v]})   f(x_{[1, v]}) \frac{\fac(z|x_{[1, v]})}{ \fac(x_{[1, v]}| z)}
  \end{align*}
This completes the proof. 
$\blacksquare.$

The above Lemma shows the assignments \eqref{psi} \eqref{psix} determines a unique anti-homomorphism on $ \textbf{H}^{\ext}_{\loc}$. 

\begin{lmm}
The anti-homomorphism $S$ satisfies the axiom $
\mu(S\otimes 1)\Delta=\mu( 1\otimes S)\Delta=\eta\epsilon$. 
\end{lmm} 
{\it Proof.}
We have
\begin{align*}
\mu(S\otimes 1)\Delta(\psi_{i}(w))
=\mu(S\otimes 1) (\psi_{i}(w)\otimes \psi_{i}(w) )
=\mu (\psi_{i}^{-1}(w)\otimes \psi_{i}(w))=1. 
\end{align*}
Similarly, $\mu(1\otimes S)\Delta(\psi_{i}(w))=1$. 

For a dimension vector $v=(v^i)_{i\in I}\in \N^I$, we define $|v|:=\sum_{i\in I }v^i$ and $
{v \choose v_1}:=\prod_{i\in I} {v^i \choose v_1^i}. $
For  $f(x_{[1, v]})\in \textbf{H}_{v}$, we have
\begin{align*}
&\mu(S\otimes 1)\Delta(f(x_{[1, v]}))\\
=&\mu(S\otimes 1) 
\sum_{\{v_1+v_2=v\}} 
\frac{\psi_{[1,v_1]}(x_{[v_1+1,v]}) f(x_{[1,v_1]}\otimes x_{[v_1+1,v]})}{\fac( x_{[v_1+1,v]}| x_{[1,v_1]})}\\
=&\mu \sum_{\{v_1+v_2=v\}} (-1)^{|v_1|}\psi_{[1, v_1]}^{-1}(x_{[1, v_1]}) 
\frac{\psi_{[1,v_1]}^{-1}(x_{[v_1+1,v]}) f(x_{[1,v_1]}\otimes x_{[v_1+1,v]})}{\fac( x_{[1,v_1]} |x_{[v_1+1,v]})}\\
=&\Big(\sum_{\{v_1+v_2=v\}} {v \choose v_1} (-1)^{|v_1|} \Big)\psi_{[1, v]}^{-1}(x_{[1, v]}) f(x_{[1, v]})
=0. 
\end{align*}
In the second equality, the denominator becomes $\fac( x_{[1,v_1]}|x_{[v_1+1,v]})$ since we switch 
$S(\psi_{[1,v_1]}(x_{[v_1+1,v]})=\psi^{-1}_{[1,v_1]}(x_{[v_1+1,v]})$ in front of $f$ using
\[
\psi_{[1, v_1]}^{-1}(x_{[v_1+1,v]}) f(x_{[1, v_1]} \otimes x_{[v_1+1,v]}) =f(x_{[1, v]} \otimes x_{[v_1+1,v]})\psi^{-1}_{[1, v_1]}(x_{[v_1+1,v]})\frac{ \fac(x_{[1, v_1]}| x_{[v_1+1,v]})} {\fac(x_{[v_1+1,v]}|x_{[1, v_1]})} . 
\]
Similarly, we compute
\begin{align*}
&\mu(1\otimes S)\Delta(f(x_{[1, v]}))\\
=&\mu(1\otimes S) 
\sum_{\{v_1+v_2=v\}}
\frac{\psi_{[1,v_1]}(x_{[v_1+1,v]}) f(x_{[1,v_1]}\otimes x_{[v_1+1,v]})}{\fac( x_{[v_1+1,v]}| x_{[1,v_1]})}\\
=& \sum_{\{v_1+v_2=v\}} {v \choose v_1} (-1)^{|v_1|}
{\psi_{[1,v_1]}(x_{[v_1+1,v]})  \psi_{[1, v_1]}^{-1}(x_{[v_1+1,v]}) f(x_{[1,v]})}\\
=&\Big(\sum_{\{v_1+v_2=v\}} {v \choose v_1} (-1)^{|v_1|} \Big) f(x_{[1, v]})=0. 
\end{align*}
This completes the proof. 
$\blacksquare.$

\subsection{The Drinfeld double of COHA}
\label{sec:doubleCoHA}
In this section, we follow the terminology from the book \cite[Chapter 3]{Joseph}. Let $A$ and $B$ be two Hopf algebras. A skew-Hopf pairing of $A$ and $B$ is an $R$-bilinear function
\[
(\cdot, \cdot): A\times B\to R
\]
that satisfies the conditions
\begin{enumerate}
\item[(a)] $(1, b)=\epsilon_B(b), (a,1)=\epsilon_A(a),$
\item[(b)] $(a, bb')=(\Delta_A(a), b\otimes b'),$
\item[(c)] $(aa', b)=(a\otimes a', \Delta_B^{\op}(b))$,
\item[(d)] $(S_A(a), b)=(a, S_B^{-1}(b))$. 
\end{enumerate}

We consider $\textbf{H}^{\ext, \coop}_{\loc}$ that satisfies  $\textbf{H}^{\ext, \coop}_{\loc}=\H^0\otimes \textbf{H}_{\loc}$ as associative algebras, but is equipped with opposed comultiplication described as follows. 
\begin{enumerate}
\item
$\phi_i(z):=(-1)^{l_i+1} \psi_i(z)$, where $l_i$ is the number of loops at vertex $i\in I$. 
\item
The $\H^0$--action on $\textbf{H}_{\loc}$ is by
 \begin{equation}
 \phi_{i}(z) g  \phi_{i}^{-1}(z):=g \frac{\fac(z|x_{[1, v]})}{ \fac(x_{[1, v]}| z)} , \text{for any $g \in \textbf{H}_{v}$}
  \end{equation}
\item The coproduct on $\textbf{H}^{\ext, \coop}_{\loc}$ is given by 
\begin{align*}
\Delta_B: \textbf{H}^{\ext, \coop}_{\loc} &\to \textbf{H}^{\ext, \coop}_{\loc} \otimes \textbf{H}^{\ext, \coop}_{\loc}, \\
\phi_{i}(w)&\mapsto \phi_i(w)\otimes \phi_i(w),\\
g(x_{[1, v]})&\mapsto 
\sum_{\{v_1+v_2=v\}}  
\frac{ g(x_{[1,v_1]}\otimes x_{[v_1+1,v]}) \phi_{[v_1+1, v]}(x_{[1,v_1]}) }{\fac( x_{[v_1+1,v]}| x_{[1,v_1]})}, 
\end{align*}
where $g(x_{[1, v]}))$ is homogeneous and lies in the component $\textbf{H}_{v, \loc}$. 
\item
The antipode $S_B$ is given by
\begin{align*}
&  \phi_i(z) \mapsto \phi_i(z)^{-1}\\
&  g(x_{[1, v]}) \mapsto (-1)^{|v|}   g(x_{[1, v]}) \phi_{[1, v]}^{-1}(x_{[1, v]}) 
\end{align*}
Similar to the proofs in the case of $\textbf{H}^{\ext}_{\loc}$, one easily verify that $\textbf{H}^{\ext, \coop}_{\loc}$  is a Hopf algebra.

\end{enumerate}
Let us take $A$ to be the Hopf algebra $\textbf{H}^{\ext}_{\loc}$ constructed in \ref{H=Hopf}, and $B$ to be the  Hopf algebra $\textbf{H}^{\ext, \coop}_{\loc}$. 
We now construct a bilinear skew-Hopf pairing between $\textbf{H}^{\ext}_{\loc}$ and $\textbf{H}^{\ext, \coop}_{\loc}$ as follows. 
\[
( \cdot , \cdot) : \textbf{H}^{\ext}_{\loc} \otimes \textbf{H}^{\ext, \coop}_{\loc} \to R
\]
\begin{itemize}
\item
$(f_v, g_w)=0$ if $v\neq w$, $(f_v, \phi_i(z))=0$, $(\psi_i(z), g_w)=0$, 
\item
$(\psi_k(u), \psi_l(w))= \frac{\fac(u|w)}{\fac(w|u)}$ for any $k, l\in I$. 
\item
$(f_{e_i}, g_{e_i}):= \Res_{x=\infty} f(x^{(i)})\cdot g(-x^{(i)}) dx$.
\end{itemize}
where $f_v\in \textbf{H}_{v}$, $g_w \in \textbf{H}_{w}^{\coop}$ and $\psi_i(z)\in \SH^0[[z]]$. The dimension vector $e_i$ is defined to be 1 at $i\in I$ and $0$ otherwise. 

We extend the pairing to the entire $\textbf{H}_{\loc, v}\times \textbf{H}_{\loc, v}^{\coop}$ using the property 
$(a, bb')=(\Delta_A(a), b\otimes b')$ and $(aa', b)=(a\otimes a', \Delta_B^{\op}(b))$.
An explicit formula for $(f_v, g_w)$ can be found in \cite{YZ2}.  

We now verify this is a skew-Hopf pairing. 
(a) is obvious; (b)(c) can be proved the same way as in \cite{YZ2}.
\begin{lmm}
$(S_A(a), b)=(a, S_B^{-1}(b))$
\end{lmm}
 {\it Proof.}
 Let $a= \psi_k(z), b=\psi_l(w)$. 
 \begin{align*}
 &(S_A(\psi_k(z)), \psi_l(w))
 = (\psi_k^{-1}(z),  \psi_l(w))
 =\frac{\fac(w|u)}{\fac(u|w)}\\
 &(\psi_k(z), S_B^{-1}(\psi_l(w)))
 = (\psi_k(z),  \psi_l^{-1}(w))
 =\frac{\fac(w|u)}{\fac(u|w)}. 
 \end{align*}
 Thus, $(S_A(\psi_k(z)), \psi_l(w))= (\psi_k(z), S_B^{-1}(\psi_l(w)))$. 
 
 Let $a=f_v, b=\psi_k(z)$.
  \begin{align*}
 (S_A(f_v), \psi_k(z))
 &=((-1)^{|v|}  \psi_{[1, v]}^{-1}(x_{[1, v]}) f(x_{[1, v]}), \psi_k(z))\\&
 =((-1)^{|v|}  \psi_{[1, v]}^{-1}(x_{[1, v]}) \otimes f(x_{[1, v]}), \psi_k(z)\otimes \psi_k(z))=0. \\
(f_v, S_B^{-1}(\psi_k(z)))
 &=(f_v,\psi_k^{-1}(z))
 =0
 \end{align*}
 Thus, $(S_A(f_v), \psi_k(z))= (f_v, S_B^{-1}(\psi_k(z)))$. 
 
  Let $a=f_v, b=g_w$.
    \begin{align*}
 (S_A(f_v), g_w)
 &=((-1)^{|v|}  \psi_{[1, v]}^{-1}(x_{[1, v]}) f(x_{[1, v]}), g_w)\\
 &=\left((-1)^{|v|}  \psi_{[1, v]}^{-1}(x_{[1, v]}) \otimes f(x_{[1, v]}), 
\sum_{\{w_1+w_2=w\}} 
\frac{ \phi_{[w_1+1, w]}(x_{[1,w_1]}) g(x_{[1,w_1]}\otimes x_{[w_1+1,w]})  }{\fac( x_{[w_1+1,w]}| x_{[1,w_1]})}\right)\\
&=\delta_{v, w}(-1)^{|v|} \Big(\psi_{[1, v]}^{-1}(x_{[1, v]}), \phi_{[1, w]}(x_{[1,w]}) \Big) \Big(f(x_{[1, v]}), g(x_{[1, v]})\Big)\\
(f_v, S_B^{-1}g_w)
&=(f_v,  (-1)^{|w|}   \phi_{[1, w]}^{-1}(x_{[1, w]}) g(x_{[1, w]}) )\\&
=\Big(\sum_{\{v_1+v_2=v\}} 
\frac{\psi_{[1,v_1]}(x_{[v_1+1,v]}) f(x_{[1,v_1]}\otimes x_{[v_1+1,v]})}{\fac( x_{[v_1+1,v]}| x_{[1,v_1]})}, 
(-1)^{|w|}   \phi_{[1, w]}^{-1}(x_{[1, w]}) \otimes g(x_{[1, w]}) 
\Big)\\
&=\delta_{v, w}(-1)^{|v|} \Big(\psi_{[1, v]}(x_{[1, v]}), \phi_{[1, w]}^{-1}(x_{[1,w]}) \Big) \Big(f(x_{[1, v]}), g(x_{[1, v]})\Big)
 \end{align*}
 Thus, we have  $(S_A(f_v), g_w)=(f_v, S_B^{-1}g_w)$, since $$(\psi_{[1, v]}^{-1}(x_{[1, v]}), \phi_{[1, w]}(x_{[1,w]}) )=(\psi_{[1, v]}(x_{[1, v]}), \phi_{[1, w]}^{-1}(x_{[1,w]}) ).$$ This completes the proof. 
$\blacksquare.$

It is known that if $A,  B$ are Hopf algebras, endowed with a skew-Hopf paring. Then, there is a unique Hopf structure on $A\otimes B$, called  the Drinfeld double of $(A, B, (\cdot, \cdot))$, determined by the following properties \cite[Lemma~3.2.2]{Joseph}
\begin{enumerate}
\item[(a)] $ (a\otimes 1) (a'\otimes 1)=aa'\otimes 1$, 
\item[(b)] $(1\otimes b) (1\otimes b')=1\otimes bb'$, 
\item[(c)] $(a\otimes 1) (1\otimes b)=a\otimes b$,
\item[(d)] \label{eq:Drinfeld double mul}
$(1\otimes b) (a\otimes 1)=\sum(a_1, S_B(b_1)) a_2\otimes b_2 (a_3, b_3)$,
\end{enumerate}
for all $a, a'\in A, b, b'\in B$. 
Denote by $\Delta^2$ the composition $\Delta\circ \Delta$. 
We use the notation that 
$\Delta_A^2(a)=\sum a_1\otimes a_2\otimes a_3$ and 
$\Delta_B^2(b)=\sum b_1\otimes b_2\otimes b_3$. 

Denote by $D(\textbf{H})$ the Drinfeld double $\textbf{H}^{\ext}_{\loc} \otimes \textbf{H}^{\ext, \coop}_{\loc}$. 
Let $Q'=(I', H') \subset Q$ be a subquiver. That is, $I'\subset I$ is a subset and $H'=\{h\in H \mid \text{the incoming and outgoing vertices of $h$ are in $I'$}\}$. 
It is clear that the shuffle algebra $\textbf{H}'$ associated to $Q'$ is a subalgebra of the shuffle algebra $\textbf{H}$ associated to $I$. 
\begin{prp}
We have $D(\textbf{H}')$ is a subalgebra of $D(\textbf{H})$. 
\end{prp}
{\it Proof.}
The Hopf structure of $\textbf{H}_{\loc}$ when restricted to $\textbf{H}'$ gives a Hopf structure on $\textbf{H}'_{\loc}$ and the restriction of the pairing is still a skew Hopf pairing. The multiplication on $D(\textbf{H}')$ by \eqref{eq:Drinfeld double mul} is compatible  with  the multiplication on  $D(\textbf{H})$. 
This completes the proof. 
$\blacksquare.$

Define $D(\SH)\subset D(\textbf{H})$ 
to be the subalgebra of $D(\textbf{H})$ generated by all the elements in  $\H^0$, $\H_{e_i}$ $\H_{e_i}^{\coop}$, as  $i$ varies in $I$.
\begin{rmk}
By construction, the Drinfeld double $D(\textbf{H})$ is very big, since the coproduct $\Delta$, and the antipode $S$ has poles, and one has to pass to the localization by localizing at the necessary poles. 
We emphasize that the definition of the subalgebra $D(\SH)$ uses only the multiplication of $D(\textbf{H})$ without localization. Note that $D(\SH)$ is not closed under the above coproduct $\Delta$. As vector spaces, we have the decomposition $D(\SH)\cong  \SH\otimes \H^0\otimes   \H^0 \otimes \SH$. 
\end{rmk}
Recall there is an algebra epimorphism 
\[
\mathcal{SH}\to \SH, 
\]
which is an isomorphism after passing to the same localization.

We define the {\it Drinfeld double $D(\mathcal{SH})$} of the spherical COHA as $D(\SH)$. 

\subsection{Relations in the Drinfeld double}

Recall $e_i$ is the dimension vector which is 1 at $i\in I$ and 0 otherwise, $i\in I$. 
Since any 1-by-1 matrices naturally commute and hence the potential function vanishes, we have $H^*_{c,\GL_{e_i}\times T}(\Rep(Q,e_i),\varphi_{\tr W})^\vee\cong R[x^{(i)}]$.
Define $\calH_i$ to be the subalgebra of  $\H$, generated by $R[x^{(i)}]$.

We compute the relations in the Drinfeld double of COHA. First note that $\textbf{SH}^{\ext}$, $\textbf{SH}^{\ext, \coop}$
are two sub-algebras of $D(\mathcal{SH})$. 
\begin{prp}
\label{prop:computeE_iF_j}
Let $E_i(u), F_j(v)$ be the generating series
\[
E_i(u):=\sum_{r\geq 0} (x^{(i)})^{r} u^{-r-1} \in \H^{\ext}[[u^{-1}]], \,\ F_j(v):=\sum_{r\geq 0} (-x^{(i)})^{r} v^{-r-1}\in \H^{\ext,  \coop}[[v^{-1}]], 
\]
for $i, j\in I$.  Then, in $D(\mathcal{SH})$, we have the relation
\begin{displaymath} 
\left\{
     \begin{array}{lr}
       [E_i(u), F_j(v) ] = \delta_{ij}
\frac{\Big( \psi_i (u) \otimes 1-1\otimes \phi_i(u)\Big)- \Big( \psi_i (v) \otimes 1-1\otimes \phi_i(v)\Big) }{u-v}, &  \text{if $\# \{\text{edge loops at $i$ }\}$ is odd} , \\
        \{E_i(u), F_j(v) \} =\delta_{ij}
\frac{\Big( \psi_i (u) \otimes 1-1\otimes \phi_i(u)\Big)- \Big( \psi_i (v) \otimes 1-1\otimes \phi_i(v)\Big) }{u-v},  & \text{if $\# \{\text{edge loops at $i$ }\}$ is even}. 
     \end{array}
   \right.
   \end{displaymath} 
   where $[a, b]:=ab-ba$ is the commutator and $\{a, b\}:=ab+ba$ is the super commutator. 
\end{prp}
{\it Proof.}
By the coproduct formula, we have
\begin{align*}
&\Delta^2(E_i(u))= 
\psi_i (x^{(i)})\otimes  \psi_i (x^{(i)}) \otimes E_i(u)+ 
\psi_i (x^{(i)})\otimes E_i(u)\otimes 1+
E_i(u)\otimes 1\otimes 1
\\
&
\Delta^2(F_j(v))=
F_j(v) \otimes \phi_j(x^{(j)})\otimes \phi_j(x^{(j)})+
1\otimes  F_j(v) \otimes \phi_j(x^{(j)})+
1\otimes 1\otimes F_j(v). 
\end{align*}
In  the multiplication formula \eqref{eq:Drinfeld double mul}, we choose $b=F_j(v)$, $a=E_i(u)$. 
Using the following pairings
\begin{align*}
&(\psi_i (x^{(i)}), 1)=1, \,\  (1, \phi_j(x^{(j)})=(-1)^{l_j+1}, \\
&(E_i(u), -F_j(v)\phi_j^{-1}(x^{(j)}))
=(\psi_i (x^{(i)})\otimes E_i(u)+
E_i(u)\otimes 1, -F_j(v) \otimes \phi_j^{-1}(x^{(j)}))
=-(-1)^{l_j+1}(E_i(u),  F_j(v)), 
\end{align*}
we have
\begin{align*}
&(1\otimes F_j(v)) (E_i(u)\otimes 1)=\sum(a_1, S_B(b_1)) a_2\otimes b_2 (a_3, b_3)
\\
=&(\psi_i (x^{(i)}), 1) \psi_i (x^{(i)}) \otimes 1 ( E_i(u), F_j(v))
+(\psi_i (x^{(i)}), 1) E_i(u)\otimes   F_j(v) ( 1, \phi_j(x^{(j)})\\
&+(E_i(u), -F_j(v)\phi_j^{-1}(x^{(j)}))  1\otimes \phi_j(x^{(j)})(1, \phi_j(x^{(j)}))\\
=&(-1)^{l_j+1} E_i(u)\otimes   F_j(v) +\Big( \psi_i (x^{(i)}) \otimes 1-1\otimes \phi_j(x^{(j)})\Big) ( E_i(u), F_j(v))\\
=&(-1)^{l_j+1} E_i(u)\otimes   F_j(v) +\delta_{ij}\Res_{x^{(i)}=\infty}
\Big( \psi_i (x^{(i)}) \otimes 1-1\otimes \phi_i(x^{(i)})\Big)  \frac{1}{u-x^{(i)}} \cdot \frac{1}{v-x^{(i)}}\\
=&(-1)^{l_j+1} E_i(u)\otimes   F_j(v) +\delta_{ij}
\frac{\Big( \psi_i (u) \otimes 1-1\otimes \phi_i(u)\Big)- \Big( \psi_i (v) \otimes 1-1\otimes \phi_i(v)\Big) }{u-v}
\end{align*}
Here the sign $(-1)^{l_j+1}$ is $1$ (resp. is $-1$), when $j$ is has odd number of edge loops  (resp. even number of edge loops). 
Thus,  we have the desired relation.  $\blacksquare.$

\begin{prp}
\label{prp:Ephi}
The relation between $E_i(u)\otimes1$ and $1\otimes \phi_j(v)$ is given as
\[
(1\otimes \phi_j(v))( E_i(u)\otimes 1)(1\otimes \phi_j^{-1}(v))= E_i(u)\frac{\fac(v|x^{(i)})}{\fac(x^{(i)}|v)}\otimes 1. 
\]
\end{prp}
{\it Proof.}
By the coproduct formula, we have
\begin{align*}
&\Delta^2(E_i(u))= 
\psi_i (x^{(i)})\otimes  \psi_i (x^{(i)}) \otimes E_i(u)+ 
\psi_i (x^{(i)})\otimes E_i(u)\otimes 1+
E_i(u)\otimes 1\otimes 1\\
& \Delta^2(\phi_j(v))=\phi_j(v)\otimes \phi_j(v)\otimes \phi_j(v). 
\end{align*}
In  the multiplication formula  \eqref{eq:Drinfeld double mul}, we choose $b=\phi_j(v), a=E_i(u)$. 
We have
\begin{align*}
&(1\otimes \phi_j(v))(E_i(u)\otimes 1)=\sum(a_1, S_B(b_1)) a_2\otimes b_2 (a_3, b_3)\\
=&(\psi_i(x^{(i)}), \phi_j^{-1}(v))E_i(u)\otimes \phi_j(v) (1, \phi_j(v))\\
=&\frac{\fac(v|x^{(i)})}{\fac(x^{(i)}|v)} E_i(u)\otimes \phi_j(v)
&\end{align*}
This is equivalent to the action
\[
(1\otimes \phi_j(v))( E_i(u)\otimes 1)(1\otimes \phi_j^{-1}(v))= E_i(u)\frac{\fac(v|x^{(i)})}{\fac(x^{(i)}|v)}\otimes 1. 
\]
This completes proof. $\blacksquare.$

\begin{prp}\label{two SH^0}
The relation between $\psi_i(u)\otimes1$ and $1\otimes \phi_j(v)$ is given as
\[
(1\otimes \phi_j(v))(\psi_i(u)\otimes 1)
=(\psi_i(u)\otimes 1)(1\otimes \phi_j(v)). 
\]
\end{prp}
{\it Proof.}
By the coproduct formula, we have
\begin{align*}
& \Delta^2(\psi_i(u))=\psi_i(u)\otimes \psi_i(u)\otimes \psi_i(u), \,\ 
 \Delta^2(\phi_j(v))=\phi_j(v)\otimes \phi_j(v)\otimes \phi_j(v). 
\end{align*}
In  the multiplication formula  \eqref{eq:Drinfeld double mul}, we choose $b=\phi_j(v), a=\psi_i(u)$. 
We have
\begin{align*}
&(1\otimes \phi_j(v))(\psi_i(u)\otimes 1)
\\
=&(\psi_i(u), \phi_j^{-1}(v)) \psi_i(u)\otimes \phi_j(v)  (\psi_i(u),\phi_j(v))
= \psi_i(u)\otimes \phi_j(v). 
\end{align*}
This completes proof. $\blacksquare.$
\subsection{Cartan doubled Yangian and Shifted Yangians}

By Proposition \ref{two SH^0}, in the Drinfeld double $D(\mathcal{SH})$, there is a commutative subalgebra which is isomorphic to $\H^0\otimes \H^0$. For fixed  $i\in I$, label the Cartan elements $\psi_{i,r}\otimes 1 \in \H^0\otimes 1  (r\geq 0)$ and $1\otimes \phi_{i,r} \in 1\otimes \H^{0} (r\geq 0)$ in $D(\mathcal{SH})$ by $\Z$ as follows
\begin{equation}
\label{eq:label by Z}
\psi_{i}^{(r)}:=\left\{
     \begin{array}{lr}
      \psi_{i,r}\otimes 1-1\otimes \phi_{i,r}, &  \text{if $r\geq 0$} , \\
      \psi_{i,-r-1}\otimes 1+1\otimes \phi_{i,-r-1}, &  \text{if $r<0$}. 
     \end{array}
   \right.
\end{equation}

When $\g$ is the Kac-Moody Lie algebra, the Cartan doubled Yangian of $\g$ is introduced in \cite{Finetal} (and recalled in Definition \ref{def Cartan double}). 
In Proposition \ref{prop:KacMoodyDouble}, we show that the 
Cartan doubled Yangian of a Kac-Moody algebra surjects to the Drinfeld double of some $\mathcal{SH}(X_{m, 0})$. 
Note that the index set of $\psi_{i}^{(r)}$ is $I\times \Z$. This explains the name ``Cartan doubled".

A  coweight $\mu$ is a $\Z$-linear $\Z$-valued function on $\bbZ^I$. In \cite{Finetal}, a shifted Yangian of a Kac-Moody algebra \cite[Definition 3.5]{Finetal} is defined to be a quotient of the Cartan doubled Yangian  by a coweight $\mu$. 
Motivated by the loc.cit. we propose the following definition. 
\begin{defn}\label{def:shiftedYangian}
For any coweight $\mu$ we define the shifted Yangian $Y_\mu$ to be the quotient of $D(\mathcal{SH})$ by the relations
\[
\psi_i^{(p)}=0,\ \text{for all $p< -\langle\mu, \alpha_i \rangle$}\ \text{and}\ \psi_i^{(-\langle\mu, \alpha_i \rangle)}=1
\]
\end{defn}
Definition of the shifted Yangian in  \cite[Definition 3.5]{Finetal} is given in the framework of  Kac-Moody Lie algebras. 
Our definition \ref{def:shiftedYangian} of the shifted Yangian $Y_\mu$ generalizes \cite[Definition 3.5]{Finetal} to the setup of an arbitrary quiver with potential. The Kac-Moody case corresponds to the case when the corresponding quiver is a triple quiver with potential as in \cite{Ginz}. See Proposition \ref{prop:KacMoodyDouble} for details. 

\section{Examples: COHA and the Drinfeld double}
\label{Examples COHA}

\subsection{Subalgebras associated to  vertices}\label{subsec:subalg_vertices}

Let $(Q,W)$ be a quiver with potential, and let $i\in I$ be a vertex of the quiver.  Recall $\calH_i$ is the subalgebra of $\calH$ generated by the dimension vector $e_i$. 
\begin{prp} In the above notation the algebra $\calH_i$ is isomorphic to the positive part of the Clifford algebra if $i$ has no edge loop, $Y^+_\hbar(\fs\fl_2)$ if $i$ has 1 edge loop. 
\end{prp}

{\it Proof.}
The COHA of this subquiver in both cases was computed in \cite[\S~2.5]{KS}. In particular it was identified with  the positive part of the infinite Clifford algebra in the first case, see \cite[\S~2.5]{KS}. In the second case it is identified with 
$Y_\hbar^+({\fs\fl}_2)$ (see also \cite{YZ3}).$\blacksquare$.

\begin{exa}
\label{exa:A1}
For  the Dynkin quiver $Q=A_1$ endowed with the potential $W=0$,  the Drinfeld double $D(\H)$ is generated by 
\[
e^{(r)}, f^{(r)}, \psi^{(s)}, r\geq 0, s\in \Z, 
\]
subject to the following relations:
\begin{align*}
& [\psi^{(s)}, \psi^{(t)} ]=0\\
& \{e^{(r_1)}, e^{(r_2)}\}=0, \,\  \{f^{(r_1)}, f^{(r_2)}\}=0,\\
& \{e^{(r_1)}, \psi^{(r_2)}\}=0,\,\  \{f^{(r_1)}, \psi^{(r_2)}\}=0,\\
& \{e^{(r_1)}, f^{(r_2)}\}=\psi^{(r_1+r_2)}, 
\end{align*}
where $\{a, b\}:=ab+ba$, and $\psi^{(s)}$ is defined in \eqref{eq:label by Z}. 
Since the $A_1$ quiver has no arrows, the conjugation action \eqref{eq:action} of $\psi(z)$ on $\H$ is given as 
\[
 \psi(z) g (\psi(z))^{-1} =\frac{\fac(z|x)}{\fac(x|z)}g = \frac{z-x}{x-z}g=-g. 
\] 
Proposition \ref{prp:Ephi} implies $\psi(z) g (\psi(z))^{-1}=-g$. Therefore, $\{\psi^{(s)}, e^{(r)}\}=0$, for $s\in \Z$.  
Similarly, $\{\psi^{(s)}, f^{(r)}\}=0$, for $s \in \Z$. Proposition \ref{prop:computeE_iF_j}  implies that $\{e^{(r_1)}, f^{(r_2)}\}=\psi_{r_1+r_2}$.  
\Omit{
The Clifford algebra presentation in \cite[P53]{X} is 
\begin{align*}
& h_s \,\ \text{is central, for all $s\in \Z$},\\
& e_{r_1}e_{r_2}+e_{r_2}e_{r_1}=0,\\
& f_{r_1}f_{r_2}+f_{r_2}f_{r_1}=0,\\
& e_{r_1}f_{r_2}+f_{r_2}e_{r_1}=\delta_{r_1, r_2}h_{r_1}.
\end{align*}
}
\end{exa}
\begin{exa}
\label{exa:Jordan quiver}
Let $Q$ be the Jordan quiver, i.e. it has  one vertex and one loop. Then for the potential $W=0$ the Drinfeld double $D(\H)$ is generated by 
\[
e^{(r)}, f^{(r)}, \psi^{(s)}, r\geq 0, s\in \Z, 
\]
subject to the Cartan doubled Yangian relations for $\mathfrak{sl}_2$ (see Definition \ref{def Cartan double}, taking $I$ to be a point). 
In particular, we have the relation
\[
[e(u), f(v)]=  \frac{\psi^+(u)-\psi^+(v) }{u-v}, 
\]
where $\psi^+(u)=\psi(u)\otimes 1-1\otimes \psi(u)$ which follows from Proposition \ref{prop:computeE_iF_j}. 
\end{exa}

\subsection{The Cartan doubled Yangian of a Kac-Moody Lie algebra}\label{sec:double_KacMoody}
Let us start with the following definition.
\begin{defn} \cite[Definition 3.1]{Finetal}
\label{def Cartan double}
Let $\g$ be a Kac-Moody Lie algebra.  The Cartan double Yangian $Y_{\infty}(\g)$ is the $\C$-algebra generated by $E_i^{(q)}, F_i^{(q)}, H_i^{(p)}$, for $i\in I,  q>0, p\in \Z$, subject to the relations
\begin{align}
& [H_i^{(p)}, H_{j}^{q}]=0, \label{HH} \tag{HH}\\
& [E_i^{(p)}, F_{j}^{(q)}]=\delta_{ij} H_i^{(p+q-1)},  \label{EF} \tag{EF}\\
&[H_i^{(p+1)}, E_j^{(q)}]-[H_i^{(p)}, E_j^{(q+1)}]=\frac{\alpha_i \cdot \alpha_j}{2}(H_i^{(p)} E_j^{(q)}+E_j^{(q)} H_i^{(p)}), \label{HE} \tag{HE} \\
&[H_i^{(p+1)}, F_j^{(q)}]-[H_i^{(p)}, F_j^{(q+1)}]=-\frac{\alpha_i \cdot \alpha_j}{2}(H_i^{(p)} F_j^{(q)}+F_j^{(q)} H_i^{(p)}),  \label{HF} \tag{HF}\\
&[E_i^{(p+1)}, E_j^{(q)}]-[E_i^{(p)}, E_j^{(q+1)}]=\frac{\alpha_i \cdot \alpha_j}{2}(E_i^{(p)} E_j^{(q)}+E_j^{(q)} E_i^{(p)}), \label{EE} \tag{EE} \\
&[F_i^{(p+1)}, F_j^{(q)}]-[F_i^{(p)}, F_j^{(q+1)}]=-\frac{\alpha_i \cdot \alpha_j}{2}(F_i^{(p)} F_j^{(q)}+F_j^{(q)} F_i^{(p)}),  \label{FF} \tag{FF}\\
&i\neq j, N=1-\alpha_i\cdot \alpha_j, \text{sym}[E_{i}^{(p_1)}, [E_i^{(p_2)},  \cdots,  [E_i^{(p_N)}, E_j^{q}]\cdots ]]=0, \notag\\
&i\neq j, N=1-\alpha_i\cdot \alpha_j, \text{sym}[F_{i}^{(p_1)}, [F_i^{(p_2)},  \cdots,  [F_i^{(p_N)}, F_j^{q}]\cdots ]]=0, \notag
\end{align}
\end{defn}

Now let $\fg$ be a symmetric Kac-Moody Lie algebra with Dynkin diagram $\Gamma$.  Let $(Q,W)$ be the tripled quiver with potential as in \cite{Ginz}. Let $\calH$ be the COHA of $(Q,W)$.

\begin{prp}\label{prop:KacMoodyDouble}
 In the above notation  we have an epimorphism
\[
Y_{\infty}(\g) \surj D(\mathcal{SH}). 
\]
\end{prp}
{\it Proof.}
Consider the following generating series of $H_i^{(p)} (p\in \Z)$
\begin{align*}
H_i^{>}(z):=&\sum_{p \geq 0} H_{i}^{(p+1)} z^{-p-1}, \,\ 
H_i^{<}(z):=1+\sum_{p \leq  0} H_{i}^{(p)} z^{-p-1}. 
\end{align*}
Define a map $Y_{\infty}(\g) \surj D(\mathcal{SH})$ by 
\begin{align*}
E_i^{(q)} &\mapsto  (x^{(i)})^{q-1} \otimes 1\in \mathcal{SH}^{\ext}\otimes 1, \\
F_i^{(q)} &\mapsto 1\otimes (-x^{(i)})^{q-1}  \in  1\otimes \mathcal{SH}^{\ext}\\
H_i^{>}(z) & \mapsto  \psi_{i}^+(z)=\psi_i(z)\otimes 1-1\otimes \psi_i(z),\\
H_i^{<}(z) & \mapsto  \psi_{i}^-(z)=\psi_i(z)\otimes 1+1\otimes \psi_i(z).
\end{align*}
By \cite[Theorem 7.1]{YZ1}, the map respects the relations of \eqref{EE}, \eqref{FF} and the two Serre relations. 
By Proposition \ref{two SH^0}, the map preserves the relation  \eqref{HH}. Proposition \ref{prop:computeE_iF_j} implies the map respects the relation \eqref{EF}. 
Using the action \eqref{eq:action} and Proposition \ref{prp:Ephi}, the map respects the relations of \eqref{HE}, \eqref{HF}. This completes proof. $\blacksquare.$

\subsection{Quiver with potential coming from toric Calabi-Yau 3-folds}\label{sec:toricCY}
With notations as in  \S~\ref{intro}, let $f: X\to Y$  be a resolution of singularities, where 
 $X$ is a smooth toric Calabi-Yau 3-fold and $Y$ is affine. 
There is a tilting generator $\calP$ of $D^b\Coh(X)$
 so that the functor  \[R\Hom_X(\calP,-):D^b\Coh(X)\to D^b(\Mod \calA_0)\]
 induces an equivalence of derived categories.
 Here $\calA_0:=\End_X(\calP)^{\op}$ is a coherent sheaf of non-commutative $\calO_Y$-algebras, and $\Mod \calA_0$ is the abelian category of coherent sheaves of (right) modules over $\calA_0$. 
Let $\mathscr{A}\subset D^b \Coh(X)$ be a heart of $D^b \Coh(X)$ corresponding to the  heart  $\Mod\calA_0$ of the standard $t$-structure of $D^b(\Mod\calA_0)$ under the functor $R\Hom_X(\calP,-)$. One example of $\mathscr{A}$ is $\Perv^{-1}(X/Y)$, the abelian category of perverse coherent sheaves in the sense of Bridgeland and Van den Bergh \cite{Br,vdB}. 
 It is known that $\calA_0$ can be identified with the Jacobian algebra of a quiver with potential $(Q,W)$. 
We denote the COHA of the pair $(Q,W)$ also be $\calH_X$.

Note that each vertex $i\in I$ of $Q$ gives rise to a simple object $S_i$ in the heart of the induced $t$-structure on $D^b\Coh(X)$. We say $S_i$ is {\it bosonic} if the simple object $S_i$  in $D^b\Coh(X)$ has the property that $\Ext^*(S_i, S_i)=H^{2*}(\PP^3)$; we say $S_i$ is {\it fermionic} if the simple object $S_i$ is a spherical object, that is, $\Ext^*(S_i, S_i)=H^*(S^3)$.  
We have the following general proposition. 

\begin{prp}\label{Prop:SimpleVertices}
For a vertex $ i\in I$ of the quiver $Q$, assume the corresponding  simple object $S_i$ is a line bundle on $j_*\calO_{\bbP^1}(m)[n]$ with $j:\bbP^1\inj X$ and $m,n\in\Z$. Then, $S_i$ is either bosonic or fermionic.
\end{prp}
{\it Proof.}
It is not difficult to calculate the self-extension of this line bundle. Without loss of generality, assume this line bundle is trivial. Taking the formal neighbourhood of $\bbP^1$ in $X$, we get a rank-2 vector bundle on $\bbP^1$ the determinant of which is $\calO(-2)$. By the assumption that $X$ is birational to its affinization $Y$, we get that this rank-2 vector bundle has no ample line sub-bundles, therefore is either $\calO\oplus\calO(-2)$ or $\calO(-1)\oplus\calO(-1)$.
 In the first case, $\calO_{\bbP^1}$ has 1-dimensional self-extension. In the second case, $\calO_{\bbP^1}$ has no non-trivial self-extensions. The entire $\Ext^*(S_i, S_i)$ is then determined by the Calabi-Yau condition.$\blacksquare.$

Therefore, when $S_i$ is bosonic, it follows from Example \ref{exa:Jordan quiver} that $D(\H_i)$ is the Cartan doubled Yangian of $\fs\fl_2$. When  $S_i$ is fermionic, it follows from Example \ref{exa:A1} that $D(\H_i)$ is the infinite Clifford algebra.

We give three examples of quiver with potentials coming from this setup.

\begin{exa}\label{ex:m_0}
Let $X=X_{m, 0}$. Then, the defining equation $xy=z^m$ of $Y_{m, 0}$ is a type $A$-singularity. 
All the simple objects are bosonic.  Indeed, let $S_m=\C^2/\Z_m$ be the type $A_{m-1}$ Kleinian surface singularity, and $\tilde{S}_m$ its crepant resolution.  Then, $X_{m, 0}\cong \tilde{S}_m\times\C$. The toric diagram of $X_{m, 0}$ is depicted in figure \ref{fig1}.  The quiver with potential is the tripled quiver of $\widehat{A_{m-1}}$ in the sense of of \cite{Ginz}.
\begin{figure}
  \centering
    \includegraphics[width=0.28\textwidth]{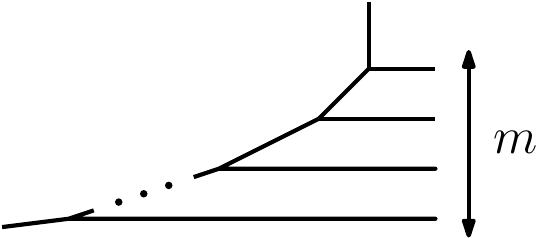}
  \caption{The toric diagram of $X_{m,0}$ consists of $m$ horizontal lines with orientation $(1,0)$ attached from the right to the sequence of lines with orientations $(k,1)$ for $k=0,\dots, n$.}
\label{fig1}
\end{figure}
\end{exa}

\begin{exa}\label{ex:m_n}
Let $X=X_{m, n}$. A toric diagram of one possible resolution of $Y_{m, n}$ is shown in figure \eqref{fig2}. Other possibilities can be obtained by considering different orderings of the left and the right $(1,0)$ lines ending on a sequence of $(k,1)$ lines. All possible choices for $Y_{2,1}$ are depicted in figure \eqref{fig3}. One can then associate a bosonic simple object to each internal line bounded by two $(1,0)$ lines ending from the same side and a fermionic simple object to each internal line bounded by $(1,0)$ lines ending from opposite directions. 
\begin{figure}
  \centering
    \includegraphics[width=0.56\textwidth]{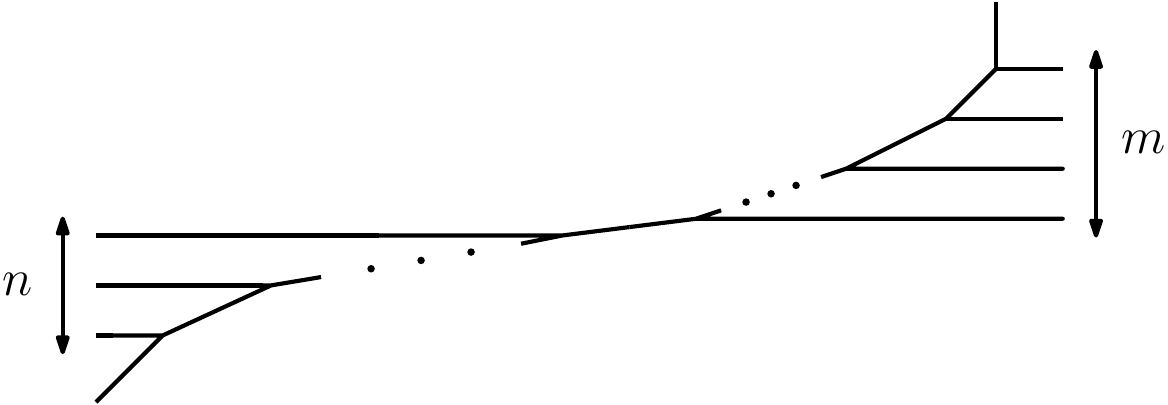}
  \caption{The toric diagram of resolved $Y_{m,n}$ consists of $m$ horizontal lines with orientation $(1,0)$ attached from the right together with $m$ horizontal lines with orientation $(1,0)$ attached from the left to the sequence of lines with orienations $(k,1)$. The figure shows one possible resolution with $m$ lines ending from the left followed by $n$ lines ending from the right. Other possible resolutions are associated to different orderings of the left and the right lines.}
\label{fig2}
\end{figure}
\begin{figure}
  \centering
    \includegraphics[width=0.53\textwidth]{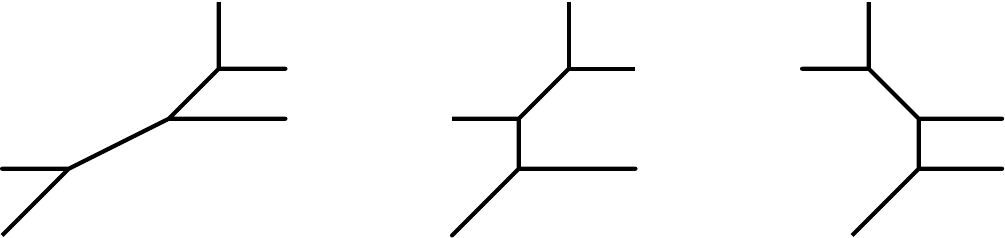}
  \caption{Three possible orderings of horisontal lines for $Y_{2,1}$ associated to the two possible choices of the root system of $\mathfrak{gl}(2|1)$. The root system associated to the left figure consists of a single bosonic and a single fermionic root whereas the left figure is associated with the root system consisting of two fermionic roots. The third possible ordering is a simple reflection of the first one.}
\label{fig3}
\end{figure} The quiver with potential is given explicitly in \cite{Nagao}
\end{exa}

However, there are more examples which have not been well studied in mathematical literature. 
\begin{exa}\label{ex:nhat}
Let $X=S\times \C$, where $S$ is an elliptically fibered $K3$-surface, with special fiber $C$ being 
a collection of $\PP^1$'s in the $\widehat{A_n}$-type configuration. That is, let $\Delta$ be the collection irreducible components of $C$ each of which is isomorphic to $\PP^1$. There is one component of $C$, removing which results a curve $C^o$ an open neighborhood of which is isomorphic to $\tilde{S}_{m}\times\C$ as in Example~\ref{ex:m_0} above. Although there is no contraction of this collection of $\PP^1$'s with rational singularity, hence the method of constructing a tilting bundle from \cite{vdB} does not apply here, the result of \cite{GKV} still suggest that there is a $t$-structure on $D^b\Coh(X)$ that contains $\calO_{C_i}(1)[-1]$ with $i\in I$ and $\calO_{C^o}$ as simple objects. It is easy to verify that all these simple objects are bosonic. 
 \end{exa}

\subsection{Shifted Yangians and perverse coherent systems}
\label{subsec:pervsys}
Let us explain how the notion of shifted Yangians appears later on in our paper.
We will consider the moduli space of stable \textit{perverse coherent system} on $X$. 
As to be explained in \ref{sec:D4 branes} and \ref{sec:D4 branes2},  the moduli space depends not only on a choice of the stability parameter $\zeta$, but also on an algebraic cycle $\chi$ of the form \[
N_0 X+\sum_{i=1}^r N_iD_i,\] where $N_i\in\mathbf N$,  and each $D_i$ is a toric divisor. Here $\chi $ determines a framing of the quiver with potential $(Q,W)$. \footnote{In  physics terminology $\chi$ should be thought of as a stack of $D6$-branes and $D4$-branes wrapping the corresponding cycles. }
We mark the framings associated to $[X]$ by $\square$ and the ones associated to each $D_i$ by $\triangle$. 
For example, when $X=X_{2,0}$, and $D_1$ is the fiber of the $0\in\PP^1$, the framed quiver is as follows 

\begin{equation}\label{eqn:orbifold_chainsaw}
\begin{tikzpicture}[scale=0.8]
        	\node at (-2.2, 0) {$\bullet$};   	\node at (-2.2, 0.5){$V_0$}; \node at (-3.5, 0){$B_3$};
	\node at (2.2, 0) {$\bullet$};  \node at (2.2, 0.5){$V_1$};\node at (3.5, 0){$\tilde{B}_3$};
		\draw[-latex,  bend left=30, thick] (-2, 0.1) to node[above]{$B_2$} (2, 0.1);
		\draw[-latex,  bend left=50, thick] (-2, 0.3) to node[above]{$B_1$} (2, 0.3);
	\draw[-latex,  bend left=30, thick] (2, -0.1) to node[above]{$\tilde{B}_1$} (-2, -0.1);
		\draw[-latex,  bend left=50, thick] (2, -0.3) to node[above]{$\tilde{B}_2$} (-2, -0.3);
				\path (2.2, 0) edge [loop right, min distance=2cm, thick, bend right=40 ] node {} (2.2, 0);
				\path (-2.2, 0) edge [loop left, min distance=2cm, thick, bend left=40 ] node {} (-2.2, 0);
\node at (0, -4) {$\triangle$};  \node at (0.5, -4) {$N_1$};  
 \node at (-1.5, -2) {$I_{13}$};   \node at (1.5, -2) {$J_{13}$};  
 \draw[->, thick] (0, -3.8) -- (-2.2, -0.2) ;
 \draw[->, thick]  (2.2, -0.2) --(0, -3.8) ;

 \node at (-2.2, -2) {$\square$}; 
  \draw[->, thick] (-2.2, -1.8) -- (-2.2, -0.2) ;
 \node at (-2.5, -1) {$\iota$}; 
  	\node at (-2.7, -2) {$N_0$}; 
    \end{tikzpicture}
\end{equation}    

 We are going to discuss the framed quiver with potential for general $\chi$ in the future publications. In this paper we only consider some special cases. In particular, in \S~\ref{sec:D4 branes2} we will discuss   how in general $\zeta$ and $\chi$ determine a coweight $\mu$. We also explain that we expect that the action of $D(\mathcal{SH})$ on cohomology of perverse coherent systems factors through $Y_\mu$.
Furthermore, we expect that the algebra $D(\mathcal{SH})$ has a  coproduct, which after passing to the quotient becomes $Y_\mu\to Y_{\mu_1}\otimes Y_{\mu_2}$ if $\mu_1+\mu_2=\mu$. This coproduct is compatible with hyperbolic restriction to fixed points of a subtorus acting on the framings, similar to \cite{RSYZ}.  Therefore, we can reduce the verification of the action of $Y_\mu$ for a general $\chi$ to the cases when the framing dimension is 1. 

\subsection{COHA of $\C^3$}
\label{subsec:C3}
Consider the case when $X=\C^3$. The quiver with potential is recalled in the Introduction. We denote the quiver by $Q_3$ and the potential  by $W_3$. 
The algebraic properties of the COHA has been studied \cite[\S~2.3]{RSYZ}, building up on earlier works \cite{SV,FT}. 
In particular, there is a subalgebra of the equivariant COHA of $(Q_3, W_3)$ which is called   the {\it equivariant spherical COHA} of $(Q_3, W_3)$ (or  ${\bf C}^3$). Recall that in  \cite[Theorem~7.1.1]{RSYZ} the equivariant spherical COHA  was identified with  $\affY^+$. 

Furthermore, using the same coproduct as in \S~\ref{sec:doubleCoHA}, one can show the reduced Drinfeld  double of $\mathcal{SH}_{\C^3}$ is isomorphic to the entire Yangian $\affY$ \cite{RSYZ}.
Similar argument as in  Proposition~\ref{prop:KacMoodyDouble} then implies that the shifted Yangians $Y_l(\vec{z})$ from \S~\ref{sec:shiftedYangian} 
are quotients of $D(\mathcal{SH}_{\C^3})$.

\subsection{COHA of the resolved conifold}

\label{subsec:conifold}
In this section, the toric Calabi-Yau $3$-fold is $X=X_{1,1}$, the resolved conifold. The corresponding quiver with potential was described in the Introduction (see \eqref{eq:X11}). Recall that we denote the quiver by $Q_{1,1}$ and the potential  by $W_{1,1}$, and we denote the corresponding COHA  by $\H_{X_{1,1}}$. By Proposition~\ref{Prop:SimpleVertices} there are positive parts of two Clifford subalgebras in $\H_{X_{1,1}}$ associated to the two vertices of the quiver $Q_{1,1}$.
The goal of this section is to show that  there is an algebra homomorphism $\H_{X_{1,1}}\to \H_{\C^3}$, where $ \H_{\C^3}$ is the COHA for $\C^3$ from   \S~\ref{subsec:C3}.

For any $n\in\bbN$ consider the open subset $\Rep(Q_{1,1}, (n,n))^0\subset \Rep(Q_{1,1}, (n,n))$ consisting of such representations that the map $b_1$ is an isomorphism. 
We will show below that 
$\oplus_{n\in \N} H^*_{c,\GL_{(n,n)}\times T}(\Rep(Q_{1,1}, (n,n))^0,\varphi_{\tr W_{1,1}})^\vee$ has a natural algebra structure, so that the restriction to open subset 
\[\oplus_{n}H^*_{c,\GL_{(n,n)}\times T}(\Rep(Q_{1,1}, (n,n)),\varphi_{\tr W_{1,1}})^\vee \to \oplus_{n}H^*_{c,\GL_{(n,n)}\times T}(\Rep(Q_{1,1}, (n,n))^0,\varphi_{\tr W_{1,1}})^\vee\] 
is an algebra homomorphism.

Observe that there is a canonical isomorphism of vector spaces 
$$\oplus_n H^*_{c,\GL_{(n,n)}\times T}(\Rep(Q_{1,1}, (n,n))^0,\varphi_{\tr W_{1,1}})^\vee\cong \oplus_n H^*_{c,\GL_{n}\times T}(\Rep(Q_{3}, n),\varphi_{\tr W_{3}})^\vee. $$
Furthermore, with the algebra structures on both sides, the above isomorphism of vector spaces is an isomorphism of algebras 
with the latter $ \oplus_n H^*_{c,\GL_{n}\times T}(\Rep(Q_{3}, n),\varphi_{\tr W_{3}})^\vee$ carrying the multiplication of COHA of $\C^3$.

\begin{prp}
Restriction to the above-defined open subset induces the algebra homomorphism $\oplus_{n\in \mathbb{N}} \H_{X_{1,1}}(n, n) \to  \H_{\C^3}$. Its  image  contains $\affY^+$.
\end{prp}

{\it Proof.}
Consider the affine space parameterizing representations of $Q_{1,1}$ on the underlying vector space $\underline{V}=(V_1,V_2)$ with maps $a_i$ from $V_1$ to $V_2$ and $b_i$ backwards for $i=1,2$. We denote such space by $\Rep(Q_{1,1}, \underline{V})$.
We first recall the multiplication of $\H_{X_{1,1}}$. 
For this, consider pairs of flags 
\[
\text{$0\to V_1\subseteq V_1'\surj V_1''\to 0$ and $0\to V_2\subseteq V_2'\surj V_2''\to 0$.}
\] Fixing one choice of such a pair, the space parameterizing all such pairs can be identified with $G/P$ where $G:= \GL(V_1')\times \GL(V_2')$, and  
$P=P_1\times P_2$ is the parabolic subgroup with $P_i=\{x\in \GL(V_i')\mid x(V_i) \subset V_i\}, i=1, 2$.
Consider $\Rep(Q_{1,1},\underline{V'})$ and the subspace $Z$ consisting of representations such that  $(V_1,V_2)$ is a sub-representation. 
We have the following diagram of correspondences
\[
\xymatrix{
&G\times_P Z \ar[ld]_(0.4){\phi} \ar@{^{(}->}[r]^(0.3){\eta}&  G\times_P  \Rep(Q_{1,1},\underline{V}') \ar[d]^{\psi}  \\
G\times_P \left({\begin{matrix}\Rep(Q_{1,1}, \underline{V}) \\ \times  \Rep(Q_{1,1},\underline{V}'')\end{matrix}}\right) & &  \Rep(Q_{1,1}, \underline{V}') 
}
\]
where $\phi$ is an affine bundle. The potential $W_{1, 1}$ defines functions  $\tr ((W_{1,1})_{\underline{V}}\boxplus (W_{1,1})_{\underline{V}''})$ on 
$\Rep(Q_{1,1}, \underline{V}) \times  \Rep(Q_{1,1},\underline{V}'')$ and $\tr((W_{1, 1})_Z)$ on $Z$ by restricting the function $\tr ((W_{1,1})_{\underline{V}'})$. 
We have $\tr ((W_{1,1})_{\underline{V}}\boxplus (W_{1,1})_{\underline{V}''})\circ \phi=\tr((W_{1, 1})_Z)$.

The map $\eta$ is a closed embedding and $\psi$ is a projection, both of which are compatible with the potential functions. 

Let $Z^0$ be the intersection of $Z$ with $\Rep(Q_{1,1}, \underline{V}')^{{0}}$. 
Note that $P$ acts on $Z^0$, and that the natural map $\eta^0:Z^0\to \Rep(Q_{1,1}, \underline{V}')^0$ is a pullback
\[\xymatrix{
Z^0\ar[r]\ar[d]& \Rep(Q_{1,1}, \underline{V}' )^0\ar[d]\\
Z\ar[r]& \Rep(Q_{1,1}, \underline{V}')}\] 
In particular, restriction to the open subsets intertwines $\eta_*$ and $(\eta^0)_*$. That is, we have the commutative diagram
\[\xymatrix{
H_{c, G\times T}^*(G\times_P Z,\varphi_{\tr W_{1,1}})^{\vee} \ar[r]^(0.4){\eta_*} \ar[d] & H_{c, G\times T}^*(G\times_P \Rep(Q_{1,1}, \underline{V}'),\varphi_{\tr W_{1,1}})^{\vee} \ar[d]\\
H_{c, G\times T}^*(G\times_P Z^0,\varphi_{\tr W_{1,1}}) ^{\vee} \ar[r]^(0.4){\eta_*^0}  & H_{c, G\times T}^*(G\times_P \Rep(Q_{1,1}, \underline{V}' )^0,\varphi_{\tr W_{1,1}})^{\vee}
}
\] 
Similarly, $\psi^0: G\times_P\Rep(Q_{1,1}, \underline{V}')^0\to \Rep(Q_{1,1},\underline{V}')^0$ is pullback of $\psi$, and hence after taking cohomology  the restrictions to the open subsets intertwines $\psi_*$ and $(\psi^0)_*$. 
\[\xymatrix{
H_{c, G\times T}^*(G\times_P \Rep(Q_{1,1}, \underline{V}'),\varphi_{\tr W_{1,1}})^{\vee} \ar[r]^(0.55){\psi_*} \ar[d] &   H_{c, G\times T}^*( \Rep(Q_{1,1}, \underline{V}' ),\varphi_{\tr W_{1,1}})^{\vee}\ar[d]\\
 H_{c, G\times T}^*(G\times_P \Rep(Q_{1,1}, \underline{V}' )^0,\varphi_{\tr W_{1,1}})^{\vee} \ar[r]^(0.55){(\psi^0)_*}  & H_{c, G\times T}^*(\Rep(Q_{1,1}, \underline{V}' )^0,\varphi_{\tr W_{1,1}})^{\vee}
}
\] 
Also, clearly the map $\phi$ induces the map  $\phi^0$ by restriction to open subsets as in the following diagram. Again after taking cohomology restrictions to the open subsets intertwines $\phi^*$ and $(\phi^0)^*$. So the following diagram commutes. 
\[
\xymatrix{
H^*_{c, G\times T}(G\times_P \left(\Rep(Q_{1,1}, \underline{V})  \times  \Rep(Q_{1,1},\underline{V}'')\right), \varphi_{\tr W_{1,1}})^{\vee} \ar[r]^(0.7){\phi^*} \ar[d]  & H^*_{c,G\times T}(G\times_P Z, \varphi_{\tr W_{1,1}})^{\vee}  \ar[d]\\
H^*_{c, G\times T}(G\times_P (\Rep(Q_{1,1}, \underline{V})^0  \times  \Rep(Q_{1,1},\underline{V}'')^0), \varphi_{\tr W_{1,1}})^{\vee} \ar[r]^(0.7){(\phi^0)^*}  & H^*_{c, G\times T}(G \times_P Z^0, \varphi_{\tr W_{1,1}})^{\vee}   
}
\]

To summarize, we have the following commutative diagram
\begin{equation}\label{diag:multi}
\xymatrix@C=1em{
&G\times_P Z^0 \ar[ld]_(0.4){\phi^0} \ar@{^{(}->}[r]^(0.4){\eta^0} \ar@{^{(}->}[d]&  G\times_P  \Rep(Q_{1,1},\underline{V}')^0 \ar[dr]^{\psi^0} \ar@{^{(}->}[d]& \\
G\times_P \left({\begin{matrix}\Rep(Q_{1,1}, \underline{V})^0 \\ \times  \Rep(Q_{1,1},\underline{V}'')^0\end{matrix}}\right) \ar@{^{(}->}[d]
&G\times_P Z \ar[ld]_(0.4){\phi} \ar@{^{(}->}[r]^(0.4){\eta}&  G\times_P  \Rep(Q_{1,1},\underline{V}') \ar[dr]^{\psi}  &\Rep(Q_{1,1}, \underline{V}')^0 \ar@{^{(}->}[d]\\
G\times_P \left({\begin{matrix}\Rep(Q_{1,1}, \underline{V}) \\ \times  \Rep(Q_{1,1},\underline{V}'')\end{matrix}}\right) & & & \Rep(Q_{1,1}, \underline{V}')  
}
\end{equation}
For $\underline{V}'=(V_1',V_2')$, we now impose a condition that $b_1: V_2' \to V_1'$ is an isomorphism. 
The space of isomorphisms $Isom(V_2', V_1')$ is a torsor over $\GL(V_1')$. 
In particular, $\GL(V_1')$ acts transitively on $\Rep(Q_{1,1},\underline{V}')^0$ by the change of basis on $V_1'$. The quotient 
\[
\Rep(Q_{1,1},\underline{V}')^0/\GL(V_1')\cong \Rep(Q_3,V_1')
\]
is canonically identified with $\Rep(Q_3,V_1')$ with $B_1=a_1\circ b_1$,  $B_2=a_2\circ b_1$, and $B_3=b_1^{-1}\circ b_2$. The action of $\GL(V_2')$ on the left hand side and the action of $\GL(V_1')$ on the right hand side are compatible. 

Similarly, the parabolic subgroup $P_1$ acts on $Z^0$ transitively. Thus, we have the identification
\[
(G\times_P Z^0)/\GL(V_1')=GL(V_2')\times_{P_2} (Z^0/P_1) \cong GL(V_2')\times_{P_2} Z_{Q_3} \cong GL(V_1')\times_{P_1} Z_{Q_3}, 
\] where $Z_{Q_3}$ is the correspondence used in the multiplication of $\calH_{Q_3}$. 

Note that the action of $P$ on $\Rep(Q_{1,1}, \underline{V})^0 \times  \Rep(Q_{1,1},\underline{V}'')^0$ factors through the projection 
$P\surj  GL(\underline{V}) \times GL(\underline{V}'')$. Therefore, we have 
\begin{align*}
&\frac{G\times_P(\Rep(Q_{1,1}, \underline{V})^0 \times  \Rep(Q_{1,1},\underline{V}'')^0)}{GL(V_1')}
\cong GL(V_2')\times_{P_2}\frac{\Rep(Q_{1,1}, \underline{V})^0 \times  \Rep(Q_{1,1},\underline{V}'')^0}{P_1}\\
\cong &GL(V_2')\times_{P_2}\frac{\Rep(Q_{1,1}, \underline{V})^0 \times  \Rep(Q_{1,1},\underline{V}'')^0}{GL(V_1) \times GL(V''_1)}\\
\cong &GL(V_2')\times_{P_2}  \Big(\Rep(Q_3,V_1) \times \Rep(Q_3,V_1'')\Big)
\cong GL(V_1')\times_{P_1}  \Big(\Rep(Q_3,V_1) \times \Rep(Q_3,V_1'')\Big)
\end{align*}
To finish the proof, it remains to notice that quotient of the top correspondence of diagram \ref{diag:multi} by $\GL(V_1')$ becomes
 \[
\xymatrix{
GL(V_1')\times_{P_1}  \left({\begin{matrix} \Rep(Q_3,V_1) \\ \times \Rep(Q_3,V_1'') \end{matrix}}\right)&GL(V_1')\times_{P_1} Z_{Q_3} \ar[l]_(0.4){\phi_{Q_3}} \ar[r]^{\psi_{Q_3} \circ \eta_{Q_3}}&  \Rep(Q_3,V_1') 
}
\]
This is the correspondence that defines multiplication of $\calH_{\C^3}$. 

The image of $\oplus_{n\in \mathbb{N}} \H_{X_{1,1}}(n, n) \to  \H_{\C^3}$ contains the spherical COHA $\mathcal{SH}_{\C^3}$. In \cite[Theorem~7.1.1]{RSYZ} $\mathcal{SH}_{\C^3}$  was identified with  $\affY^+$. Thus, the image contains $\affY^+$. 

This finishes the proof. 
$\blacksquare.$

\section{Shifted Yangians of $\widehat{{\mathfrak{gl}}(1)}$  and Costello conjectures}\label{sec:shiftedYangian}

\subsection{Definition of the shifted Yangians of $\widehat{{\mathfrak{gl}}(1)}$}
Let $\hbar_1, \hbar_2, \hbar_3$ be the deformation parameters with $\hbar_1+\hbar_2+\hbar_3=0$. 
Set $\sigma_2:= \hbar_1\hbar_2+\hbar_1\hbar_3+\hbar_2\hbar_3$, and $\sigma_3=\hbar_1\hbar_2\hbar_3$. 
Let $\vec{z}=(z_1, z_2, \cdots, z_{|l|})$ be the parameters ($l\in \Z$). 
\begin{defn}
\label{l shifted}
Let $Y_{l}(\vec{z})$ be a $\C[\hbar_1, \hbar_2]$-algebra generated by $D_{0, m} (m\geq 1)$, $e_{n}, f_{n}$ $(n\geq 0)$ with the relations 
\begin{align}
&[D_{0, m}, D_{0, n}]=0 (m, n\geq 1)\\
&[D_{0, m}, e_n]=- e_{m+n-1} (m\geq 1, n\geq 0) \label{eq:he} \\
&[D_{0, m}, f_{n}]=- f_{m+n-1} (m\geq 1, n\geq 0) \\
&3[e_{m+2}, e_{n+1}]-3[e_{m+1}, e_{n+2}]-[e_{m+3}, e_n]+[e_m, e_{n+3}] \notag\\
&\phantom{1234567891} +\sigma_2([e_{m+1}, e_n]-[e_m, e_{n+1}])+\sigma_3 (e_me_n+e_ne_m)=0  \label{ee}\\
&3[f_{m+2}, f_{n+1}]-3[f_{m+1}, f_{n+2}]-[f_{m+3}, f_n]+[f_m, f_{n+3}] \notag\\
&\phantom{1234567891} +\sigma_2([f_{m+1}, f_n]-[f_m, f_{n+1}])-\sigma_3 (f_mf_n+f_nf_m)=0 \label{ff} \\
& \Sym_{\mathfrak{S}_3}[e_{i_1}, [e_{i_2}, e_{i_3+1}]]=0, \label{eee}  \\
& \Sym_{\mathfrak{S}_3}[f_{i_1}, [f_{i_2}, f_{i_3+1}]]=0. \label{fff} \\
&[e_m, f_n]= h_{m+n} \label{eq:ef=h}\\
&1-\sigma_3 \sum_{n\geq 0} h_n z^{-(n+1)} =
\left\{
     \begin{array}{lr}
       \prod_{i=1}^l (z-z_i) \cdot \psi(z),   & l\geq 0, \\
      \prod_{i=1}^{-l} \frac{1}
{ (z-z_i)} \cdot \psi(z)
, &  l<0.
     \end{array}
   \right.
\label{rel:comm}
\end{align}
where the equality \eqref{rel:comm} means the coefficients of $z^{-i}$ on both sides are equal for each $i\geq 1$. Here $\psi(z):=\exp\left(-\sum_{n\geq 0} D_{0, n+1}  \varphi_n(z)\right)$ and the function $\varphi_n(t)$ is a formal power series in $t$ depending on $\hbar_1, \hbar_2, \hbar_3$. 
It is given by the following formula. 
\begin{align*}
&\exp( \sum_{n\geq 0}(-1)^{n+1} a^{n} \varphi_n(z))
=\frac
{(z+a-\hbar_1)
(z+a-\hbar_2)
(z+a-\hbar_3)}{(z+a+\hbar_1)(z+a+\hbar_2)(z+a+\hbar_3)}, 
\end{align*}
where $a$ is any element in $\C$. 
\end{defn}
When $l=0$, this is the Yangian of $\widehat{{\mathfrak{gl}}(1)}$. We denote it by $\affY$. 
When $l>0$ ($l<0$), we refer to it as the $l$-positively ($-l$-negatively) shifted Yangian in parameters $\vec{z}$. 
The positively shifted Yangian $Y_{l}(\vec{z})$ is first defined in \cite[Section 6]{KN}.
\begin{rmk}
The shifted Yangian $Y_{l}(\vec{z})$ in Definition \ref{l shifted} has a parameter $\vec{z}$, while the shifted Yangian $Y_\mu$ in Definition \ref{def:shiftedYangian} has no such parameters. 

Set $\psi(z)=1-\sigma_3 \sum_{k\geq 0} \psi_{k} z^{-(k+1)}$. When $\vec{z}=0$,  relation \eqref{rel:comm} becomes 
\begin{equation}\label{eq:rel:comm}
1-\sigma_3 \sum_{n\geq 0} h_n z^{-(n+1)} =z^l( 1-\sigma_3 \sum_{k\geq 0} \psi_k z^{-(k+1)}). 
\end{equation}
When $l\geq 0$,  relation \eqref{eq:rel:comm} means 
\[
h_{0}=\psi_{l}, \,\ h_1=\psi_{l+1}, \,\ h_2=\psi_{l+2}, \,\ \cdots, \,\ h_{n}=\psi_{l+n}, \,\ \cdots
\]
When $l< 0$,  relation \eqref{eq:rel:comm} means
\[
h_{0}=0, \,\ h_1=0, \,\ \cdots, h_{-l-1}=0, \,\ h_{-l}=\psi_{0},\,\  h_{-l+1}=\psi_{1}, \,\  \cdots, \,\ h_{n}=\psi_{l+n}, \,\ \cdots
\]
One can define $\psi_{n}, (n\in \mathbb{Z})$ such that $\psi_{n}=0$ for $n<0$. This defines 
$h_{n} (n\in \mathbb{Z})$ with the condition that $h_{n}=\psi_{l+n}$, for all $n\in \Z$. 

So the generators of the shifted Yangian $Y_l(\vec{0})$, can be taken as $e_{n}, f_{n}$ $(n\geq 0)$, $h_{n}$ $(n\in \mathbb{Z})$, such that $h_{n}=\psi_{l+n}$, for $n\geq -l$, quotient by the relation
\[
h_{p}=0,\ \text{for all $p< -l$}\ \text{and}\ \ h_{-l}=\psi_0.  
\] Such presentation has the same form as the one in Definition \ref{def:shiftedYangian}. 
\end{rmk} 

\begin{rmk}
In the above presentation, we have the correspondence compared with the notation in \cite[Section 6]{KN}. 
\[
\hbar_1\mapsto \hbar,\ \hbar_1\mapsto \bold{t},\ \hbar_3\mapsto -\hbar-\bold{t},\ D_{0, n} \mapsto \frac{D_{0, n} }{\hbar},\ e_n \mapsto e_n,\ f_{n}\mapsto f_{n+l}.
\]
We also take the central elements $c_0=c_1=\cdots=0$ in the central extension $\SH^{\bold{c}}$ of $\affY$, 
so that the factor $\frac{(1-(\hbar+\bold{t})x)(1+wtx) }{1-(\hbar+(1-w)\bold{t})x}$ in \cite{KN} does not show up in the formula of $\psi(z)$.  
\end{rmk}
\subsection{Relation to the Costello conjectures} \footnote{We thank one of the referees for pointing out an error in the earlier version. }

\label{Costello Conj}

We explain how Theorem \ref{thm:intro} implies a special case of Costello conjectures  \cite{Costello}.

Let $X_l=\widetilde{\C^2/(\Z/l\Z)}$ be a crepant resolution of $\C^2/(\Z/l\Z)$. 
We have the Hilbert scheme  $X_l^{[k]}$ of $k$ points on $X_l$, which is smooth and quasi-projective. 
We denote the space of quasi-maps from $\PP^1$ to $X_l^{[k]}$ of degree $d$ by 
$\mathcal{M}_{k, d, l}$. 

In \cite{Costello} the author defined the algebra $A_{l, \hbar, \epsilon}$  which is a deformation of the enveloping algebra $U(\text{Diff}_\epsilon(\C)\otimes \mathfrak{gl}(l))$. Here $\text{Diff}_\epsilon(\C)$ denotes the algebra of differential operators on the line $\C$ with $[\partial_z, z]=\epsilon$, and $\hbar$ is the deformation parameter. 

 It is conjectured in \cite[P73]{Costello} that $A_{l, \hbar, \epsilon}$ acts on 
$\bigoplus_{d}H^*_{T}(\mathcal{M}_{k, d, l}, \mathcal{P})$. 
The space  $\mathcal{M}_{k, d, l}$ can be realized as the critical locus 
of a regular  function (``potential") defined similarly to the one in \cite[Section 4.3]{O17}. The sheaf of vanishing cycles of this function is the perverse constructible sheaf $\mathcal{P}$ on $\mathcal{M}_{k, d, l}$ 
\cite[P72]{Costello}. Here $\hbar, \epsilon$ are equivariant parameters corresponding to the 2-dimensional torus $T=\C^*_\hbar\times \C^*_\epsilon$ action.

It is explained in \cite[Theorem 1.6.1]{Costello} that the algebra $A_{l, \hbar, \epsilon}$ is a certain
large $N$ limit of quantized Higgs branch algebra $\mathcal{O}_{\hbar}(\mathcal{M}_{N, l}^{\epsilon})$ of the Jordan quiver gauge theory, i.e.
\[
A_{l, \hbar, \epsilon}=\lim_{N\to \infty} \mathcal{O}_{\hbar}(\mathcal{M}_{N, l}^{\epsilon}). 
\]
Equivalently,  $\mathcal{O}_{\hbar}(\mathcal{M}_{N, l}^{\epsilon})$ is the Coulomb branch algebra of a $\hat{A}_{l-1}$ quiver gauge theory. 

The algebra $A_{l, \hbar, \epsilon}$ is expected to be the deformed double current algebra $U_{\hbar}(\text{Diff}_\epsilon(\C)\otimes \mathfrak{gl}(l))$, which is the imaginary $1$-shifted affine Yangian of $\mathfrak{gl}(l)$. 

We recall the spherical cyclotomic rational Cherednik algebra (see \cite[Section 5]{KN}), which is the Coulomb branch algebra of the Jordan quiver gauge theory. 

Let $\Gamma_N:= \mathfrak{S}_N\ltimes (\Z/l\Z)^N$ be the wreath product of the symmetric group and the cyclic group of order $l$. 
Denote a fixed generator of the $i$-th factor of $ (\Z/l\Z)^N$ by $\alpha_i$. The group $\Gamma_N$ acts on $\C\langle \xi_1, \cdots, \xi_{N}, \eta_1, \cdots, \eta_{N}\rangle$ by 
\[
\alpha_i(\xi_i)=\epsilon \xi_i, \,\
\alpha_i(\xi_j)=\xi_j, 
\alpha_i(\eta_i)=\epsilon^{-1}\eta_i, 
\alpha_i(\eta_j)=\eta_j\,\  (i\neq j)
\]with the obvious $\mathfrak{S}_N$-action. Here $\epsilon$ denotes a primitive $l$-th root of unity. 
\begin{defn}
The \textit{cyclotomic rational Cherednik algebra} $\bold{H}^{cyc}_{N, l}$ for $\mathfrak{gl}(N)$ is the quotient of the algebra $\C[\hbar_1, \hbar_2, c_1, \cdots, c_{l-1}] \langle \xi_1, \cdots, \xi_{N}, \eta_1, \cdots, \eta_{N}\rangle \rtimes \Gamma_N$ by the relations
\begin{gather}
[\xi_i, \xi_j]=0=[\eta_i, \eta_j], \,\ (i, j=1, \cdots, N)\\
[\eta_i, \xi_j]=\left\{
                \begin{array}{ll}
                  -\hbar_1+\hbar_2 \sum_{k\neq i} \sum_{m=1}^{l-1} s_{ik} \alpha_i^m \alpha_k^{-m} +\sum_{m=1}^{l-1} c_m \alpha_i^m& \text{if $i=j$}\\              
                 -\hbar_2 \sum_{m=0}^{l-1} s_{ij} \epsilon^m \alpha_i \alpha_j^{-m}& \text{if $i\neq j$}
                \end{array}
              \right.
\end{gather}
Let $e_{\Gamma_N}$ be the idempotent for the group $\Gamma_N$. The spherical cyclotomic rational  Cherednik algebra is defined as
\[
\bold{SH}^{cyc}_{N, l}=e_{\Gamma_N} \bold{H}^{cyc}_{N, l} e_{\Gamma_N}. 
\]
\end{defn}

\begin{thm}\cite[Theorem 1.1, Theorem 6.14, Theorem 1.5]{KN}
\begin{enumerate}
\item
Let $\mathcal{A}_{\hbar_1}$ be the quantized Coulomb branch algebra of the gauge theory $(G, \bold{N})=(\GL(N), \mathfrak{gl}(N)\oplus (\C^N)^{\oplus l})$. 
If $l>0$ then $\mathcal{A}_{\hbar_1}$ is isomorphic to $\bold{SH}^{cyc}_{N, l}$ with an  explicitly given  correspondence between parameters in both algebras in {\it loc.cit.}.
\Omit{
The parameters $\bold{t}, \hbar$ of the quantized Coulomb branch algebra and the Cherednik algebra are the same, the others are matched by 
\[
z_k=-\frac{1}{l} \left(
(l-k) \hbar+\sum_{m=1}^{l-1} \frac{c_m (1-\epsilon^{mk})}{1-\epsilon^m}
\right)
\]
when $l>0$. }
\item
Let $\mathcal{A}_{\hbar_1}$ be the quantized Coulomb branch algebra for $\dim(V)=N$, $\dim(W)=l$. Then, we have a surjective homomorphism of algebras 
$
\Psi: Y_{l}(\vec{z})\to \mathcal{A}_{\hbar_1} $. 
\Omit{
given by 
\begin{align*}
&D_{0, m} \mapsto \sum_{i=1}^N \bar{B}_m(w_i)-\bar{B}_m(-(i-1)\bold{t}), \,\  (m\geq 1)\\
&e_n\mapsto E_1[(w+\hbar)^n], \,\, 
f_{n+l} \mapsto F_1^{(l)} [(w+\hbar)^n], \,\ (n\geq 0). 
\end{align*}
where $F_1^{(l)} [(w+\hbar)^n]=F_1^{(0)}[(w+\hbar)^n \prod_{k=1}^l (w-z_k)]=\sum_{i=0}^l (-1)^{l-i} e_{l-i}(\vec{z}-\hbar) F_1^{(0)}[(w+\hbar)^{i+n}]$. }
\Omit{\item
There is a faithful embedding of $\bold{SH}^{cyc}_{N, l}$ to the ring $\tilde{\mathcal{A}}_{\hbar_1}$ of localized difference operators such that 
\begin{gather}
e_{\Gamma_N} \left(
\sum_{i=1}^N ( l^{-1} \xi_i\eta_i +\bold{t} \sum_{j<i} s_{ji})^n
\right)e_{\Gamma_N}=\sum_{i=1}^N \omega_i^n, \,\ (n>0),  \\
e_{\Gamma_N} \left( \sum_{i=1}^N \xi_i^l \right)e_{\Gamma_N} =E_1[1], \,\ 
e_{\Gamma_N} \left( \sum_{i=1}^N (l^{-1}\eta_i)^l \right)e_{\Gamma_N} =F_1[1]. 
\end{gather}}
\end{enumerate}
\end{thm}
It follows from the above theorem that  we have an epimorphism
\[
Y_{l}(\vec{z})\to  \bold{SH}^{cyc}_{N, l}, \text{for all $N$, such that $e_0\mapsto e_{\Gamma_N} \left( \sum_{i=1}^N \xi_i^l \right)e_{\Gamma_N}$. } 
\] 
The relation of the shifted Yangian $Y_{l}(\vec{z})$ and the spherical cyclotomic rational Cherednik algebra $\bold{SH}^{cyc}_{N, l}$ is
\begin{equation}\label{shifted Yangian and cyc RCA}
Y_{l}(\vec{z})=\lim_{N\to \infty} \bold{SH}^{cyc}_{N, l}. 
\end{equation}
Here according to \cite[\S~13.4.2.~Remark]{Costello} the limit is expected to be made precise using Deligne categories,  similar to the tensor categories extending those of representations of the symmetric groups $S_n$ to complex values of $n$.  Then the index $N$ becomes a central element in the limit algebra. See also \cite{EKR} for an analogue of this limit process via Deligne categories. It is expected that the generators and relations in \cite{GY} give another presentation of this algebra.

In this paper, we focus on the case where $l=1$, thus $X_l=\C^2$, and $X_l^{[k]}$ is the Hilbert scheme $\Hilb^k(\C^2)$. 
The above-mentioned algebra $A_{1, \hbar, \epsilon}$ \cite{KN} is then isomorphic to  $Y_{1}(z_1)$. 

Let $X_{1, 1}$ be the resolved conifold. The PT moduli space of $X_{1, 1}$ space is identified with the space of $\calO(-1)$-twisted  quasi-maps from $\PP^1$ to $\Hilb^k(\C^2)$   \cite[\S 4.3.19 Exercise 4.3.22]{O17}. Recall that the PT moduli space parameterizes the complexes $\mathcal{O}_{X_{1, 1}} \to \mathcal{F}$ of sheaves on $X_{1, 1}$, with the condition that the sheaf $\mathcal{F}$ is pure 1-dimensional (i.e., has no 0-dimensional subsheaves),  and $\text{Cokernel}(s)$ is 0-dimensional. 
The sheaf $\mathcal{F}$ is shown to be an extension of the structure sheaf of a curve in $X_{1,1}$ by a zero-dimensional sheaf supported on this curve \cite[\S~1.3]{PT07}. 
Note that $X_{1,1}$ is the total space of a rank-2 vector bundle on $\PP^1$, the projection to $\PP^1$ gives a map from the curve supporting $\mathcal{F}$  to $\PP^1$. Pushing-forward via the projection $\pi:X_{1,1}\to \PP^1$, the sheaf $\mathcal{F}$ gives a vector bundle on $\bbP^1$ the rank of which is equal to the degree of the map from this curve to $\PP^1$. The action of $\pi_*\calO_{X_{1,1}}$ on the vector bundle $\pi_*\mathcal{F}$ is equivalent to a pair of commuting $\mathcal{O}(-1)$-valued Higgs fields. The section of $\mathcal{F}$ is equivalent to a section of the vector bundle $\pi_*\mathcal{F}$, which generates the bundles under the action of $\pi_*\calO_{X_{1,1}}$ hence of the two Higgs fields.  
This  is the same as a quasi-map from $\PP^1$ to the Hilbert scheme of points on $\C^2$.

By Theorem \ref{thm:intro}, we have an action of the shifted Yangian $Y_{1}(z_1)$ on the cohomology of the moduli space of perverse coherent systems  on $X_{1, 1}$ with stability parameter on the PT-side of the imaginary root hyperplane. 
This can be viewed as an action of $A_{l, \hbar, \epsilon}$ on $\bigoplus_{d}H^*_{T}(\mathcal{M}_{k, d, l}, \mathcal{P})$, when $l=1$.

\section{Hilbert scheme of $\C^3$ and shifted Yangians}
\label{sec:Hilb}
In this section we prove Theorem~\ref{thm:intro}(1). 
\subsection{The Hilbert scheme of points on $\C^3$}

Consider the following framed quiver $\tilde Q_3$ with potential
\[
\begin{tikzpicture}
  \draw ($(0,0)$) circle (.08);
  \node at (0, -1.5) {$\square$};
    \node at (0.4, -1.5) {$1$};
     \node at (0.3, 0) {$n$};
       \node at (1.2, 0) {$B_3$};   
       \node at (0, 1.2) {$B_2$};
       \node at (-1.2, 0) {$B_1$};  
       \node at (0.2, -0.8) {$I$};  
       
  \draw [<-] (0, -0.1) -- (0, -1.4);
   \draw[->,shorten <=7pt, shorten >=7pt] ($(0,0)$)
   .. controls +(40:1.5) and +(-40:1.5) .. ($(0,0)$);
     \draw[->,shorten <=7pt, shorten >=7pt] ($(0,0)$)
   .. controls +(90+40:1.5) and +(90-40:1.5) .. ($(0,0)$);
 \draw[->,shorten <=7pt, shorten >=7pt] ($(0,0)$)
   .. controls +(180+40:1.5) and +(180-40:1.5) .. ($(0,0)$);
   \node at (0, -2.3){ The potential is $W_3=\tilde W_3=B_3([B_1, B_2])$};
\end{tikzpicture}
\]

A representation of the quiver $Q_3$ of dimension $n$ is  stable framed if the following additional property is satisfied
\[
\C\langle B_1, B_2, B_3\rangle I(\C)=\C^n.
\] 
Here $\C\langle B_1, B_2, B_3\rangle$ is the algebra of non-commutative polynomials in the variables
$B_1, B_2, B_3$. 

The set of stable framed representations of dimension $n$ is denoted by $\calM(n, 1)^{st}$. It  consists of triples of $n\times n$ matrices $B_1, B_2, B_3\in \End(\C^n)$, together with a cyclic vector $v=I(1)\in \C^n$. 
The group $\GL_n$ acts by conjugation. Cyclicity of $v$ means that it generates $\C^n$ under the action of $B_i's$. 

Note that  the critical locus $\Crit(\tilde W)$ of the potential $\tr \tilde W$ in $\calM(n, 1)^{st}$ consists of  triples of commuting matrices $(B_1,B_2,B_3)$  satisfying the property that $\Image (I)$ is a cyclic vector of $\C^n$ under the three matrices. Therefore, the 	quotient stack $\Crit(\tilde W)/\GL_n
$ is a scheme isomorphic to $\Hilb^n(\C^3)$. Here $\Hilb^n(\C^3)$ is the Hilbert scheme of $n$-points on $\C^3$ 
defined by
\[
\Hilb^n(\C^3)=\{J  \subset  \C[x_1, x_2, x_3]\mid \text{$J$ is an ideal, and} \,\ \dim(\C[x_1, x_2, x_3]/J)=n\}. 
\]
The isomorphism $\Crit(\tilde W)/\GL_n
\cong \Hilb^n(\C^3)$ is given such as follows. 
To an ideal $J \subset  \C[x_1, x_2, x_3]$ let us choose an isomorphism $\C[x_1, x_2, x_3]/J\cong \C^n$. Then the linear maps $B_1, B_2, B_3\in \End(\C^n)$ are given by multiplications by $x_1, x_2, x_3$ mod $J$. Furthermore $I\in \Hom(\C, \C^n)$ is $I(1)=1$ mod $J$. 
Different choices of basis of $\C[x_1, x_2, x_3]/J$ give isomorphic representations.
Conversely, if we have $(B_1, B_2, B_3, I)\in \Crit(\tilde W)/\GL_n$, then, the ideal $J$ is defined as the  kernel of the following map
\[
\phi: \C[x_1, x_2, x_3] \to \C^n, f\mapsto f(B_1, B_2, B_3) I(1). 
\] Note that $\phi$ is surjective by the framed stability condition. Hence $\dim(\C[x_1, x_2, x_3]/J)=n$. 

There is an action of $(\C^*)^3$ on $\Hilb^n(\C^3)$ induced by rescaling of the coordinates $x_1, x_2, x_3$. Fixed points $(\Hilb^n(\C^3))^{(\C^*)^3}$ of this action are in one-to-one correspondence with  
3-dimensional partitions of $n$.

By \cite[Theorem 5.2.1]{RSYZ}, we have an action of $\mathcal{H}^{(Q_3, W_3)}$ on 
\[
\bigoplus_{n\geq 0} H^\ast_{c,\GL_n\times T\times \C^*_{fr}}(\mathcal{M}(n, 1)^{st}, \varphi_{\tilde W_3}\C)^\vee
\cong \bigoplus_{n\geq 0}  H_{c, T\times \C^*_{fr}}^\ast(\Hilb^n(\C^3), \varphi_{\tilde W_3}\C)^\vee
\]
given by natural correspondences. It also induces an action of the spherical subalgebra $\mathcal{SH}^{(Q_3, W_3)}$.

\subsection{Action of the shifted Yangian}
\label{sub:Action of the shifted Yangian}
Denote by $R_T$ the cohomology ring $H_{c, T}^\ast(\pt)^\vee$, and  $K_T$ its  field of fractions. 
In this section, we will construct an action of the shifted Yangian $Y_{-1}(z_1)$ on 
\[
\bigoplus_{n\geq 0} H_{c, T\times \C^*_{fr}}^\ast(\Hilb^n(\C^3), \varphi_{\tilde W_3}\C)^\vee\otimes_{R_T} K_T. 
\]
In general the  scheme $\Hilb^n(\C^3)$ is singular which makes considerations more complicated than in the case $n=2$.
However, we have the following isomorphism of cohomology groups
\[
H_{c, T\times \C^*_{fr}}^\ast(\Hilb^n(\C^3), \varphi_{\tilde W_3}\C)^\vee \cong H_{c,\GL_n\times T\times \C^*_{fr}}^\ast(\mathcal{M}(n, 1)^{st}, \varphi_{\tilde W_3}\C)^\vee, 
\]
where $\mathcal{M}(n, 1)^{st}$ and $\mathcal{M}(n, 1)^{st}/\GL_n$ are smooth. 

Let $V$ denote the standard  coordinate vector space of dimension $n+1$. Fix $\xi\subset V$, a one dimensional subspace and let $V_2:=V/\xi$ be the quotient. 
Consider the following correspondence 
\[
\xymatrix@R=1em@C=0.5em{
	&\mathcal{M}(n, n+1, 1)^{st}\ar[ld]_{p}\ar[rd]^{q}&\\
	\mathcal{M}(n, 1)^{st}&&\mathcal{M}(n+1, 1)^{st}
}\]
where 
\[
\mathcal{M}(n, n+1, 1)^{st}=\{(B_1, B_2, B_3, I)\in  \mathcal{M}(n+1, 1)^{st} \mid B_i(\xi)\subset \xi\}. 
\]
The parabolic subgroup $P:=\{x\in \GL_{n+1}=\GL(V)\mid x(\xi)\subset \xi\}$ acts on $\mathcal{M}(n, n+1, 1)^{st}$. 
The map $p$ is given by $(B_1, B_2, B_3, I)\mapsto (B_1, B_2, B_3, I)$ mod $\xi$, where $I$ mod $\xi$ is the composition of 
$I: \C\to V$ with the projection $V \surj V_2=V/\xi$. The map $q$ is the natural inclusion. 

It induces a correspondence
\[
\xymatrix@R=1em@C=0.5em{
	&\mathcal{M}(n, n+1, 1)^{st}/P \ar[ld]_{p}\ar[rd]^{q}&\\
	\mathcal{M}(n, 1)^{st}/\GL_n& &\mathcal{M}(n+1, 1)^{st}/\GL_{n+1}
}\]
The subspace $\xi\subset V$ gives a tautological line bundle on the correspondence $\mathcal{M}(n, n+1, 1)^{st}/P$, which will be denoted by $L$. 
Recall that the Levi subgroup of $P$ is $\C^*\times\GL_n$ where 
$\xi$ above is the standard weight 1 representation of the $\C^*$-factor. Hence, $c_1(L)$ coincide with the equivariant variable of the $\C^*$-factor.
For any rational function $g(x)$ in one variable with coefficients in $H^*_T(\pt)$ we define the class $g(c_1(L))$ in a localization of $H_{c,P\times T}^\ast(\mathcal{M}(n, n+1,1)^{st}, \varphi_{\tilde W_3}\C)^\vee$. Here the localization is taken with respect to $H_{c,T\times\C^*}(\pt)\cong H_{c,T}(\pt)[x]$.

Define the actions of the operators 
\begin{align}\label{eqn:EFdef}
E(g): H&_{c,\GL_n\times T\times \C^*_{fr}}^\ast(\mathcal{M}(n, 1)^{st}, \varphi_{\tilde W_3}\C)^\vee\otimes_{R_T} K_T  \\& \to H_{c,\GL_{n+1}\times T\times \C^*_{fr}}^\ast(\mathcal{M}(n+1, 1)^{st}, \varphi_{\tilde W_3}\C)^\vee\otimes_{R_T} K_T, \\
F(g): H&_{c,\GL_{n+1}\times T\times \C^*_{fr}}^\ast(\mathcal{M}(n+1, 1)^{st}, \varphi_{\tilde W_3}\C)^\vee\otimes_{R_T} K_T  \\& \to H_{c,\GL_n\times T\times \C^*_{fr}}^\ast(\mathcal{M}(n, 1)^{st}, \varphi_{\tilde W_3}\C)^\vee\otimes_{R_T} K_T
\end{align}
 by the following convolutions: 
\[
E(g):=q_*(g(c_1(L))\cup p^*), \,\  F(g):=p_*(g(c_1(L))\cup q^*). 
\]
Note that the fixed point loci with respect to the $T$-actions are finite for all three spaces above. Hence the push-forward  $p_*$ is well-defined after passing to the localization despite of the fact that $p$ is not a proper map.

\begin{rmk}
Note that $T$ is a Calabi-Yau subtorus of $(\C^*)^3$, i.e. it preserves the standard holomorphic volume form. The fixed-point subset of the critical locus for the action of the  torus $(\C^*)^3$ coincides with the one for the subtorus  $T$. This claim can be found e.g. in \cite[Page 25]{MNOP}. A proof   in the case of $\C^4$ that topologically these two fixed-point loci agree can be found in \cite[Lemma 3.1]{CK}. Scheme theoretical statement is shown in \cite[Lemma 3.6]{CK}. Both arguments rely in an essential way on the finite length condition of the defining ideal. There are few examples in which the corresponding fact is shown for ideals of curves \cite[Lemma 2.1, 2.2]{CK2}
or PT stable pairs \cite[Prop. 2.6]{CK2} are shown.
\end{rmk}

Let $\lambda$ be a 3-dimensional Young diagram with $n$ boxes, which we denote as a partition $\lambda\vdash n$. 
The dual of the compactly supported cohomology group $$H_{c,\GL_n\times T\times \C^*_{fr}}^\ast(\mathcal{M}(n, 1)^{st}, \varphi_{\tilde W_3}\C)^\vee\otimes_{R_T} K_T$$ has a basis given by $\{\lambda\mid \lambda \vdash n\}$. 
Let $(\lambda+\blacksquare)\vdash (n+1)$ be the Young diagram obtained by adding a box $\blacksquare$ to $\lambda$. 
We use the notation $\langle\lambda|(E(g))| \lambda+\blacksquare\rangle$ for the coefficient of $\lambda+\blacksquare$ in the expansion of $E(g)(\lambda)$. Similar convention is used for $\langle \lambda+\blacksquare|(F(g))| \lambda\rangle$.

We denote by $\square_{i, j, k}$ the box in the $3$-dimensional Young diagram $\lambda$ which has coordinates $(i, j, k)$. Here we follow the French convention 
in writing the coordinates of boxes. For example, all  boxes has non-negative coordinates and the corner box has coordinates $(0,0,0)$.
Recall that $(\C^*)^3$ naturally acts on $\C^3$ and $H^*_{c,(\C^*)^3}(\pt)=\C[\hbar_1,\hbar_2,\hbar_3]$ with $\hbar_i$ being the equivariant variable of the $i$th $\C^*$-factor of $(\C^*)^3$. We also have the 2-dimensional torus $T\subseteq(\C^*)^3$ preserving the canonical line bundle, and $H^*_{c,T}(\pt)=\C[\hbar_1,\hbar_2,\hbar_3]/(\hbar_1+\hbar_2+\hbar_3=0)$.
We also have an $\C^*_{fr}$ acting on the framing, and let $\chi$ be the $\C^*_{fr}$-equivariant variable. 
For any box $\square_{i, j, k}$, we consider the following element $x_{\square_{i, j, k}}:=\chi+i\hbar_1+j\hbar_2+k \hbar_3\in H^*_{c,T\times\C^*_{fr}}(\pt)$. In the proof, we show this element is equal to the equivariant Chern root of some vector bundle.

\begin{prp}\label{prop:ef formula}
	The matrix coefficients of the operators $E(g), F(g)$ in the  basis of fixed points are as follows:
	\begin{align*}
	\langle\lambda|(E(g))| \lambda+\blacksquare\rangle=
	&  \Res_{z=x_{\blacksquare}}g(z)
	\frac{1}{z-x_{\square_{000}}}
	\prod_{\square\in \lambda}
	\frac
	{z-x_{\square}}{
		(z-x_{\square}-\hbar_1)(z-x_{\square}-\hbar_2)(z-x_{\square}-\hbar_3)
	}
	\\
	\langle \lambda+\blacksquare|(F(g))| \lambda\rangle=&g(z)  \prod_{\square\in \lambda} \frac{(z-x_{\square}+\hbar_1)(z-x_{\square}+\hbar_2)(z-x_{\square}+\hbar_3)} {z-x_{\square}}|_{z=x_{\blacksquare}}	\end{align*}
Here for a box $\square$ or $\blacksquare$ in a 3-dimensional Young diagram, $x_\square$ or $x_\blacksquare$ stands for the $T$-weight of the corresponding box, which is an element in $H^*_{T\times \C^*_{fr}}(\pt)$. 
\end{prp}

{\it Proof.}
Both formulas can be  derived from \eqref{eqn:JefferyKirwan} and \eqref{eqn:convolution_localization}. In order to use \eqref{eqn:JefferyKirwan} we recall that the correspondence $G\times_P\mathcal{M}(n, n+1, 1)^{st}=\GL_{n+1}\times_P\mathcal{M}(n, n+1, 1)^{st}$ has a free $G$-action. The Levi subgroup of $P$ is $\C^*\times \GL(n)$. We denote the equivariant variable corresponding to the $\C^*$-action by $-\tilde{z}$. Then at  the point $(\lambda, \lambda+\blacksquare) \in G\times_P\mathcal{M}(n, n+1, 1)^{st}$ the tautological bundle $L$ becomes the standard 1-dimensional representation of $\C^*$, and hence the $T\times \C^*\times \C^*_{fr}$ equivariant $c_{1}(L)|_{(\lambda, \lambda+\blacksquare)}=-\tilde{z}+x_{\blacksquare}$. 
Let $z=-\tilde{z}+x_{\blacksquare}$. 
Now 
\eqref{eqn:JefferyKirwan} (residue at $\tilde{z}=0$ amounts to residue at $z=x_{\blacksquare}$) gives
\begin{align*}
&\langle\lambda|(E(g))| \lambda+\blacksquare\rangle
=\Res_{z=x_{\blacksquare}}g(z) \frac{e(T_{\lambda} (\mathcal{M}(n,1)^{st}/\GL_n )} {e(T_{\lambda, \lambda+\blacksquare} (\mathcal{M}(n, n+1, 1)^{st}/P)}
\end{align*}
Here $e$ denotes  the $T\times \C^*\times \C^*_{fr}$-equivariant Euler class, where $\C^*$ is the above-mentioned factor of  the Levi subgroup of $P$.

For the lowering operator $F(g)$, it is easier to use the formula \eqref{eqn:convolution_localization} and only consider the $T\times \C^*_{fr}$-equivariance. Note that $E(g)$ and $F(g)$ have different forms: $E(g)$ is given by a residue formula, and $F(g)$ is given by a restriction. This is for convenience of verifying their commutation relations \eqref{eq:ef=h}, see \S\ref{subsub:ef=h}. 

Now  \eqref{eqn:convolution_localization} gives
\begin{align*}
&\langle \lambda+\blacksquare|(F(g))| \lambda\rangle
=g(c_1(L))|_{(\lambda, \lambda+\blacksquare)} \frac{e(T_{\lambda+\blacksquare} (\mathcal{M}(n+1,1)^{st}/\GL_{n+1} )}{e(T_{\lambda, \lambda+\blacksquare} (\mathcal{M}(n, n+1, 1)^{st}/P)}
\end{align*}
Here $e$ stands for the $T\times \C^*_{fr}$-equivariant Euler class.

Next  we calculate all the tangent spaces at the fixed points. 
 We start by calculating the tangent space for $\mathcal{M}(n+1,1)^{st}/\GL_{n+1} $.
Let $V$ be the tautological bundle of rank $n+1$. 
For simplicity, for any vector bundle $X$, we write the dual bundle as $X^\vee$. Hence, for two vector bundes $X,Y$ we have $
\Hom(X, Y)=X^\vee\otimes Y
$.

The tangent complex $T(\mathcal{M}(n+1,1)^{st}/\GL_{n+1} )$ is given by 
\[
\End(V)\to \Big(
\C^3 \otimes \End(V)+\Hom(W, V)\Big). 
\]
Its  class in the Grothendieck group of $\mathcal{M}(n+1,1)^{st}/\GL_{n+1} $ is then given by 
\begin{align*}
T(\mathcal{M}(n+1,1)^{st}/\GL_{n+1}) 
= \C^3\otimes V^\vee\otimes V+W^\vee\otimes V-\End(V). 
\end{align*}

The tangent complex of the correspondence $\mathcal{M}(n, n+1, 1)^{st}/P$ has the following class in the Grothendieck group, where  we do not include the terms coming from obstruction theory. 
\[
T(\mathcal{M}(n, n+1, 1)^{st}/P)=\C^3\otimes P+W^\vee\otimes V-P.
\] 

Now formula \eqref{eqn:JefferyKirwan}  gives the following
\begin{align*}
\langle\lambda|(E(g))| \lambda+\blacksquare\rangle
=&\Res_{z=x_{\blacksquare}}g(c_{1}(L))|_{(\lambda, \lambda+\blacksquare)} \frac{e(T_{\lambda} (\mathcal{M}(n,1)^{st}/\GL_n )} {e(T_{\lambda, \lambda+\blacksquare} (\mathcal{M}(n, n+1, 1)^{st}/P)} |_{(\lambda, \lambda+\blacksquare)}   \\
=&\Res_{z=x_{\blacksquare}}g(z) \frac{e(\C^3\otimes V_2^\vee\otimes V_2+W^\vee\otimes V_2-\End(V_2))}{e(\C^3\otimes P+W^\vee\otimes V-P)}\\
=&\Res_{z=x_{\blacksquare}}g(z)\frac{e(V_2^\vee\otimes \xi)}{e(\C^3 \otimes V_2^\vee\otimes \xi)} \frac{1}{e(W^\vee\otimes \xi)}
\end{align*}
Here as before $z$ is the equivariant Chern root of $\xi$ and recall that $V_2=V/\xi$. We take the basis of  $V_2$ to be $\{{\square} \mid \square\in \lambda\}$. Then the equivariant Chern root of the line spanned by $\square$ is the $T\times \C^*\times \C^*_{fr}$-weight of $\square$, which is  $x_{\square}$ in our notation. 
The three loops $B_1, B_2, B_3$ has weights $\hbar_1, \hbar_2, \hbar_3$ respectively. 
Let $\chi$ be the Chern roots of $W$. 
Thus, we have
\begin{align*}
\langle\lambda|(E(g))| \lambda+\blacksquare\rangle=&\Res_{z=x_{\blacksquare}} \big(g(z) 
\prod_{\square\in \lambda} \frac{z-x_{\square}}{(z-x_{\square}-\hbar_1)(z-x_{\square}-\hbar_2)(z-x_{\square}-\hbar_3)
} \frac{1}{z-\chi}\big)
\end{align*}

Similarly, by  \eqref{eqn:convolution_localization} we have
\begin{align*}
\langle \lambda+\blacksquare|(F(g))| \lambda\rangle
=&g(c_1(L))|_{(\lambda, \lambda+\blacksquare)}  \frac{e(T_{\lambda+\blacksquare} (\mathcal{M}(n+1,1)^{st}/\GL_{n+1} )} {e(T_{\lambda, \lambda+\blacksquare} (\mathcal{M}(n, n+1, 1)^{st}/P)}|_{(\lambda, \lambda+\blacksquare)}   \\
=&g(c_{1}(L))|_{(\lambda, \lambda+\blacksquare)}  \frac{e(\mathbb{C}^3\otimes V^\vee\otimes V+W^\vee\otimes V-\End(V))}{e(\mathbb{C}^3\otimes P+W^\vee\otimes V-P)}\\
=&g(c_{1}(L))|_{(\lambda, \lambda+\blacksquare)} \frac{e(\mathbb{C}^3\otimes \xi^\vee\otimes V_2)}{e(\xi^\vee\otimes V_2)}\\
=&g(x_{\blacksquare})
	\prod_{\square\in \lambda} \frac
	{(x_{\blacksquare}-x_{\square}+\hbar_1) (x_\blacksquare-x_{\square}+\hbar_2) (x_\blacksquare-x_{\square}+\hbar_3)}{x_\blacksquare-x_{\square}}
\end{align*}
Recall the line bundle $L$ is defined using the subspace $\xi\subset V$, with the quotient $V/\xi=V_2$.  Therefore, the $T\times \C^*_{fr}$-weight of $\xi$ at the fixed point ${(\lambda, \lambda+\blacksquare)}$ is given by $x_{\blacksquare}$.

This completes the proof. 
$\blacksquare$

Now in the definition \eqref{eqn:EFdef} above we take the function $g(x)$ to be $x^i$. For $e_i\in Y_{-1}(z_1)$, we define the action of $e_i$ on $H_{c, T\times \C^*_{fr}}^\ast(\Hilb^*(\C^3), \varphi_{\tilde W_3}\C)^\vee\otimes_{R_T}K_T$ by the operator $E(x^i)$ above. Similarly, for $f_i\in Y^-_{-1}(\vec{z})$, we define the action of $f_i$ by the operator $F(x^i)$ above. The parameter $z_1$ acts by $\chi$.  

We define the action of $\psi(z)$ on $H_{c, T\times \C^*_{fr}}^\ast(\Hilb^*(\C^3), \varphi_{\tilde W_3}\C)^\vee\otimes_{R_T}K_T$
via the Chern polynomial of  the following class in the Grothendieck group \[(q_1q_2)V+(q_1q_3)V + (q_2q_3)V-
q_1V - q_2V -q_3V.\] Here $q_1, q_2,q_3$ are the three coordinate lines of $\C^3$, considered as $T$-representations. We refer to \cite[\S~4.1]{RSYZ} for the definition of Chern polynomial.  
\begin{prp}
Let $\lambda \vdash n$ be a 3d partition. The eigenvalue of $\psi(z)$ on $[\lambda]$ equals
\[
\prod_{\square\in \lambda}
\frac{z-x_{\square}+\hbar_1}{z-x_{\square}-\hbar_1}
\frac{z-x_{\square}+\hbar_2}{z-x_{\square}-\hbar_2}
\frac{z-x_{\square}+\hbar_3}{z-x_{\square}-\hbar_3}.
\]
\end{prp}
\begin{defn}\label{addible and removable}
For a 3d Young diagram $\lambda$, an  addible box is a box in another 3d Young diagram $\lambda'$, such that $\lambda'$ is obtained by adding this box to $\lambda$. 
A removable box in $\lambda$ is a box such that $\lambda$ with this box removed is still a 3d Young diagram. 
\end{defn}

\begin{lmm}\label{lem:poles}
Assume $\hbar_1+\hbar_2+\hbar_3=0$.  Then
 $h(z)|_{\lambda}=\frac{1}{z-x_{\square_{000}}}\psi(z)|_{\lambda}$  is multiplication by a rational function, the poles of which are at the addible and removable boxes of $\lambda$. 
\end{lmm}
{\it Proof.}
Follows from a careful examination of the cancellations. 
$\blacksquare$

\subsection{Checking the relations of the shifted Yangian}\label{subsec:C3Yangian}

\subsubsection{Checking the relations \eqref{ee} \eqref{ff} \eqref{eee} \eqref{fff}}
This follows from the epimorphism 
\begin{align*}
& Y_{-1}^+(z_1) \to \mathcal{SH}^{(Q_3, W_3)}, \,\ e_{n} \mapsto  x^n
, \\
&Y_{-1}^-(z_1)\mapsto  \mathcal{SH}^{(Q_3, W_3)} , \,\  (-1)^{n-1} f_{n-1}\mapsto  x^n, 
 \end{align*}
 and the compatibility of the COHA action and the $Y_{-1}^+(z_1), Y_{-1}^-(z_1)$ actions on the Hilbert scheme. 
 Recall that $x$ is the algebraic generator of $\mathcal{SH}_1^{(Q_3, W_3)}= H_{T}(\pt)[x]$. 

\subsubsection{Checking the relation \eqref{eq:he}}
This is exactly the same as the affine Yangian case. 

Let $e(y)$ be the generating series $\sum_{n=0}^{\infty} e_n y^{n+1}$, where $y$ is a formal variable. 
The relation  \eqref{eq:he} is equivalent to the following relation (see \cite{KN})
\begin{equation}\label{eq:he2}
\psi(z) e(y)= \left(
e(y) \psi(z) \frac{(z-y^{-1}-\hbar_1)(z-y^{-1}-\hbar_2)(z-y^{-1}-\hbar_3)}
{(z-y^{-1}+\hbar_1)(z-y^{-1}+\hbar_2)(z-y^{-1}+\hbar_3)}\right)_+
\end{equation}
where $( )_{+}$ denotes the part with positive powers in $y$. Here $\frac{(z-y^{-1}-\hbar_1)(z-y^{-1}-\hbar_2)(z-y^{-1}-\hbar_3)}
{(z-y^{-1}+\hbar_1)(z-y^{-1}+\hbar_2)(z-y^{-1}+\hbar_3)}$ is regarded as an element in $\mathbb{C}[\hbar_1, \hbar_2][[z^{-1}]]$.

The action of the series $e(y)$ on $H_{c, T\times \C^*_{fr}}^\ast(\Hilb^*(\C^3), \varphi_{\tilde W_3}\C)^\vee\otimes_{R_T}K_T$ is given by the operator $E(\frac{1}{y^{-1}-x})$.
By Proposition \ref{prop:ef formula} the matrix elements of $e(y)$ are given by
	\begin{align*}
	\langle\lambda|(e(y))| \lambda+\blacksquare\rangle=
	&
	\frac{y}{1-x_{\blacksquare} y}
	\Res_{z=x_{\blacksquare}}
	\frac{1}{z-x_{\square_{000}}}
	\prod_{\square\in \lambda}
	\frac
	{z-x_{\square}}{
		(z-x_{\square}-\hbar_1)(z-x_{\square}-\hbar_2)(z-x_{\square}-\hbar_3)
	}.
	\end{align*}	
Let $A:=\Res_{z=x_{\blacksquare}}
	\frac{1}{z-x_{\square_{000}}}
	\prod_{\square\in \lambda}
	\frac
	{z-x_{\square}}{
		(z-x_{\square}-\hbar_1)(z-x_{\square}-\hbar_2)(z-x_{\square}-\hbar_3)
	}
$.
Applying both sides of \eqref{eq:he2} to the partition $\lambda$, and looking at the coefficient of $\lambda+\blacksquare$, we obtain
\begin{align*}
\langle\lambda|(\psi(z) e(y))| \lambda+\blacksquare\rangle=\psi(z)|_{\lambda+\blacksquare}  \langle\lambda|(e(y))| \lambda+\blacksquare\rangle
=\psi(z)|_{\lambda+\blacksquare}   \frac{y}{1-x_{\blacksquare} y}A.
\end{align*}	
We have following inductive formula of $\psi(z)$
\[
\psi(z)|_{\lambda+\blacksquare}=\psi(z)|_{\lambda}\cdot 
\frac{z-x_{\blacksquare}-\hbar_1}{z-x_{\blacksquare}+\hbar_1}
\frac{z-x_{\blacksquare}-\hbar_2}{z-x_{\blacksquare}+\hbar_2}
\frac{z-x_{\blacksquare}-\hbar_3}{z-x_{\blacksquare}+\hbar_3}. 
\]
Therefore, the left hand side of \eqref{eq:he2}  becomes 
\begin{align*}
&\langle\lambda|(\psi(z) e(y))| \lambda+\blacksquare\rangle\\
=&\psi(z)|_{\lambda}\cdot 
\frac{z-x_{\blacksquare}-\hbar_1}{z-x_{\blacksquare}+\hbar_1}
\frac{z-x_{\blacksquare}-\hbar_2}{z-x_{\blacksquare}+\hbar_2}
\frac{z-x_{\blacksquare}-\hbar_3}{z-x_{\blacksquare}+\hbar_3}  \frac{y}{1-x_{\blacksquare} y} A.
\end{align*}	
Applying the right hand side of \eqref{eq:he2} to the partition $\lambda$, we obtain the right hand side is
\begin{align*}
\left( \frac{y}{1-x_{\blacksquare} y}  \psi(z)|_{\lambda} \frac{(z-y^{-1}-\hbar_1)(z-y^{-1}-\hbar_2)(z-y^{-1}-\hbar_3)}
{(z-y^{-1}+\hbar_1)(z-y^{-1}+\hbar_2)(z-y^{-1}+\hbar_3)} \right)_{+}A.
\end{align*}
The equality \eqref{eq:he2} follows from applying the operator $\Res_{y=(x_{\blacksquare})^{-1}} y^i$ on both sides.

\subsubsection{Checking the relation \eqref{eq:ef=h}}
\label{subsub:ef=h}
Let us introduce the operator $h_{i+j}:=[e_i, f_j]$ on $H_{c, T\times \C^*_{fr}}^\ast(\Hilb^*(\C^3), \varphi_{\tilde W_3}\C)^\vee\otimes_{R_T} K_T$. 
By the formula of Proposition \ref{prop:ef formula}, it is clear that 
$
[e_i, f_j]=[e_{i'}, f_{j'}], \,\ \text{if $i+j=i'+j'$}. 
$
Define the generating series 
\[
h(z):= 1-\sigma_3\sum_{i\geq 0} h_i z^{-i-1}, 
 \]
 where $\sigma_3=\hbar_1\hbar_2\hbar_3$. 
 
 The relation \eqref{eq:ef=h} is equivalent to  (see \cite[Proposition 1.5]{T})
\[
\sigma_3 \cdot (w-z)[e(z), f(w)]=h(z)-h(w). 
\]
 Then, the relation  \eqref{eq:ef=h} becomes that 
$
h(z)= \frac{1}{z-x_{\square_{000}}}
\psi(z).
$

Let $\lambda\vdash n$. By Proposition \ref{prop:ef formula}, we have
\begin{align*}
&\langle\lambda+\blacksquare|(f_j)| \lambda\rangle \langle\lambda|(e_i)| \lambda+\blacksquare\rangle \\
=& \Res_{z=x_{\blacksquare}} z^{i+j} \frac{1}{z-x_{\square_{000}}} \prod_{\square\in \lambda} \frac
{(z-x_{\square}+\hbar_1)(z-x_{\square}+\hbar_2)(z-x_{\square}+\hbar_3)}
{
(z-x_{\square}-\hbar_1)(z-x_{\square}-\hbar_2)(z-x_{\square}-\hbar_3)}
\end{align*}
On the other hand, using again Proposition \ref{prop:ef formula}, 
we have
\begin{align*}
&\langle\lambda-\blacksquare|(e_i)| \lambda\rangle
\langle\lambda|(f_j)| \lambda-\blacksquare\rangle \\
=
&\Res_{z=x_{\blacksquare}}z^{i+j}
\frac{1}{z-x_{\square_{000}}}
\prod_{\square\in \lambda-\blacksquare}
\frac
{(z-x_{\square}+\hbar_1)(z-x_{\square}+\hbar_2)(z-x_{\square}+\hbar_3)}
{
(z-x_{\square}-\hbar_1)(z-x_{\square}-\hbar_2)(z-x_{\square}-\hbar_3)}\\
=&-
\Res_{z=x_{\blacksquare}}
z^{i+j}
\frac{1}{z-x_{\square_{000}}}
\prod_{\square\in \lambda}
\frac
{(z-x_{\square}+\hbar_1)(z-x_{\square}+\hbar_2)(z-x_{\square}+\hbar_3)}
{
(z-x_{\square}-\hbar_1)(z-x_{\square}-\hbar_2)(z-x_{\square}-\hbar_3)}\\
\end{align*}
This implies that 
\begin{align*}
[e_i, f_j]_{\lambda}
=&-\sum_{\text{addible boxes}} 
  \Res_{z=x_{\blacksquare}}z^{i+j}
\frac{1}{z-x_{\square_{000}}}
\prod_{\square\in \lambda}
\frac
{(z-x_{\square}+\hbar_1)(z-x_{\square}+\hbar_2)(z-x_{\square}+\hbar_3)}
{
(z-x_{\square}-\hbar_1)(z-x_{\square}-\hbar_2)(z-x_{\square}-\hbar_3)}\\&
-\sum_{\text{removable boxes}}
\Res_{z=x_{\blacksquare}}z^{i+j}
\frac{1}{z-x_{\square_{000}}}
\prod_{\square\in \lambda}
\frac
{(z-x_{\square}+\hbar_1)(z-x_{\square}+\hbar_2)(z-x_{\square}+\hbar_3)}
{
(z-x_{\square}-\hbar_1)(z-x_{\square}-\hbar_2)(z-x_{\square}-\hbar_3)}\\
=& -\sum_{\text{addible boxes $\cup$ removable boxes}} 
\Res_{z=x_{\blacksquare}}z^{i+j}
\frac{1}{z-x_{\square_{000}}}
\prod_{\square\in \lambda}
\frac
{(z-x_{\square}+\hbar_1)(z-x_{\square}+\hbar_2)(z-x_{\square}+\hbar_3)}
{
(z-x_{\square}-\hbar_1)(z-x_{\square}-\hbar_2)(z-x_{\square}-\hbar_3)}\\
=&  \Res_{z=\infty}z^{i+j}
\frac{1}{z-x_{\square_{000}}}
\prod_{\square\in \lambda}
\frac
{(z-x_{\square}+\hbar_1)(z-x_{\square}+\hbar_2)(z-x_{\square}+\hbar_3)}
{
(z-x_{\square}-\hbar_1)(z-x_{\square}-\hbar_2)(z-x_{\square}-\hbar_3)}\\
=& \Res_{z=\infty}z^{i+j}\frac{1}{z-x_{\square_{000}}}\psi(z)|_{\lambda}. 
\end{align*}
where the last two equalities follow from the residue theorem and Lemma \ref{lem:poles}.

 Then we have
 \begin{align}\label{eq:psi}
h(z)= \frac{1}{z-x_{\square_{000}}}
\prod_{\square\in \lambda}
\frac
{(z-x_{\square}+\hbar_1)(z-x_{\square}+\hbar_2)(z-x_{\square}+\hbar_3)}
{
(z-x_{\square}-\hbar_1)(z-x_{\square}-\hbar_2)(z-x_{\square}-\hbar_3)}. 
 \end{align}
We understand the equality \eqref{eq:psi} as the equality of all coefficients for powers $z^{-(i+1)}, i\geq 0$ the power series 
expansions of   rational functions.

\begin{rmk}
More generally,  consider the quiver with potential given by 3 loops quiver, with framing vector space $\C^l$. 
\[
\begin{tikzpicture}
  \draw ($(0,0)$) circle (.08);
  \node at (0.05, -1) {$\square$};
    \node at (0.4, -1) {$l$};
     \node at (0.3, 0) {$n$};
     
    \draw [<-] (0.05, -0.1) -- (0.05, -0.9);
   \draw[->,shorten <=7pt, shorten >=7pt] ($(0,0)$)
   .. controls +(40:1.5) and +(-40:1.5) .. ($(0,0)$);
     \draw[->,shorten <=7pt, shorten >=7pt] ($(0,0)$)
   .. controls +(90+40:1.5) and +(90-40:1.5) .. ($(0,0)$);
 \draw[->,shorten <=7pt, shorten >=7pt] ($(0,0)$)
   .. controls +(180+40:1.5) and +(180-40:1.5) .. ($(0,0)$);
\end{tikzpicture}
\]
Let $\mathcal{M}(n, l)^{st}$ be the set of stable representations of the above quiver, let $\Coh_l(\C^3,n)=\mathcal{M}(n, l)^{st}/\GL_n$. Let $T_l$ be the maximal torus of $\GL_l$ acting by change of basis on the framing. 
By \cite[Theorem 5.1.1]{RSYZ}, we have an action of $\mathcal{SH}^{(Q_3, W_3)}$ on 
\[
\bigoplus_{n\geq 0} H_{c,\GL_n\times T_l\times T}^\ast(\mathcal{M}(n, l)^{st}, \varphi_{\tilde W_3}\C)^\vee
\cong \bigoplus_{n\geq 0}  H_{c, T_l\times T}^\ast(\Coh_l(\C^3, n), \varphi_{\tilde W_3}\C)^\vee
\]
given by natural correspondences. 
The same calculation as above shows that this action extends to an action of the shifted Yangian $Y_{-l}(z_1, z_2, \cdots, z_l)$ on $H_{c, T_l\times T}^\ast(\Coh_l(\C^3, *), \varphi_{\tilde W_3}\C)^\vee\otimes_{R_T} K_T$.
\end{rmk}

\section{PT moduli space of the resolved conifold}
\label{sec:PTconifold}
In this section we prove Theorem~\ref{thm:intro}(2).

\subsection{The COHA action}\label{subsec:COHAaction}
In this section, we are going to prove that on cohomology of  the moduli space associated with a stability condition chosen on the PT side of the imaginary root hyperplane, there is an action of the equivariant COHA  $\H_{\C^3}$. Furthermore, the above action can be lifted to the Drinfeld double $D(\mathcal{SH}_{\C^3})$. The latter gives rise to the action of the 1-shifted affine Yangian of $\fg\fl(1)$.

Let $(Q, W)$ be the quiver with potential from \eqref{eq:X11}. 
Let $\tilde{Q}$ be the framed quiver (see \eqref{eq: ext Q}) obtained by adding to the quiver $Q$ a framing vertex $\infty$ and an extra edge $i_\infty$ going from $\infty$ to $0$. Here $0$ denotes one of the vertices  of $Q$. 
We extend the potential for $\tilde{Q}$ by the formula $\tilde{W}=W$. 
For $\tilde{\zeta}=(\zeta_0,\zeta_1, \zeta_{\infty})$ a triple of real numbers 
and a triple of vector spaces $(V_0, V_1, V_{\infty})$ associated to the following 3 vertices of extended quiver $\tilde{Q}$, a representation $F=(a_1, a_2, b_1, b_2, \iota)$ of $\tilde{Q}$ with dimension vector $(\dim(V_0), \dim(V_1), \dim(V_{\infty}))$,
we define the slope of $F$ to be
\[
\theta_{\tilde{\zeta}}(F):=\frac{\zeta_0 \dim(V_0)+\zeta_{1} \dim(V_1)+\zeta_{\infty} \dim(V_{\infty})}{\dim(V_0)+\dim(V_1)+\dim(V_{\infty})}. 
\]
A representation $F$ of $\tilde{Q}$  is said to be $\theta_{\tilde{\zeta}}$-(semi-)stable if 
\[
\theta_{\tilde{\zeta}}(F')< (\leq)\theta_{\tilde{\zeta}}(F)
\] 
 for any nonzero proper $\tilde{Q}$-subrepresentation $F'$ of $F$.

For any pair of real numbers $\zeta=(\zeta_0, \zeta_1)$, define $\zeta_{\infty}:=-\zeta_0 \dim(V_0)-\zeta_1\dim(V_1)$ so that we have a triple $\tilde{\zeta}=(\zeta_0,\zeta_1, \zeta_{\infty})$. 
 We say the representation $F$ with dimension vector $(\dim(V_0), \dim(V_1), 1)$ is $\zeta$-stable, if it is $\theta_{\tilde{\zeta}}$-stable. 
 Note that in the special case $\zeta_0=\zeta_1=-1$, 
 a representation of $\tilde{Q}$ is $\theta_{\tilde{\zeta}}$-stable 
  if every $\tilde{Q}$-subrepresentation of $(V_0, V_1, V_{\infty})$ containing the framing $V_{\infty}$ is the entire $(V_0, V_1, V_{\infty})$. For this choice of $\zeta$ we also call a $\zeta$-stable representation  {\it cyclicly stable}.

Let $\fM_\zeta(v_0, v_1)$ denote the moduli space of $\zeta$-stable $\tilde{Q}$-representations, such that $\dim(V_0)=v_0, \dim(V_1)=v_1, \dim(V_{\infty})=1$. 
The space $\coprod_{v_0,v_1}\fM_\zeta(v_0, v_1)$ is isomorphic to  the moduli space  $\fM_\zeta$ of $\zeta$-stable perverse coherent systems in the sense of Nagao and Nakajima \cite[\S~1.3]{NN}.
We refer to $\zeta=(\zeta_0, \zeta_1)$ as the stability parameter, and the 2-dimensional real vector space $\{\zeta=(\zeta_0, \zeta_1)\mid \zeta_i\in\R\}$ the space of stability conditions.

Let us recall some terminologies from the Introduction. A  stability parameter $\zeta$ is called \textit{generic} if every $\zeta$-semi-stable point is in fact $\zeta$-stable. 
The locus of non-generic stability conditions is a union  of hyperplanes (which in this case is a union of lines) called {\it walls}. The complement of the hyperplane arrangement is the union of open connected components which we call {\it chambers}. The moduli space $\fM_\zeta$ in the case of generic $\zeta$  depends only on the chamber that contains $\zeta$.  This is the case we are mostly interested in.

The chamber structure of the space of stability conditions for $X_{1,1}$ considered here is described in \cite[Figure 1]{NN}. In particular, one wall is given by the imaginary root hyperplane $\zeta_0+\zeta_1=0$. We consider the chambers above this imaginary root hyperplane, i.e. those where $\zeta_0+\zeta_1>0$. Those walls are given by 
\[
(m+1)\zeta_0+m\zeta_1=0,\,\ \text{ where $m\in \Z_{\geq 0}$.}
\] 
Consider  
a generic stability condition $\zeta_{m}$ which is above this wall but close to it. We use a different description of 
the moduli space $\fM_{\zeta_{m}}$ more relevant for our purpose. 
For this, we consider the  following quiver $\tilde{Q}_\zeta$ with potential $\tilde{W}_{\zeta}$ (see \cite[Figure 9]{NN}). Note that $\tilde{Q}_\zeta$ is also obtained from $Q$ by adding the vertex $\infty$ which by an abuse of terminology we still refer to as the framing,  but the set of arrows depends on $m$. We consider representation of $\tilde{Q}_\zeta$ on a triple of vector spaces $(V_0, V_1, V_{\infty})$ associated to the 3 vertices as before. 
Again a representation of $\tilde{Q}_\zeta$ is said to be {\it cyclicly stable} if every $\tilde{Q}_\zeta$-subrepresentation of $(V_0, V_1, V_{\infty})$ containing the framing $V_{\infty}$ is the entire $(V_0, V_1, V_{\infty})$ \cite{NN}.
The moduli space   $\fM_m$ of   cyclically stable $\tilde{Q}_\zeta$-representations with $V_{\infty}=\C$ is again isomorphic to $\fM_{\zeta_{m}}$ by \cite[Theorem 4.13]{NN}:
\[
\fM_{\zeta_{m}}(v_0, v_1)
=
\fM_m((m-1)v_0-mv_1, mv_0-(m+1) v_1). 
\]

    $$\begin{tikzpicture}[scale=0.8]
        	\node at (-2.2, 0) {$\bullet$};   	\node at (2.2, 0) {$\bullet$};  
		\draw[-latex,  bend left=30, thick] (-2, 0.1) to node[above]{$a_2$} (2, 0.1);
		\draw[-latex,  bend left=50, thick] (-2, 0.3) to node[above]{$a_1$} (2, 0.3);
	\draw[-latex,  bend left=30, thick] (2, -0.1) to node[above]{$b_1$} (-2, -0.1);
		\draw[-latex,  bend left=50, thick] (2, -0.3) to node[above]{$b_2$} (-2, -0.3);
					\node at (-2.6, 0) {$V_1$};\node at (2.6, 0) {$V_2$};
\node at (2.2, -4) {$\square$};  
\draw[->, thick] (2.6, -3.8) -- (2.6, -0.4) ;
 \node at (2.3, -2) {$\cdots$};  
 \draw[->, thick] (2, -3.8) -- (2, -0.4) ;
  \node at (4, -2) {$q_1, \cdots, q_m$};  
 
    \draw[-latex,  bend right=30, thick] (-2.2, -0.2) to node[above]{} (2, -3.8) ;   
     \node at (-1.6, -2) {$\cdots$};   
      \node at (-3.6, -2) {$p_1, \cdots, p_{m+1}$};  
     \draw[-latex,  bend right=30, thick] (-2.6, -0.6) to node[above]{} (2, -4.2) ; 
       \node at (0, -4.7) {The quiver $\tilde{Q}_{\zeta_{m}}$ with potential $\tilde{W}_{\zeta}$ given by};
       \node at (0, -5.4){$a_1b_1a_2b_2-a_1b_2a_2b_1+p_1b_1q_1+p_2(b_1q_2-b_2q_1)+\cdots+p_m(b_1 q_m-b_2q_{m-1})-p_{m+1}b_2q_m$};
    \end{tikzpicture}
   $$      
   
 There is a 4-dimensional torus $(\C^*)^4$ acting by scaling $(a_1,a_2, b_1, b_2)$. However, the 1-dimensional subtorus $T_0=\{(t,t,t^{-1},t^{-1})\in (\C^*)^4\}$ acts trivially. Then we have the action of  the 3-dimensional subtorus of  $(\C^*)^4$ with the first coordinate being $1\in\C^*$. Note that this 3-dimensional torus maps isomorphically to $(\C^*)^4/T_0$ under the quotient map. We write an element in this subtorus as $(1,t,q,h)\in (\C^*)^4$. We consider $T\subseteq (\C^*)^4$ consisting of $(1,t,q,h)$ with $tqh=1$. The action of $T$ on $X$ preserves the canonical line bundle. In particular, the induced action on the moduli spaces preserves the potential.

 The weights of $p_i$ and $q_i$ are determined by the conditions
 \[
 wt(p_i)=wt(q_i^{-1}), wt(p_{i+1})=wt(q_i^{-1}) t^{-1}.   
 \] 
Therefore, we have 
\[
wt(b_1)=1, wt(b_2)=t, wt(a_1)=q, wt(a_2)=t^{-1}q^{-1}, wt(p_i)=t^{1-i}, wt(q_i)=t^{i-1}. 
\]

The moduli space $\fM_m$ contains an open subset $\fM_m^0$ singled out by the condition that $b_1$ is surjective. The natural flat pull-back to the open subset  induces the isomorphism $H_{c,T\times \C^*_{fr}}^*(\fM_m,\varphi_{\tr \tilde{W}_{\zeta}}\C)^\vee\otimes_{R_T} K_T\to H_{c,T\times \C^*_{fr}}^*(\fM_m^0,\varphi_{\tr \tilde{W}_{\zeta}}\C)^\vee\otimes_{R_T} K_T$.

\begin{thm}\label{thm:conifoldYangian}
\begin{enumerate}
\item
There is a natural action of $\H_{\C^3}$ on $H_{c,T\times \C^*_{fr}}^*(\fM_m^0,\varphi_{\tr \tilde{W}_{\zeta}}\C)^\vee$. 
\item
This action extends to an action of $Y_1(z_1)$ on $H_{c,T\times \C^*_{fr}}^*(\fM_m^0,\varphi_{\tr \tilde{W}_{\zeta}}\C)^\vee\otimes_{R_T}K_T$.
\end{enumerate}
\end{thm}

We prove here only part 1 of the statement.
In order to  prove part 2 it suffices to find the commutation relations of the raising and lowering operators acting on $ H_{c,T\times \C^*_{fr}}^*(\fM_m^0,\varphi_{\tr \tilde{W}_{\zeta}}\C)^\vee\otimes_{R_T}K_T$, which we leave until \S~\ref{subsec:conifold_operators}.

{\it Proof.}
We keep the notation of Section \ref{subsec:conifold}. To avoid confusion, here in the proof, we denote by $W$  in place of $V_\infty$ the vector space at the framing vertex of the quiver. 
Let 
\[
Surj(V_2', V_1')=\{x\in \Hom(V_2', V_1')\mid \text{$x$ is surjective}\}
\]
be the set of surjective morphisms from $V_2'$ to $V_1'$.
Denote $\underline{V}'=(V_1', V_2')$.
Let
\[
\calM(\underline{V}', W)^{0}:=\Big(\Surj(V_2', V_1')\oplus \Hom(V_2', V_1')\oplus \Hom(V_1', V_2')^{\oplus 2} \oplus \Hom (V_1', W)^{\oplus m+1} \oplus \Hom (W, V_2')^{\oplus m}\Big)
\]
Note that the quotient $Surj(V_2', V_1')/ \GL(V_1')$ is a Grassmannian.
Consider the following space: 
\begin{align*}
\GL(V_1')\times_{P_1} \calM(\underline{V}', W)^{0}
&=\{(\xi, b_1, b_2, a_1, a_2, p_1, \cdots, p_{m+1}, q_1, \cdots, q_m) \mid \\
&
\xi\subset V_1',
 b_1: V_2'\surj V_1' \,\ \text{is surjective}, b_2: V_2'\to V_1', a_i: V_1'\to V_2', i=1, 2, \\
 &p_1, \cdots, p_{m+1}\in \Hom(V_1', W), \text{and}\,\ 
 q_1, \cdots q_m\in \Hom(W, V_2')
 \}. 
\end{align*}

Denote by $\calM(\underline{V}', W)^{st, 0} \subset \calM(\underline{V}', W)^{0}$ the stable locus. 

We define the following spaces
\begin{align*}
&\widetilde{\calM}=\{(\xi_1, \xi_2, b_1, b_2, a_1, a_2, p_1, \cdots, p_{m+1}, q_1, \cdots, q_m)\mid \xi_2\subset V_2', \xi_1\subset V_1', b_1: V_2'\surj V_1',  \\
&\phantom{1234567890}\text{such that $b_1(\xi_2)\subset \xi_1$, $b_1|_{\xi_2}: \xi_2\to \xi_1$ is an isomorphism}, b_2: V_2'\to V_1', a_i: V_1'\to V_2', i=1, 2\},  \\
 &\phantom{1234567890}\ \text{and}\,\ p_1, \cdots, p_{m+1}\in \Hom(V_1', W), \text{and}\,\ 
 q_1, \cdots q_m\in \Hom(W, V_2')
 \}. \\
&G\times_P Z^{0}=\{(\xi_1, \xi_2, b_1, b_2, a_1, a_2, p_1, \cdots, p_{m+1}, q_1, \cdots, q_m) \in \widetilde{\calM} \mid  b_2(\xi_2)\subset \xi_1, a_i(\xi_1)\subset \xi_2, i=1, 2,\\
&\phantom{1234567890}
 \text{and}\,\ p_k(\xi_1)=0, \text{for $k=1, \cdots, m+1$} \}. 
\end{align*}

We have a closed embedding $
i: G\times_P Z^{0}\inj \widetilde{\calM}$ and an affine bundle projection $\pi: \widetilde{\calM}\to  \GL(V_1')\times_{P_1} \calM(\underline{V}', W)^{0}$. 

We use the following notations
\[
\underline{V}=(\xi_1, \xi_2), \underline{V}'=(V_1', V_2'), \text{and $\underline{V}''=\underline{V}'/ \underline{V}=(V_1'', V_2'')$}. 
\]
For an element $ (\xi_1, \xi_2, b_1, b_2, a_1, a_2, p_1, \cdots, p_{m+1}, q_1, \cdots, q_m)\in G\times_P Z^{0}$, restricting the element on the subspace $\underline{V}\subset \underline{V}'$ gives an element in $G\times_P \Rep(Q_{1,1}, \underline{V})^{0}$, taking quotient by the subspace $\underline{V}\subset \underline{V}'$ gives an element in $G\times_P\calM(\underline{V}'', W)^{0}$. Denote the induced map by \[
\phi:   G\times_P Z^0\to G\times_P\big(\Rep(Q_{1,1}, \underline{V})^{0}  \times  \calM(\underline{V}'', W)^{0}\big).\]
We then have the following correspondence
\begin{equation}\label{diag:action1}
\xymatrix{
&G\times_P Z^0  \ar[ld]_(0.4){\phi} \ar@{^{(}->}[r]^{i} &  \widetilde{\calM}  \ar[r]^(0.3){\pi}  &  \GL(V_1')\times_{P_1}\calM(\underline{V}', W)^{0} \ar[d]^{\psi}\\
G\times_P \left({\begin{matrix} \Rep(Q_{1,1}, \underline{V})^{0} \\ \times  \calM(\underline{V}'', W)^{0} \end{matrix}}\right) 
&&&   \calM(\underline{V}', W)^{0}
}
\end{equation}
We have the following isomorphisms
\begin{align*}
&\frac{G\times_P (\Rep(Q_{1,1}, \underline{V})^{0} \times  \calM(\underline{V}'', W)^{st, 0} )}{\GL(V_1')}\\
\cong  & GL(V_2')\times_{P_2} \frac{(\Rep(Q_{1,1}, \underline{V})^{0} \times  \calM(\underline{V}'', W)^{st, 0} )}{P_1}\\
\cong & GL(V_2')\times_{P_2} \frac{(\Rep(Q_{1,1}, \underline{V})^{0} \times  \calM(\underline{V}'', W)^{st, 0} )}{\GL(V_1)\times \GL(V_1'') }\\
\cong &GL(V_2')\times_{P_2}  (\Rep(Q_3, V_1) \times  (\calM(\underline{V}'', W)^{st, 0} / \GL(V_1'') )
\end{align*}
Taking the quotient of the correspondence \eqref{diag:action1} by $\GL(V_1')$, we have the correspondence
\begin{equation}\label{diag:action2}
\xymatrix@C=1em{
&GL(V_2')\times_{P_2} (Z^0/P_1)  \ar[ld]_(0.4){\phi} \ar@{^{(}->}[r]^(0.6){i} &  \widetilde{\calM}/ \GL(V_1') \ar[r]^{\pi}  & \calM(\underline{V}', W)^{0}/P_1 \ar[d]^{\psi}\\
GL(V_2')\times_{P_2}  \left({\begin{matrix} \Rep(Q_3, V_1) \times \\ \calM(\underline{V}'', W)^{0}/ \GL(V_1'')  \end{matrix}}\right) 
&&&   \calM(\underline{V}', W)^{0}/ \GL(V_1')
}
\end{equation}
The action is defined by 
\begin{equation}\label{act}
\psi_* \circ (\pi^{*})^{-1} \circ i_* \circ  \phi^*, 
\end{equation}
on the stable locus of \eqref{diag:action2}, 
where $(\pi^{*})^{-1}$ exists since $\pi$ is an affine bundle. 
$\blacksquare$

\subsection{The tangent space}
\label{subsec:tangent space PT}
In this section, we describe the tangent spaces to the torus fixed points in the open locus $\fM_m^0\subset \fM_m$,  as well as the correspondence used in \eqref{diag:action2}. Again in this section to avoid confusion we denote by $W$ the vector space at the framing vertex of the quiver in place of $V_\infty$.

Recall that $\fM_m^0(\underline{V})$ is the quotient $\calM(\underline{V}, W)^{st, 0}/(\GL(V_1)\times \GL(V_2))$, where $\calM(\underline{V}, W)^{st, 0}$ is the cyclically stable locus of 
\[
\calM(\underline{V}, W)^{0}=\Big(\Surj(V_2, V_1)\oplus \Hom(V_2, V_1)\oplus \Hom(V_1, V_2)^{\oplus 2} \oplus \Hom (V_1, W)^{\oplus m+1} \oplus \Hom (W, V_2)^{\oplus m}\Big).
\]

The quotient $Surj(V_2, V_1)/ \GL(V_1)$ is a Grassmannian,  denoted by $\Grass(v_1,V_2)$, 
parameterizing quotients $b_1:V_2\surj V_1$ with the dimension of $V_1$ being $v_1$. 
Therefore, $\calM(\underline{V}, W)^{0}/\GL(V_1)$ is an affine bundle over the Grassmannian with fiber 
\[
F:= \Hom(V_2, V_1)\oplus \Hom(V_1, V_2)^{\oplus 2} \oplus \Hom (V_1, W)^{\oplus m+1} \oplus \Hom (W, V_2)^{\oplus m}. 
\]

The tangent space of $\fM_m^0(\underline{V})$ is an extension of the above vector space $F$ by the  tangent space of $\Grass(v_1,V_2)$, taking into account of the $\GL(V_2)$-action. 
Again we  denote by 
\[
0\to \overline{V_2}\to V_2 \to V_1 \to 0.
\]  
the tautological sequence of vector bundles on $\Grass(v_1,V_2)$. Then, the tangent space of  $\Grass(v_1,V_2)$ is $\sHom(\overline{V_2}, V_1)$. 
The tangent bundle $\calT$ to $\calM(\underline{V}, W)^{0}/\GL(V_1)$ has a filtration with associated graded given by
\begin{align*}
Gr(\calT)=&\sHom(\overline{ V_2}, V_1)\\
\oplus&(\sHom(V_2,  V_1)\oplus \sHom( V_1, V_2)^{\oplus 2}\oplus \sHom(V_1, W)^{\oplus (m+1)} \oplus \sHom(W, V_2)^{\oplus m})
\end{align*}
The action of $\GL(V_2)$ induces a Lie algebra map from $\End(V_2)$ to the sheaf of sections of $\calT$ endowed with the standard Lie bracket on vector fields. 

The tangent bundle of the quotient $\calM(\underline{V}, W)^{0}/(\GL(V_1)\times \GL(V_2))$ is given by the cokernel of the above action map. 

We now compute the tangent space of the correspondence $Z^{st, 0}/(P_1\times P_2)$ in \eqref{diag:action2}. 
In the above notation let us consider the space 
\begin{align*}
&\widetilde{\Fl}:=\{(\xi_1, \xi_2, b_1)\mid \xi_2\subset V_2', \xi_1\subset V_1', b_1: V_2'\surj V_1'\\
&\phantom{1234567890}\text{such that $b_1(\xi_2)\subset \xi_1$, $b_1|_{\xi_2}: \xi_2\to \xi_1$ is an isomorphism}\}/\GL(V_1')
\end{align*}
which parameterizes diagrams of vector spaces as follows
\[
\xymatrix{
0\ar[r]&\xi_1 \ar[r] & V_1' \ar[r] & V_1 \ar[r]& 0\\
0\ar[r]&\xi_2 \ar[u]^{b_1}_{\cong} \ar[r] & V_2'  \ar@{->>}[u]^{b_1}\ar[r] & V_2 \ar@{->>}[u]^{b_1} \ar[r]&0\\
&& \overline{V_2'}  \ar@{^{(}->}[u]& \overline{V_2}\ar@{^{(}->}[u] &
}
\]
where $\dim(V_i)=v_i$ and $\dim(V_i')=v_i+1$, for $i=1, 2$. 
Note that $(G\times_P Z^{0})/\GL(V_1')$ is an affine bundle over $\widetilde{\Fl}$.

The fiber of the affine bundle is isomorphic to 
\[
P(V_2',   V_1')\oplus P(  V_1', V_2')^{\oplus 2}\oplus \sHom(  V_1', W)^{\oplus (m+1)}  \oplus \sHom(W, V_2')^{\oplus m}
\]
where
\begin{align*}
P(V_2',   V_1'):=\{x\in \sHom(V_2',   V_1')\mid x(\xi_2)\subset \xi_2 \}, \\
P(V_1',   V_2'):=\{x\in \sHom(V_1',   V_2')\mid x(\xi_1)\subset \xi_1 \}. 
\end{align*}

Let  $\Fl(v_1+1,1,V_2')$ be the flag variety of pairs consisting of 
a quotient $b_1:V_2' \surj V_1'$ with $\dim(V_1')=v_1+1$ and a line $\xi_1 \subseteq V_1'$.  Then the base space $\widetilde{\Fl}$ maps to $\Fl(v_1+1,1,V_2')$, where  the map is given by forgetting the subspace $\xi_2$. This map realizes $\widetilde{\Fl}$ as an affine bundle with each fiber isomorphic to $\Hom(\xi_2, \overline{V_2'})$.

The tangent bundle of the flag variety $\Fl(v_1+1,1,V_2')$ has the same class in the Grothendieck group as 
\[
\sHom(\overline{V_2}, V_1')\oplus \sHom(\xi, V_1). 
\]
Hence there is a filtration on the tangent sheaf of $(G\times_P Z^{0})/\GL(V_1')$ with associated graded being 
\begin{align*}
&\sHom(\overline{V_2}, V_1')\oplus \sHom(\xi, V_1)\\
\oplus &P(V_2',   V_1')\oplus P(  V_1', V_2')^{\oplus 2}\oplus \sHom(  V_1', W)^{\oplus (m+1)}  \oplus \sHom(W, V_2')^{\oplus m}\\
\oplus &\Hom(\xi_2, \overline{V_2'}). 
\end{align*}
Note that the affine bundle  $\pi: \widetilde{\calM}\to  \GL(V_1')\times_{P_1} \calM(\underline{V}', W)^{0}$ (see \eqref{diag:action1}) has fiber isomorphic to $\Hom(\xi_2, \overline{V_2'})$.

In the action defined by \eqref{act}, we only focus on the contribution of the tangent sheaf of $(G\times_P Z^{0})/\GL(V_1')$ without the term $\Hom(\xi_2, \overline{V_2'})$.

As before, the action of $\GL(V_2')$ induces a Lie algebra map 
\[
\End(V_2')\to \calT_{(G\times_P Z^{0})/\GL(V_1')}.
\]
The tangent bundle of the quotient $Z^{st, 0}/(P_1\times P_2)=(G\times_P Z^{0})/(\GL(V_1') \times \GL(V_2'))$ is again given by the cokernel of the above action map.

\subsection{Fixed points}\label{subsec:pyramid}
In this subsection we will use the terminology and some results of \cite{NN}.
The following are the pictures for the empty room configurations (ERC) for the finite type
pyramid partitions with length 3 and length 4 (\cite[Figure~12]{NN}). 

\begin{equation}
\begin{tikzpicture}[scale=0.5]
\draw[fill=black!40!white] (-2, 0) circle (1cm);
 \draw[fill=black!40!white] (0, 0) circle (1cm);
 \draw[fill=black!40!white] (2, 0) circle (1cm);
  
 \draw[fill=white] (-1, 0) circle (1cm); 
 \draw[fill=white] (1, 0) circle (1cm);
    
\draw[fill=black!40!white] (1, 1) circle (1cm);
\draw[fill=black!40!white] (1, -1) circle (1cm);
\draw[fill=black!40!white] (-1, -1) circle (1cm);
\draw[fill=black!40!white] (-1, 1) circle (1cm);

   \draw[fill=white] (0, 1) circle (1cm); \draw[fill=white] (0, -1) circle (1cm); 
   
   \draw[fill=black!40!white] (0, 0) circle (1cm);
   \draw[fill=black!40!white] (0, -2) circle (1cm);
      \draw[fill=black!40!white] (0, 2) circle (1cm);    
      \node at (0, -5) {The ERC for finite type };
        \node at (0, -5.8) {pyramid partitions with length 3};
\end{tikzpicture}
\,\ \,\ \,\ \,\ 
\begin{tikzpicture}[scale=0.5]
      \draw[fill=black!40!white] (-3, 0) circle (1cm);
   \draw[fill=black!40!white] (-1, 0) circle (1cm);
      \draw[fill=black!40!white] (1, 0) circle (1cm);
   \draw[fill=black!40!white] (3, 0) circle (1cm);

   \draw[fill=white] (0, 0) circle (1cm); 
   \draw[fill=white] (-2, 0) circle (1cm); 
   \draw[fill=white] (2, 0) circle (1cm); 

\draw[fill=black!40!white] (0, -1) circle (1cm);
\draw[fill=black!40!white] (-2, -1) circle (1cm);
\draw[fill=black!40!white] (2, -1) circle (1cm);
\draw[fill=black!40!white] (0, 1) circle (1cm);
\draw[fill=black!40!white] (-2, 1) circle (1cm);
\draw[fill=black!40!white] (2, 1) circle (1cm);

\draw[fill=white] (1, 1) circle (1cm);
\draw[fill=white] (1, -1) circle (1cm);
\draw[fill=white] (-1, -1) circle (1cm);
\draw[fill=white] (-1, 1) circle (1cm);

\draw[fill=black!40!white] (-1, 0) circle (1cm);
\draw[fill=black!40!white] (-1, -2) circle (1cm);
\draw[fill=black!40!white] (-1, 2) circle (1cm);
\draw[fill=black!40!white] (1, 0) circle (1cm);
\draw[fill=black!40!white] (1, -2) circle (1cm);
\draw[fill=black!40!white] (1, 2) circle (1cm);

   \draw[fill=white] (0, 0) circle (1cm); 
   \draw[fill=white] (0, -2) circle (1cm); 
   \draw[fill=white] (0, 2) circle (1cm); 
   
      \draw[fill=black!40!white] (0, -3) circle (1cm);
   \draw[fill=black!40!white] (0, -1) circle (1cm);
      \draw[fill=black!40!white] (0, 1) circle (1cm);
   \draw[fill=black!40!white] (0, 3) circle (1cm);
   \node at (0, -4.5) {The ERC for finite type };
     \node at (0, -5.3) {pyramid partitions with length 4};
\end{tikzpicture}
\end{equation}
In general, for ERC for the finite type pyramid partitions with length $m$, there are $1\times m$ black stones on
the first layer, $1\times (m-1)$ white stones on the second layer,  $2 \times (m-1)$ black stones on the
third, $2 \times (m-2)$ white stones on the fourth, and so on until we reach $m\times 1$ black stones.

The following definition can be found in \cite{CJ}. 
\begin{defn}\label{pyramid partition}
A finite type pyramid partition of length $m$ is a finite subset $\Pi$ of 
the ERC of length $m$  in which, for every stone
in $\Pi$, the stones directly above it are also in $\Pi$.
\end{defn}

\begin{prp}\cite[Proposition 4.14]{NN}
The set of $T$-fixed points in $\fM_m$ is isolated and parameterized by finite type pyramid partitions of length $m$.
\end{prp}

By our conventions on the torus action, the weights of the black stones on the top layer in the empty room partition are
\[
1, t, t^2, \cdots, t^{m-1}. 
\]
The weights of the black stones on the layer 3 are
\[
qt, ht;  qt^2, ht^2;  \cdots, qt^{m-1}, ht^{m-1}. 
\]
The weights of the black stones on the layer 5 are
\[
q^2t^2, qh t^2, h^2t^2, q^2t^3, qh t^3, h^2t^3, \cdots, q^2t^{m-1}, h^2t^{m-1}. 
\]
The last one is the $2(m-1)+1=(2m-1)$th layer, where the weights of the black stones  are
\[
q^{m-1}t^{m-1}, q^{m-2}ht^{m-1}, q^{m-3}h^2t^{m-1}\cdots, h^{m-1}t^{m-1}. 
\]
By the definition of the torus weight of $b_1$, for each pair of black and white stones as in Figure~\eqref{a pair}, the weight of the  white stone is the same as that of the black stone. 

Let us consider a pair consisting of a black stone and a white stone such that the black stone is right above of the white stone (as in the the following picture ):
\begin{equation}\label{a pair}
\begin{tikzpicture}[scale=0.5]
        \draw (0, 3) circle (1cm);
\draw[fill=black!40!white] (0, 2) circle (1cm);
\end{tikzpicture}
\end{equation}
We say that  a  pair of stones as in \eqref{a pair} is {\it removable}, if after removing the pair, we still get a pyramid partition. 
Similarly, we say that a pair of stones as in \eqref{a pair}  is {\it addible} to a pyramid partition, if they are not in the pyramid partition but after adding the pair to the pyramid partition, we still get one. 

Let us  analyze  what the  addible and removable pairs are for a given  pyramid partition. For this we use the following terminology in order to describe the relative position of one stone with respect to another one:
\[
\begin{tikzpicture}
       \node at (0, 1.5) {front};
       \draw[-latex, thick] (0, 0.2) to  (0, 1);      
       \draw[-latex, thick] (0, -0.2) to  (0, -1);
       \node at (0, -1.5) {back }; \node at (0, -2) {($t$)};
        \draw[-latex, thick] (-0.2, 0) to  (-1, 0);
        \node at (-1.5, 0) {left};
         \node at (-1.5, -0.4) {($q$)};
         \draw[-latex, thick] (0.2, 0) to  (1, 0);
        \node at (1.5, 0) {right};
          \node at (1.5, -0.4) {($h$)};
       \end{tikzpicture}
\]
Furthermore, in this  terminology the word \textit{above} will mean {\it up with respect to the paper surface}, and \textit{below} will mean {\it down with respect the paper surface}.

A pair of black and white stones is addible, if the following conditions hold for the black stone in this pair:
\begin{enumerate}
\item(Black condition) There is a black stone in front of it  in the pyramid partition. 
\item(White condition) If there are white stones above it in the empty room partition (see the following picture), then these white stones have to be in the pyramid partition. 
\[
\begin{tikzpicture}[scale=0.5]
\draw[fill=black!40!white] (0, 0) circle (1cm);
\end{tikzpicture} \,\ \text{or} \,\ 
\begin{tikzpicture}[scale=0.5]
\draw[fill=black!40!white] (0, 0) circle (1cm);
\draw[fill=white] (1, 0) circle (1cm);
\end{tikzpicture} \,\ \text{or} \,\ 
\begin{tikzpicture}[scale=0.5]
\draw[fill=black!40!white] (0, 0) circle (1cm);
\draw[fill=white] (1, 0) circle (1cm);
\draw[fill=white] (-1, 0) circle (1cm);
\end{tikzpicture}
\]
\end{enumerate}
For a pair of black and white stones let us analyze the  above two conditions for the black stone in this pair. In order to avoid confusion we refer to the black stone in the pair as  \textit{black-one}.
If the black condition holds for the black-one, then there is another black stone, which we call  \textit{black-zero}, which is positioned in front of the  black-one. Moreover, the white condition holds automatically for black-zero. 
Unless one of the white stones above black-zero is the end of a chain of the whites, the white condition also holds for  black-one.

Similarly, 
a pair of black and white stones as in \eqref{a pair} is removable, if the pair is at the end of a chain like this
\[
\begin{tikzpicture}[scale=0.5]
        \draw (0, 3) circle (1cm);
    \draw (0, 1) circle (1cm);
\draw[fill=black!40!white] (0, 4) circle (1cm);
\draw[fill=black!40!white] (0, 2) circle (1cm);
\draw[fill=black!40!white] (0, 0) circle (1cm);
\node at (0, -2) {$\vdots$};
   \draw (0, -4) circle (1cm);
   \draw[fill=black!40!white] (0, -5) circle (1cm);
   \draw [->, in=50,  thin] (4, -3) to node[anchor=west, xshift=-0.3cm, yshift=-0.5cm]
   {No black stone below it} (0, -3.5);
\end{tikzpicture}
\] and there is no black stone below the last white stone. 

\begin{prp}
Assume $t+q+h=0$. Poles of the function 
\begin{align*}
h(z)=&(-1)^{|\{\text{black only stones}\}|}(-1)^{m+1}(z-\chi-mt)\\
&\prod_{b\in {\text{black stones and white stones}}} \frac{(z-x_b+t)(z-x_b+q)(z-x_b+h)}{(z-x_b-t)(z-x_b-q)(z-x_b-h)}\\&
 \prod_{b\in \text{black stones only} }  \frac{(z-x_b)(z-x_b+q)(z-x_b+h)}{(z-x_b-t)} 
\end{align*}are at the addible and removable pairs. 
\end{prp}
{\it Proof.} By definition the function $h(z)$ is product of factors coming from each stone.  Let us analyze the contribution of each chain of stones. 
From the above discussion of addible and removable pairs, in order to prove the proposition it suffices to show that $h(z)$ has the following zeros and poles. Here the formula of $h$ is given as a product over all the stones in $\lambda$. We calculate the contribution of one chain of stones in the product, by taking the product of the corresponding factors over all the stones in this chain.

Zeros:
\begin{enumerate}
\item For a chain of black stones, we want the product of factors from this chain to have  two zeros above the end stone. That is, assume the end stone is $x$, then we want two zeros at $x-h$ and $x-q$ respectively. 
\item For a chain of white stones we want the product of factors from this chain to have two zeros,
 under the stones in front of the end stone of the chain. That is, let the end stone of this chain be $x$, then the two zeros are at $x+t+q$ and $x+t+h$.  ( The stone in front of $x$ is $x+t$. The two places $x+t+q$ and $x+t+h$, where we want the zeros to be, are the two stones underneath $x+t$. 
\item Similarly, a chain of black stones should contribute one zero at the beginning stone. 
\end{enumerate}
Poles:
\begin{enumerate}
\item A chain of black stones (with or without white stones) should contribute a pole at the end stone of the chain. 
\item A chain of black stones and white stones should contribute a pole at the end stone. 
\end{enumerate}

Now we verify that  $h(z)$ does have the aforementioned zeros and poles from the contribution of each chain of stones.
Suppose we have a chain of black stones (with marked weights) as on the following picture
$$\begin{tikzpicture}[scale=0.5]
\draw[fill=black!40!white] (0, 4) circle (1cm);
\node at (0, 4) {$x$};
\draw[fill=black!40!white] (0, 2) circle (1cm);
\node at (0, 2) {$xt$};
\draw[fill=black!40!white] (0, 0) circle (1cm);
\node at (0, 0) {$x t^2$};
\node at (0, -2) {$\vdots$};
\draw[fill=black!40!white] (0, -4) circle (1cm);
\node at (0, -4) {$x t^{a-1}$};
    \end{tikzpicture}$$ 
The function $h(z)$ has the following factor
\begin{align*}
h(z)
    &=\prod_{i=0}^{a-1}  \frac{(z-x-it)(z-x-it+q)(z-x-it+h)}{(z-x-it-t)}\\
    &=\frac{z-x}{z-x-at}\prod_{i=0}^{a-1}(z-x-it+q)(z-x-it+h).    \end{align*}
    Thus, the pole of $h(z)$ is $xt^{a}$. 
    It is clear that $xt^a$ is an addible place. 
    
Suppose we have a chain of black and white stones (with marked weights) as on the following picture
$$  \begin{tikzpicture}[scale=0.5]
        \draw (0, 3) circle (1cm);
    \draw (0, 1) circle (1cm);
\draw[fill=black!40!white] (0, 4) circle (1cm);
\node at (0, 4) {$x$};
\draw[fill=black!40!white] (0, 2) circle (1cm);
\node at (0, 2) {$xt$};
\draw[fill=black!40!white] (0, 0) circle (1cm);
\node at (0, 0) {$x t^2$};
\end{tikzpicture}
$$
To illustrate the idea of the proof, we compute this case in the explicit example as in the picture. This can be made in general. 
In this case, the function $h(z)$ has the following factor
\begin{align*}
 &\frac{(z-x)(z-x+q)(z-x+h)}{(z-x-t)} 
 \frac{(z-x-t+t)(z-x-t+q)(z-x-t+h)}{(z-x-t-t)(z-x-t-q)(z-x-t-h)}\\
  &\frac{(z-x-2t+t)(z-x-2t+q)(z-x-2t+h)}{(z-x-2t-t)(z-x-2t-q)(z-x-2t-h)}
 \\
 =&\frac{(z-x)^2(z-x-2t+q)(z-x-2t+h)}{(z-x-2t)(z-x-3t)}
\end{align*}
 Thus, the poles of $h(z)$ are $x+2t$ and $x+3t$. 
    It is clear that $xt^3$ is an addible place, and $xt^2$ is a removable place. 
$\blacksquare$. 

\subsection{Raising and lowering operators}
\label{subsec:conifold_operators}
We define the raising and lowering operators on
\[
\bigoplus_{\underline{V}}H_{c,T\times \GL_W}^*(\fM_m^0(\underline{V}, W),\varphi_{\tr \tilde{W}_{\zeta}}\C)^\vee \otimes_{R_T} K_T
\cong \bigoplus_{\underline{V}} H_{c,\GL_{\underline{V}} \times T\times \GL_W}^\ast(\mathcal{M}(\underline{V}, W)^{st, 0}, \varphi_{\tr \tilde{W}_{\zeta}}\C)^\vee\otimes_{R_T} K_T
\]
similar as in \S\ref{sub:Action of the shifted Yangian}

We follow the notations in \S\ref{subsec:tangent space PT}. 
Let us fix a pair of vector spaces $\underline{V}'=(V_1', V_2')$. 
Let $(\xi_1, \xi_2)\subset \underline{V}'$ be a pair of one-dimensional subspaces, and let $\underline{V}:=\underline{V}'/\xi$ be the quotient. 

Consider the following correspondence 
\[
\xymatrix@R=1em@C=0.5em{
	&Z^{st, 0}/(P_1\times P_2)\ar[ld]_{p}\ar[rd]^{q}&\\
	\mathcal{M}(\underline{V}, W)^{st, 0}/\GL_{\underline{V}}&&\mathcal{M}(\underline{V}', W)^{st, 0}/ \GL_{\underline{V}'}
}\]

The subspace $(\xi_1, \xi_2)\subset \underline{V}'$ gives a tautological line bundle on the correspondence $Z^{st, 0}/(P_1\times P_2)$, which will be denoted by $L$. 

Define the  operators 
\begin{align}\label{eqn:EFdef}
&e(g): H_{c, T\times \GL_W}^\ast(\fM_m^0(\underline{V}, W)^{st}, \varphi_{\tilde W_3}\C)^\vee\otimes_{R_T}K_T  \to 
H_{c, T\times \GL_W}^\ast(\fM_m^0(\underline{V}', W)^{st}, \varphi_{\tilde W_3}\C)^\vee\otimes_{R_T}K_T, \\
&f(g): H_{c, T\times \GL_W}^\ast(\fM_m^0(\underline{V}', W)^{st}, \varphi_{\tilde W_3}\C)^\vee\otimes_{R_T}K_T  \to H_{c, T\times \GL_W}^\ast(\fM_m^0(\underline{V}, W)^{st, 0}, \varphi_{\tilde W_3}\C)^\vee\otimes_{R_T}K_T
\end{align}
 by the following convolutions: 
\[
e(g):=q_*(g(c_1(L))\cup p^*), \,\  f(g):=p_*(g(c_1(L))\cup q^*). 
\]
Taking into account  weights and using the canonical isomorphism $\Hom(V_2, V_1)\simeq V_2^\vee\otimes V_1$, the matrix coefficients of  raising operator can be calculated using the Appendix \S~\ref{app:residue} in the following way.
\begin{align*}
&\frac{e\left( {\begin{matrix} t V_2^\vee\otimes  V_1\oplus qV_1^\vee\otimes  V_2\oplus h V_1^\vee\otimes  V_2\oplus\overline{V_2}^\vee\otimes  V_1 \\ \oplus \oplus_{i=1}^m (t^{i-1} W^\vee\otimes  V_2)
\oplus\oplus_{i=1}^{m+1} (t^{1-i} V_1^\vee\otimes  W)-\End(V_2)\end{matrix}}\right)}
{e\left( { \begin{matrix} t P(V_2', V_1')\oplus q P(V_1', V_2')\oplus h P(V_1', V_2')\oplus \overline{V_2}^\vee\otimes V_1'\oplus \xi^\vee\otimes  V_1\\
\oplus \oplus_{i=1}^m (t^{i-1} W^\vee\otimes  V_2')
\oplus\oplus_{i=1}^{m+1} (t^{1-i}V_1^\vee\otimes  W) -\End(V_2') \end{matrix}}
\right)} \\
=&\frac{1} {e\left(t V_2^\vee\otimes  \xi +q  V_1^\vee\otimes  \xi+h  V_1^\vee\otimes  \xi\right)} \frac{e\left(V_1^\vee\otimes  \xi+\xi^\vee\otimes  V_2\right)}{e\left(\xi^\vee\otimes  V_1\right)} \frac{1}{e\left(\oplus_{i=1}^m (t^{i-1} W^\vee\otimes  \xi) \right)}.
\end{align*}
We denote a pyramid partition by $\lambda$, and we denote by $\lambda+\blacksquare$ the one obtained from $\lambda$ by adding a pair \eqref{a pair} denoted by $\blacksquare$. 
Let $\{x_w \mid w=\text{white stone}\}$ be the Chern roots of $V_1$, and 
let $\{x_{b} \mid b=\text{black stone}\}$ be the Chern roots of $V_2$. 
Let $\chi$ be the Chern root of $W$. 
Then the above formula gives
\begin{align*}
\langle\lambda|e_i|\lambda+\blacksquare\rangle=
\Res_{z=\blacksquare}&(-1)^{|\{\text{black only stones}\}|}
z^i\prod_{b\in \text{black stones}} \frac{z-x_{b}}{z-x_{b}-t} \\&\prod_{w \in \text{white stones}} \frac{1}{(z-x_{w}-q)(z-x_{w}-h)}
\prod_{i=1}^{m} \frac{ 1}{z-\chi-(i-1)t }
\end{align*}

Similarly, the lowering operator is given by
\begin{align*}
&\frac{e\left(  { \begin{matrix}t V_2'^\vee\otimes  V_1'\oplus q V_1'^\vee\otimes  V_2'\oplus h V_1'^\vee\otimes  V_2'\oplus\overline{V_2}^\vee\otimes  V_1' \\
\oplus \oplus_{i=1}^m (t^{i-1} W^\vee\otimes  V_2')
\oplus\oplus_{i=1}^{m+1} (t^{1-i}V_1'^\vee\otimes  W)-\End(V_2') 
\end{matrix}}\right)}
{e\left({\begin{matrix}  t P(V_2', V_1')\oplus q P(V_1', V_2')\oplus h P(V_1', V_2')\oplus \overline{V_2}^\vee\otimes V_1'\oplus \xi^\vee\otimes  V_1\\
\oplus \oplus_{i=1}^m (t^{i-1} W^\vee\otimes  V_2')
\oplus\oplus_{i=1}^{m+1} (t^{1-i}V_1^\vee\otimes  W) -\End(V_2')\end{matrix}}\right)}\\
=&\frac{e\left(t\xi^\vee\otimes  V_1 +q \xi^\vee\otimes  V_2+h \xi^\vee\otimes V_2\right)}{e\left(\xi^\vee\otimes  V_1\right)}  e(\oplus_{i=1}^{m+1} (t^{1-i} \xi^\vee\otimes  W) )
\end{align*}
The matrix coefficient is given by 
\begin{align*}
&\langle\lambda+\blacksquare|f_j|\lambda\rangle \\
=&z^j\prod_{b\in \text{black stones} }(z-x_b+q)(z-x_b+h) 
\prod_{w\in \text{white stones} } \frac{z-x_w+t}{z-x_w}
\prod_{i=1}^{m+1} (-z+\chi-(1-i)t)|_{z=\blacksquare}
\end{align*}

Then we  have
\begin{align*}
h(z)=(-1)^{|\{\text{black only stones}\}|}&\prod_{b\in \text{black stones} } \frac{(z-x_b)(z-x_b+q)(z-x_b+h)}{ z-x_b-t} \\
&\prod_{w \in \text{white stones}} \frac{z-x_w+t}{(z-x_w)(z-x_{w}-q)(z-x_{w}-h)}
  \\&\prod_{i=1}^{m} \frac{ 1}{z-\chi-(i-1)t }\prod_{i=1}^{m+1} (-z+\chi-(1-i)t)
\\
=(-1)^{|\{\text{black only stones}\}|}&(-1)^{m+1}(z-\chi-mt)\\
&\prod_{b\in {\text{black stones and white stones}}} \frac{(z-x_b+t)(z-x_b+q)(z-x_b+h)}{(z-x_b-t)(z-x_b-q)(z-x_b-h)}\\
& \cdot \prod_{b\in \text{black stones only} }  \frac{(z-x_b)(z-x_b+q)(z-x_b+h)}{(z-x_b-t)} 
\end{align*}
Similar to the proof in \S~\ref{subsec:C3Yangian}, the operators $e_i, f_j$, and $h(z)$ satisfy the relations of 
$Y_{1}(z_1)$, where $z_1=\chi+mt$. 
This concludes the proof of Theorem~\ref{thm:conifoldYangian}(2).


\section{Shift of the real roots for the resolved orbifold} 
\label{sec:D4 branes}

Let us now move to the discussion of sheaves supported on toric divisors inside of the toric Calabi-Yau $3$-fold $X=X_{m,n}$. As we mentioned in the Introduction, in physics language these divisors correspond to  configurations of D4-branes. We also mentioned that they are expected to give rise  to shifts of real roots of the affine Yangian. The shifts are determined by the intersection numbers of the corresponding divisor with the projective lines $\bbP^1$'s inside $X=X_{m,n}$. As we will see in the examples below this is indeed the case. 

For any  embedded smooth rational curve $\bbP^1\subseteq X$ its formal neighbourhood is isomorphic to the one of the embedded smooth rational curve $\bbP^1$ in either $X_{2,0}$ or $X_{1,1}$. 
For any toric divisor $D\subseteq X$, the formal completion at $\bbP^1$ is isomorphic to the one for the standard toric divisors in either $X_{2,0}$ or $X_{1,1}$. Therefore, we study these cases as the basic building blocks. 

\subsection{Chainsaw quiver}
\label{Chainsaw}
Consider the Calabi Yau 3-fold $T^*\PP^1\times \C=\Tot(\calO_{\PP^1}(-2)\oplus \calO_{\PP^1})$. We take the effective divisor $D$ to be fiber of 
$T^*\PP^1\times \C \to \PP^1$ over the north pole or over the south pole of $\PP^1$. This is illustrated by the  toric diagram \eqref{Intro:N or S pole}.

Let us discuss one of these diagrams (say $D$ is the fiber of 
$T^*\PP^1\times \C \to \PP^1$ over the north pole of $\PP^1$), since the other is similar. 
Note that the $(x,z)$-coordinate hyperplane divisor $\C^2_{xz}\subseteq \C^3_{xyz}$ is invariant under the natural $\Z/2$-action, hence defines a 
divisor in $\C^3/(\Z/2)$. The strict transform under the resolution of singularity $T^*\PP^1\times \C\to \C^3/(\Z/2)$ is the divisor $D$, the fiber of 
$T^*\PP^1\times \C \to \PP^1$ over the north pole of $\PP^1$. Therefore we can describe the moduli space of stable framed sheaves on $T^*\PP^1\times \C$ supported on $D$ such as follows.

Start with the quiver with potential $(Q_3^{fr}, W_{3}^{fr})$, which is quiver description for the
moduli space of framed sheaves on $\C^3$ that supported on the divisor $\C^2_{xz}$. 
Note that the framing $(Q_3^{fr}, W_{3}^{fr})$ has to do with the toric divisor $\C^2_{xz}$, which is different from the framing $(\tilde Q_3, \tilde W_{3} )$ considered in \S~\ref{sec:Hilb}. 
\[
\begin{tikzpicture}[scale=0.9]
    \draw[-, thick] (0,0) -- (2, 0);
    \draw[-, thick] (0,0) -- (0, 1.5);
    \draw[-, thick] (0, 0) -- (-1, -1);
    \filldraw (0,0) node[anchor=north, yshift=-0.5cm, xshift=0.1cm] {$N_{13}$};
    \end{tikzpicture}
 \,\ \,\  \,\ \,\   \,\ \,\  \,\ \,\  
\begin{tikzpicture}
  \draw ($(0,0)$) circle (.08);

     \node at (0.3, 0) {$K$};
       \node at (1.2, 0) {$B_3$};   
       \node at (0, 1.2) {$B_2$};
       \node at (-1.2, 0) {$B_1$};  

    \node at (0, -1.5) {$\square$};    
      \node at (0, -1.8) {$N_{13}$};  
  \draw [<-] (0.1, -0.1) -- (0.1, -1.4);
    \draw [->] (-0.1, -0.1) -- (-0.1, -1.4);
           \node at (0.4, -0.8) {$I_{13}$};  
           \node at (-0.4, -0.8) {$J_{13}$};  
               
   \draw[->,shorten <=7pt, shorten >=7pt] ($(0,0)$)
   .. controls +(40:1.5) and +(-40:1.5) .. ($(0,0)$);
     \draw[->,shorten <=7pt, shorten >=7pt] ($(0,0)$)
   .. controls +(90+40:1.5) and +(90-40:1.5) .. ($(0,0)$);
 \draw[->,shorten <=7pt, shorten >=7pt] ($(0,0)$)
   .. controls +(180+40:1.5) and +(180-40:1.5) .. ($(0,0)$);
              \node at (0, -2.5) {The potential is $W_{3}^{fr}=B_2\big([B_1, B_3]+I_{13}J_{13}\big). $};  
\end{tikzpicture}
\]

Explicitly, the group $\Z/2\Z$ acts on $\C^3$ in such a way that its generator $\sigma\in \Z/2\Z$ is represented  by the matrix
\[
\sigma \mapsto 
\begin{bmatrix}
   -1 & 0 & 0\\
   0     & -1 & 0 \\
   0& 0& 1
\end{bmatrix}
\]
The torus $(\C^*)^2\subset (\C^*)^3$ acts on $\C^3$ by rescaling the coordinates by $t_1,t_2,t_3$ subject to the Calabi-Yau condition $t_1t_2t_3=1$. The above formula gives an embedding of  the  group $\Z/2\Z$ into  the subtorus $(\C^*)^2$. 
Under the $\Z/2\Z$ action the weights of the arrows of $Q_3^{fr}$ are given by 
\begin{align*}
&B_1 \rightsquigarrow -1, \,\ B_2 \rightsquigarrow -1, \,\ B_3 \rightsquigarrow 1\\ 
&I_{13}  \rightsquigarrow +1, \,\ J_{13}  \rightsquigarrow -1.
\end{align*}
We decompose the vector spaces $K,  N_{13}$ into eigenspaces with eigenvalues $1, -1$. 
\begin{align*}
&K=K^0\oplus K^1\\ 
&N_{13}=N_{13}^0\oplus N_{13}^1
\end{align*}
Setting $N_{13}^1$ to be zero, we obtain the 
 framed quiver $Q$ from \eqref{eqn:orbifold_chainsaw} (with the square node removed).\Omit{
\begin{equation}\label{eqn:orbifold_chainsaw}
\begin{tikzpicture}[scale=0.8]
        	\node at (-2.2, 0) {$\bullet$};   	\node at (-2.2, 0.5){$V_0$}; \node at (-3.5, 0){$B_3$};
	\node at (2.2, 0) {$\bullet$};  \node at (2.2, 0.5){$V_1$};\node at (3.5, 0){$\tilde{B}_3$};
		\draw[-latex,  bend left=30, thick] (-2, 0.1) to node[above]{$B_2$} (2, 0.1);
		\draw[-latex,  bend left=50, thick] (-2, 0.3) to node[above]{$B_1$} (2, 0.3);
	\draw[-latex,  bend left=30, thick] (2, -0.1) to node[above]{$\tilde{B}_1$} (-2, -0.1);
		\draw[-latex,  bend left=50, thick] (2, -0.3) to node[above]{$\tilde{B}_2$} (-2, -0.3);
				\path (2.2, 0) edge [loop right, min distance=2cm, thick, bend right=40 ] node {} (2.2, 0);
				\path (-2.2, 0) edge [loop left, min distance=2cm, thick, bend left=40 ] node {} (-2.2, 0);
\node at (0, -4) {$\square$};  \node at (0.5, -4) {$W$};  
 \node at (-1.5, -2) {$I_{13}$};   \node at (1.5, -2) {$J_{13}$};  
 \draw[->, thick] (0, -3.8) -- (-2.2, -0.2) ;
 \draw[->, thick]  (2.2, -0.2) --(0, -3.8) ;
    \end{tikzpicture}
\end{equation}   }
with potential $W$ given by
\[
W=-(B_2 \tilde{B}_1\tilde{B}_3-\tilde{B}_2\tilde{B}_3B_1+\tilde{B}_2B_1 B_3-B_2B_3\tilde{B}_1)+B_2I_{13} J_{13}
\]
Framed stability condition implies that the image of $I_{13}$ is invariant under $B_1,B_3$ and $\tilde{B}_1,\tilde{B}_3$ and generates the whole space.

We now use the dimensional reduction of a quiver with potential and a cut (see \cite[Section 4.8]{KS}, \cite[Appendix]{D}). We take the cut to be the set of arrows consisting  of $B_2, \tilde{B_2}$. Below we present the dimensionally reduced quiver with relations in our case, 
referring the reader to \cite[Appendix]{YZ4} for a description of the COHA for the dimensionally reduced quiver with relations for a general cut. 

Let $Q':=Q \setminus \{B_2, \tilde{B_2}\}$ be the quiver obtained by removing the two arrows $B_2, \tilde{B}_2$. A picture of the quiver obtained is \eqref{eqn:chainsaw}. Then, $\Rep(Q, V_0, V_1)=\Rep(Q', V_0, V_1)\times (\Hom(V_0, V_1)\oplus \Hom(V_1, V_0))$. 
We have the following diagram
\[
\xymatrix{
Z\times \mathbf{A}^n \ar@{^{(}->}[r]  \ar[d]&\Rep(Q, V_0, V_1)  \ar[d]^{\pi}\ar[rd]&\\
Z \ar@{^{(}->}[r] & \Rep(Q', V_0, V_1) \ar[r]& \pt
}
\]
where the affine space $\mathbf{A}^n$ is $\Hom(V_0, V_1)\oplus \Hom(V_1, V_0)$, and $Z$ is  the algebraic variety
\begin{align*}
Z=&\{z\in  \Rep(Q', V_0, V_1)\mid \tr(W)(z, l)=0, \text{for all $l\in \Hom(V_0, V_1)\oplus \Hom(V_1, V_0)$}\}\\
=&\{z\in  \Rep(Q', V_0, V_1)\mid \frac{\partial W}{\partial B_2}(z)=0, \frac{\partial W}{\partial \tilde{B}_2}(z)=0\}\\
=& \Rep(\C[Q']/J, V_0, V_1), 
\end{align*}
where $J$ is the ideal generated by
\begin{equation}\label{rel:I}
\frac{\partial W}{\partial B_2}=-(\tilde{B}_1\tilde{B}_3-B_3\tilde{B}_1)+I_{13} J_{13}, \,\ 
\frac{\partial W}{\partial \tilde B_2}=-(-\tilde{B}_3B_1+B_1 B_3).
\end{equation}
Then \cite[Section 4.8]{KS} and \cite[Appendix]{D} give an isomorphism
\[
H^*_{c, \GL(V_0)\times  \GL(V_1)\times T\times \GL(N_{13})}(\Rep(Q, V_0, V_1), \varphi_{tr W})\cong H^*_{c, \GL(V_0)\times  \GL(V_1)\times T\times \GL(N_{13})}(Z\times \mathbf{A}^n,  \C), 
\]
where $\GL(N_{13})$ is the flavor group acting on $\Rep(Q, V_0, V_1)$. 
Note that the variety $Z$ is the representation variety of the quiver $Q'$ with relations on the vector spaces $(V_0,V_1)$.
The quiver $Q' $ with relation \eqref{rel:I} is known as a {\it chainsaw quiver} (defined in \cite{FR}).

The following definition of Yangian of $\mathfrak{gl}(n)$ can be found in \cite{N12} and \cite{FPT}. 

\begin{defn}
The Yangian of $\mathfrak{gl}(n)$ is generated by 
\[
\{e_{i}^{(r)}, f_{i}^{(r)}, g_{j}^{(r)} \mid 1\leq i < n, 1\leq j \leq n, r \geq 1\}. 
\]  
Define the following generating series of the generators
\[
e_{i}(z):=\sum_{r\geq 1} e_{i}^{(r)} z^{-r}, 
f_{i}(z):=\sum_{r\geq 1} f_{i}^{(r)} z^{-r},  
g_{j}(z):=1+\sum_{r\geq 1} g_j^{(r)}z^{-r}. 
\]
They are subject to the following relations \cite[Lemma 2.48]{FPT}. 
\begin{align}
&[g_i(z), g_j(w)]=0; \label{rel:1}\\
& (z-w)[g_i(z), e_j(w)]=(\delta_{i, j}-\delta_{i, j+1}) g_i(z) (e_j(z)-e_j(w));\\
&(z-w)[g_i(z), f_j(w)]=(\delta_{i, j+1}-\delta_{i, j}) (f_j(z)-f_j(w))  g_i(z) ;\\
&[e_{i}(z), f_{j}(w)]=0, \text{if $i\neq j$};\\
&(z-w) [e_i(z), f_i(w)]=\frac{g_{i+1}(w)}{g_i(w)} -\frac{g_{i+1}(z)}{g_i(z) }; \label{rel:2}\\
&(z-w) [e_i(z), e_i(w)]=-(e_i(z)-e_i(w))^2;\\
&(z-w) [e_i(z), e_{i+1}(w)]=-e_{i}(z) e_{i+1}(w)+e_{i}(w)e_{i+1}(w)-[e_{i+1}^{(1)}, e_i(w)]+[e_{i+1}^{(1)}, e_i(z)];\\
& [e_i(z), e_{j}(w)]=0, \text{if $|i-j|>1$};\\
& [e_i(z_1), [e_i(z_2), e_j(w)]]+[e_i(z_2), [e_i(z_1), e_j(w)]]=0, \text{if $|i-j|=1$};\\
&(z-w) [f_i(z), f_i(w)]=(f_i(z)-f_i(w))^2;\\
&(z-w) [f_i(z), f_{i+1}(w)]=f_{i+1}(w) f_{i}(z)-f_{i+1}(w)f_{i}(w)+[ f_i(w), f_{i+1}^{(1)}]-[f_i(z), f_{i+1}^{(1)}];\\
& [f_i(z), f_{j}(w)]=0, \text{if $|i-j|>1$};\\
& [f_i(z_1), [f_i(z_2), f_j(w)]]+[f_i(z_2), [f_i(z_1), f_j(w)]]=0, \text{if $|i-j|=1$}.  \label{rel:3}
\end{align}
\end{defn}
Let $\bold{r}=(r_1, \cdots, r_n)\in \N^n$,  we define the $\bold{r}$-shifted Yangian of $\widehat{\mathfrak{gl}(n)}$ to be the algebra generated by 
\[
\{e_{\bar{i}}^{(r)}, f_{\bar{i}}^{(r)}, g_{\bar{j}}^{(r)} \mid \bar{i}, \bar{j}\in \Z/n\Z, r \geq 1\}, 
\]  
subject to the same relations \eqref{rel:1}-\eqref{rel:3}, except the relation \eqref{rel:2} is modified to be
\begin{equation}
(z-w) [e_{\bar{i}}(z), f_{\bar i}(w)]=\frac{g_{\bar{i+1}}(w)}{g_{\bar{i}}(w)} -(z^{r_{i+1}-r_i}) \frac{g_{\bar{i+1}}(z)}{g_{\bar{i}}(z) }.
\end{equation}

We have identified the cohomology $H^*_{c, \GL(V_0)\times  \GL(V_1)\times T \times \GL(N_{13})}(\Rep(Q, V_0, V_1)^{st}, \varphi_{tr W}\C)$ with the cohomology of the chainsaw quiver
$H^*_{c, \GL(V_0)\times  \GL(V_1)\times T \times \GL(N_{13})}(Z^{st}\times \mathbf{A}^n,  \C)$ via the dimensional reduction. 
It is shown in \cite[Theorem 4.3]{N} that the shifted quantum toroidal algebra of type $A$ acts on the K-theory of the chainsaw quiver variety.   
Following the same calculation, with $K$-theory replaced by cohomology, we obtain
\begin{thm}
The $(1, 0)$-shifted Yangian of  $\widehat{\mathfrak{gl}(2)}$ acts on 
\[
\bigoplus _{V_0, V_1} H^*_{c, \GL(V_0)\times  \GL(V_1)\times T \times \GL(N_{13})}(\Rep(Q, V_0, V_1)^{st}, \varphi_{tr W}\C)^\vee \otimes_{R_T} K_T. 
\]
\end{thm}

\subsection{Few more examples}
\subsubsection{Nakajima quiver varieties}
Consider the CY 3-fold $T^*\PP^1\times \C=\Tot(\calO_{\PP^1}(-2)\oplus \calO_{\PP^1})$, and the effective divisor $\mathcal{O}(-2)$ on $\PP^1$. 
This is illustrated in the following toric diagram. 
\[\begin{tikzpicture}[scale=0.9]
    \draw[-, thick] (0,1) -- (2, 1) ;
    \draw[-, thick] (0,0) -- (2, 0);
    \draw[-, thick] (0,0) -- (0, 1);
    \draw[-, thick] (0, 1) -- (-1, 2);
    \draw[-, thick] (0, 0) -- (-1, -1);
         \filldraw (0,0.5) node[anchor=west, xshift=-0.8cm] {$1$};
    \end{tikzpicture}
    \]
 The quiver description is given by the following figure:
 
    $$ \begin{tikzpicture}[scale=0.8]
        	\node at (-2.2, 0) {$\bullet$};   	\node at (-3.5, 0) {$B_3$};   	
	\node at (2.2, 0) {$\bullet$};  \node at (3.5, 0) {$\tilde{B}_3$};  
		\draw[-latex,  bend left=30, thick] (-2, 0.1) to node[above]{$B_2$} (2, 0.1);
		\draw[-latex,  bend left=50, thick] (-2, 0.3) to node[above]{$B_1$} (2, 0.3);
	\draw[-latex,  bend left=30, thick] (2, -0.1) to node[above]{$\tilde{B}_1$} (-2, -0.1);
		\draw[-latex,  bend left=50, thick] (2, -0.3) to node[above]{$\tilde{B}_2$} (-2, -0.3);
				\path (2.2, 0) edge [loop right, min distance=2cm, thick, bend right=40 ] node {} (2.2, 0);
				\path (-2.2, 0) edge [loop left, min distance=2cm, thick, bend left=40 ] node {} (-2.2, 0);
\node at (-2.2, -2) {$\square$};  
 \draw[->, thick] (-2.1, -1.8) -- (-2.1, -0.2) ; \node at (-1.7, -1) {$I_{12}$};   
 \draw[->, thick]  (-2.3, -0.2) --(-2.3, -1.8) ; \node at (-2.7, -1) {$J_{12}$};  
    \end{tikzpicture}
 $$  
The potential is given by the formula:
\[
B_2 \tilde{B}_1\tilde{B}_3-\tilde{B}_2\tilde{B}_3B_1+\tilde{B}_2B_1 B_3-B_2B_3\tilde{B}_1+B_3 I_{12}J_{12}. 
\]

The framed stability condition says that any subspace containing the image of $I_{12}$ and invariant  under the action of $B_1,B_2,\tilde{B}_1,\tilde{B}_2$ coincides with the whole  space (i.e. the image of $1$ is a cyclic vector). Then 
using the dimension reduction with respect to the loops, the corresponding moduli spaces are the Nakajima quiver varieties with framing $1$. It is known that there is an action of the affine Yangian of ${\mathfrak{sl}}_2$ (without shifting).

\subsubsection{Blowup of $\C^2$}
Let $X$ be the resolved conifold $\Tot(\calO_{\PP^1}(-1)\oplus \calO_{\PP^1}(-1))$. As the effective divisor we take the total space of one of the line bundle summands over $\PP^1$. We denote them by  $\calO(-1)_1$ and $\calO(-1)_2$ respectively. This is illustrated in the following toric diagram. 
    \[
\begin{tikzpicture}[scale=0.9]
    \draw[-, thick] (0,0) -- (2, 0);
    \draw[-, thick] (0,0) -- (-1, 1) ;
     \draw[-, thick] (0,0) -- (-1, -1) ;
         \draw[-, thick] (2,0) -- (3, 1) ;
     \draw[-, thick] (2,0) -- (3, -1) ;
         \filldraw (1,0) node[anchor=north, yshift=0.8cm] {$1$};
    \end{tikzpicture}
    \]
    
The corresponding quiver is
\begin{equation*}
\begin{tikzpicture}[scale=0.8]
        	\node at (-2.2, 0) {$\bullet$};   		
	\node at (2.2, 0) {$\bullet$};  
		\draw[-latex,  bend left=30, thick] (-2, 0.1) to node[above]{$a_2$} (2, 0.1);
		\draw[-latex,  bend left=50, thick] (-2, 0.3) to node[above]{$a_1$} (2, 0.3);
	\draw[-latex,  bend left=30, thick] (2, -0.1) to node[above]{$b_1$} (-2, -0.1);
		\draw[-latex,  bend left=50, thick] (2, -0.3) to node[above]{$b_2$} (-2, -0.3);				
\node at (0, -4) {$\square$};  
 \node at (-1.5, -2) {$i$};   \node at (1.5, -2) {$j$};  
 \draw[->, thick] (0, -3.8) -- (-2.2, -0.2) ;
 \draw[->, thick]  (2.2, -0.2) --(0, -3.8) ;
    \end{tikzpicture}
\end{equation*}   
with the potential 
\[
a_1b_1a_2 b_2-a_1b_2a_2b_1+b_1ij. 
\]
Applying the dimensional reduction to the above quiver with potential with the cut  consisting of  the edge $b_1$, we obtain the quiver with relations from
\cite{NY}, which described stable framed  sheaves on the blowup of $\C^2$. 
We postpone the study of shifted Yangian action to the future.

\section{What to expect for general toric Calabi-Yau $3$-folds}
\label{sec:D4 branes2}
\label{subsec:Expectation}

 \subsection{Moduli spaces and shifted Yangians}\label{subsec:shiftedYangian_expectation}
 In this subsection we explain how to put examples from Sections \ref{sec:Hilb}, \ref{sec:PTconifold},  \ref{sec:D4 branes} into a more general framework.

Fix a quiver with potential $(Q, W)$. Assume $Q$ is symmetric.  (It follows from this assumption that the quadratic form $\chi_{R}|_{\Gamma_l}$ defined in \cite{KS} is symmetric.) In \cite[Section 6.2, Question 6.2 ]{KS}, Kontsevich-Soibelman conjecture that there exits a ${\bf Z}$-graded Lie algebra $\g$ (BPS Lie algebra), such that the (non-equivariant) COHA of $(Q,W)$ is isomorphic to the universal enveloping algebra of the current algebra $\g\otimes \C[\mathbb{T}]$ of $\g$. This conjecture was proved later by Davison and Meinhardt \cite{DM}. In the special case when $(Q, W)$ is the ``tripled" quiver of a simply-laced Kac--Moody Lie algebra as in \cite{Ginz}, it follows from the recent paper \cite{D2} of Davison that the zeroth piece of the perverse filtration on the COHA of $(Q, W)$ is isomorphic to the universal enveloping algebra of $\g$. The latter  contains the upper-triangular subalgebra of the Kac--Moody Lie algebra. As a consequence, if the tripled quiver $Q$ has exactly one  loop at each vertex, then the $(i,j)$-th entry of the Cartan matrix of the Kac-Moody Lie algebra is equal to the number of arrows of $Q$ between $i$ and $j$ if $i\neq j$. The diagonal entries are all equal to $2$. 
  
We do not know an explicit description of the root system of the BPS Lie algebra for general symmetric quivers with potential. 
Nevertheless, assuming furthermore that the potential $W$ is homogeneous,  the three examples ~\ref{ex:m_0}~\ref{ex:m_n} and ~\ref{ex:nhat}, 
 suggest that the root system of the BPS Lie algebra should contain a sub root system, the Cartan matrix of which is given as
 \begin{equation}
 \label{super Cartan}
\langle \alpha_i, \alpha_j \rangle = \begin{cases}
-\#\{h: i\to j, h\in H\} &\text{if $i\neq j$; }\\
2 & \text{if  $i=j$, and $i$ has one  loop; }\\
0 & \text{if  $i=j$, and $i$ has no loops. }
\end{cases}
\end{equation}
This sub root system is that of a Lie subalgebra of the BPS Lie algebra, corresponding  to the spherical subalgebra of the COHA. 
The Drinfeld double of the entire COHA is expected to be a Cartan doubled Yangian of the root system of the entire BPS Lie algebra.

Now go back to geometry. Let $X$ be a toric CY 3-fold as in \S~\ref{intro}. We are in the setting of \S~\ref{sec:toricCY}, that is, $f:X\to Y$ is a resolution with $Y$ affine and the  fibers of $f$ are at most 1-dimensional. As mentioned in the introduction of \cite[pg.2]{Nagao} that such affine toric CY 3 fold $Y$ can be classified in terms of the lattice polygon in $\R^2$. The classification is given by the family  $Y=Y_{m, n}$, together with finitely many exceptional cases.  When $Y=Y_{m, n}$, the corresponding quiver is symmetric \cite[Section 1.2]{Nagao}.

Let $\Modf\calA_0\subseteq \Mod\calA_0$ be the full subcategory consisting of coherent sheaves of $\calA_0$-modules whose  set-theoretical supports  as coherent sheaves on $Y$ are zero-dimensional.  Similarly, we have the abelian category $\mathscr{A}_f$, which is equivalent to $\Modf\calA_0$ under the functor $R\Hom_X(\calP,-)
$.

Let $K_0(\mathscr{A})_\Z$ be the Grothendieck group of $\mathscr{A}_f$. 
Let  $\{S_i\}_{i\in I}$ be a collection of pairwise distinct simple objects in $\mathscr{A}_f$.
For each $S_i$, let $[S_i]\in K_0(\mathscr{A})_{\Z}$ be its  class in the Grothendieck group, which we call a {\it simple root}. Using this terminology, the set $\{[S_i]\}_{i\in I}$ of simple roots is  an integer basis of $K_0(\mathscr{A})_{\Z}$. Thus the latter can be thought of as a {\it root lattice}. 
We call a simple root $[S_i]$ is {\it bosonic} if the simple object $S_i$  in $D^b\Coh(X)$ is bosonic, that is $\Ext^*(S_i, S_i)=H^{2*}(\PP^3)$. 
We call a simple root $[S_i]$ is {\it fermionic} if the simple object $S_i$ is fermionic,  that is, $\Ext^*(S_i, S_i)=H^*(S^3)$.  
This terminology is explained in \S~\ref{sec:toricCY}.

 In the three examples ~\ref{ex:m_0}~\ref{ex:m_n} and ~\ref{ex:nhat}, the Cartan matrix \eqref{super Cartan} can be alternatively described as the following pairing on $K_0(\mathscr{A})_{\Z}$
\begin{equation*}
\langle [S_i], [S_j] \rangle = \begin{cases}
-\dim( \Ext^1(S_i, S_j)) &\text{if $i\neq j$; }\\
2 & \text{if  $i=j$, and $i$ is bosonic; }\\
0 & \text{if  $i=j$, and $i$ is fermionic. }
\end{cases}
\end{equation*}
Indeed, in all these cases, the algebra $\calA_0$ has a grading coming from the $\C^*$-action on $X$ which contracts $X$ to its special fiber, so that the tilting bundle is $\C^*$-equivariant. The quiver with potential $(Q,W)$, whose Jacobian algebra is $\calA_0$, can be chosen so that the vertices are labeled by $I$, and  the number of arrows  from $i$ to $j$ is $\dim \Ext^1(S_i,S_j)$ for $i,j\in I$ \cite[Theorem~3.1 and proof]{B}. Hence, the matrix \eqref{super Cartan} can be expressed in this way. 
In particular, this gives rise to  a free abelian group $K_0(\mathscr{A})_{\Z}$ endowed with a basis and an integer pairing. 
In general, this resembles a part of some ``Cartan data", possibly that of a sub root system of the BPS Lie algebra.

\begin{exa}
In Example~\ref{ex:m_n},
reading off the sequence of bosonic and fermionic simple objects by following the diagram \eqref{fig2} from the bottom to the top (or equivalently from the top to the bottom that produces a reflected Dynkin diagram), one recovers different root systems of the  $\mathfrak{gl}(m|n)$ algebra \cite{FSS}. The above-defined root system thus coincides with the root system of  $\widehat{\mathfrak{gl}}(m|n)$. This configuration has been considered by various authors including \cite{LY,R,Ue}.

 In Example~\ref{ex:nhat}, the generalized root system is expected to be of ``double affine" type $A_n$. The construction of \S~\ref{sec:doubleCoHA} is expected to give a definition of double affine Yangian of type $A_n$ \cite{Costello2}.

\end{exa}

The moduli space of perverse coherent system of $X$ in the sense of \S~\ref{subsec:modulipervcoh} will depend on a choice of the stability parameter $\zeta$, as well as the algebraic cycle $\chi$. 
Here the stability parameter is understood in the same sense as in \S~\ref{sec:PTconifold}. We denote by $\mathbb X$ the space of stability parameters of the perverse coherent system. It can be identified with $K_0(D^b \Coh(X))_\R^\vee$. We say that $\zeta \in \mathbb X$ is generic if it lies in the interior of a chamber. Each generic $\zeta$ determines an abelian subcategory $\mathscr{A}_\zeta$, which depends only on the chamber containing $\zeta$. In the case when $X=X_{m,n}$ this abelian subcategory comes from a tilting bundle explicitly constructed in \cite[\S~1]{Nagao}.

The above-mentioned parameter  $\chi$ is an algebraic cycle of the form $N_0 X+\sum_{i=1}^r N_iD_i$ where $N_i\in\mathbf N$,  and each $D_i$ is a toric divisor. Taking its class in the Borel-Moore homology with integral coefficients, we get 
$[\chi]\in H^{BM}_*(X,\Z)$. As we only consider algebraic cycles, we have $[\chi]\in H^{BM}_6(X,\Z)\oplus H^{BM}_4(X,\Z)$.

We can think of  $[\chi]$ as a coweight (i.e. an integer functional on the root lattice) in the following way. The 0-th homology of $X$ is one dimensional, let $[\pt]$ be its basis. 
Let $\{[C_1], \cdots, [C_{m+n-1]}\}$ be the classes of  curves $C_i\simeq \PP^1$ which give a basis of $H_2(X,{\bf Z})$.
Let 
\[
ch: K_0(\mathscr{A}_\zeta)\to H_*(X)
\]
be the homological Chern character map \cite[Section 5.9]{CG}. 
Each  simple object $S_i \in K_0(\mathscr{A}_\zeta)$ gives rise to a class $ch([S_i])\in H_*(X)$.  
The shift is then determined by the pairing
\begin{equation}\label{eq:pairing}
\langle \chi, ch[S_i]\rangle:=-[\chi] \cap ch[S_i] , \text{for all $i\in I$}, 
\end{equation}
where $\cap: H_{6-*}^{\BM}(X)\otimes H_*(X)   \to \Z$ is the intersection pairing \cite[Proposition 2.6.18]{CG}. 
Since the set of simple objects $\{S_i\mid i\in I\}$ gives rise to the set of simple roots spanning the root lattice, 
the homology class $\chi$ gives rise to a coweight via the pairing \eqref{eq:pairing}.

Let us illustrate the above discussion in the  following three examples. 
\begin{enumerate}
\item 
Consider the PT moduli space of $X_{1, 1}$. 
There are two simple objects 
\[
\{S_0, S_1\}=\{\calO_{\PP^1}(m), \calO_{\PP^1}(m+1)[-1]\}
\]

They generate the abelian category (heart of the corresponding $t$-structure) $\mathscr{A}_m^+$. Their classes in  $K_0( \mathscr{A}_m^+)$ gives  the two simple roots. 
Consider the short exact sequence $
0\to \calO\to \calO(m) \to \C_0^{\oplus m} \to 0$ of sheaves on $\PP^1$, 
where $\C_0$ is the skyscraper sheaf at $0\in\PP^1$.  Thus, in $K_0( \mathscr{A}_m^+)$, we have $
[\calO_{\PP^1}(m)]=[\calO]+m[\C_0].$
We calculate the homological Chern character map using the devissage principle \cite[Proposition~5.9.3]{CG}, we have 
\begin{align*}
&ch[S_0]=ch[\calO_{\PP^1}(m)]=[C]+m[\pt],\\ 
&ch[S_1]=ch[\calO_{\PP^1}(m+1)[-1]]=-[C]-(m+1)[\pt], 
\end{align*} where $[C]$ is the $\PP^1$ class in $H_2(X)$. 
In this case $\chi=[X_{1, 1}]$. The intersection pairing is given by $[X_{1, 1}] \cap C=0, [X_{1, 1}] \cap [\pt]=1$. 
Note that the imaginary root $\delta=ch[S_0]+ch[S_1]$.  This implies that 
\[
\langle \chi,  \delta \rangle=-\chi\cap \delta=1. 
\]
Thus, we have the $+1$ shift of the imaginary root, cf. in \S\ref{sec:PTconifold}. 

\item Consider $X=X_{2, 0}$, and $D$ be the fiber of the vector bundle $X_{2, 0} \to \PP^1$ over either  the north pole or  the south pole of $\PP^1$, cf. \S\ref{Chainsaw}. 
There are two simple objects 
\[
\{S_0, S_1\}=\{\calO_{\PP^1}(-m), \calO_{\PP^1}(-m-1)[-1]\}
\] in $K_0( \mathscr{A}_m)$, which corresponds to the two simple roots. 
Consider the short exact sequence $
0\to \calO(-m)\to \calO \to \C_0^{\oplus m} \to 0$, 
where $\C_0$ is the skyscraper sheaf at $0\in\PP^1$. 
Thus, in $K_0( \mathscr{A}^-_m)$, we have $
[\calO_{\PP^1}(-m)]=[\calO]-m[\C_0].$
Applying the homological Chern character map, we have 
\begin{align*}
&ch[S_0]=ch[\calO_{\PP^1}(-m)]=[C]-m[\pt],\\ 
&ch[S_1]=ch[\calO_{\PP^1}(-m-1)[-1]]=-[C]+(m+1)[\pt], 
\end{align*}
where $[C]$ is the $\PP^1$ class in $H_2(X)$. 
In this case $\chi=[D]$. The intersection pairing is given by 
\[
[D] \cap [C]=1, [D] \cap [\pt]=0. 
\] This implies that 
\[
\langle  \chi, ch[S_0] \rangle=- \chi\cap ch[S_0]=-1, 
\langle  \chi, ch[S_1] \rangle =-\chi\cap ch[S_1]=1. 
\]
Thus, the shifts of the two simple roots are $-1$, $1$ respectively, cf. \S\ref{Chainsaw}. 
\item
Consider the DT moduli space of $\C^3$. In this case $\chi=[\C^3]$. The imaginary root $\delta$ is given by the class $[\pt]$. 
The intersection pairing is given by 
\[
\langle  \chi, \delta \rangle=-[\C^3] \cap \delta=-1. 
\]
Thus we have the $-1$ shift of the imaginary root, cf. \S\ref{sec:Hilb}. 
\end{enumerate}

\subsection{Cartan doubled Yangian from geometry}
\label{subsec:D6}
\label{sec:cartanDoubled}

Let $D^b \Coh(X)$ be the bounded derived category of coherent sheaves on $X$, and $D^{\Z_2} \Coh(X)$ be the derived category of 
$2$-periodic complexes of coherent sheaves on $X$ (see \cite{Br2} for the definition). 
Then we have a functor 
\[
F: D^b \Coh(X)\to D^{\Z_2} \Coh(X)
\]
by taking direct sum of odd and even degree complexes. 
We expect that there is a geometric construction of the Cartan doubled Yangian $Y_{\infty}$,  constructed by taking cohomology of certain moduli stack of objects in $D^{\Z_2} \Coh(X)$ endowed with Hall multiplication, similarly to \cite{Br2}. 

Let $\mathscr{A}\subset D^b \Coh(X)$ be the heart of $D^b \Coh(X)$ corresponding to a tilting bundle as in \S~\ref{intro}. Denote by $\calH_{\mathscr{A}}$ the COHA associated to $\mathscr{A}_f$ constructed as in \S~\ref{intro}. 
We expect the choice of the  t-structure defines an algebra embedding $\calH_{\mathscr{A}} \to  Y_{\infty}$, as well as a triangular decomposition  \[Y_{\infty}\cong\calH_{\mathscr{A}}\otimes\calH_0\otimes\calH_{\mathscr{A}[1]},\]
where the subalgebra $\calH_0$ is a polynomial algebra $\C[\hbar_1,\hbar_2][\psi_{i,k}]_{i\in I, k\in \Z}$, with infinite variables labeled by $I \times \Z$. 

For each $i\in I$, let $\langle S_i \rangle\subseteq \mathscr{A}$ be the Serre subcategory of $\mathscr{A}$ generated by $S_i$. 
Then, the Cartan doubled Yangian $Y_{\infty}^{\alpha_i}$, associated to  $D^{\Z_2}(\langle S_i \rangle)$,  is a subalgebra of $Y_{\infty}$.

Let $\delta$ be the imaginary root defined in the same way as in \cite[Chapter 5]{Kac}. That is, 
$\delta$ is a root, but $\delta$ is not in the Weyl group orbit of the simple roots. We expect that a quotient of $Y_{\infty}$ is isomorphic to $\affY$, which corresponds to the root $\delta$ of $Y_{\infty}$.

We are interested in representations of $Y_{\infty}$ on the cohomology of the moduli spaces of perverse coherent systems on $X$. In the case $X=X_{m,n}$, everything can be spelled out in the language of quivers with potential using \cite{Nagao}.  
Let $\fM_\zeta^\chi$ be the moduli space of stable perverse coherent systems. There is a potential function $W$ defining the symmetric obstruction theory of $\fM_\zeta^\chi$, which carries an action of $T\times G_{fr}$. Let us consider the vector space $V:=H^*_{c, T\times G_{fr}}(\fM_\zeta^\chi,\varphi_W\C)^\vee\otimes_{R_T} K_T$. 
Then we expect the following.

\begin{conj}
 The algebra $Y_{\infty}$ is isomorphic to $D(\mathcal{SH})$ constructed in \S~\ref{sec:doubleCoHA}. There is an action of $Y_\infty$ on $V$ in agreement with the general philosophy of \cite{So}. Moreover, this action factors through an action of the shifted affine Yangian, with the shift determined by $[\chi]$ and $\zeta$ as above.
\end{conj}

\subsection{A braid group action}

Starting from the heart of the $t$-structure coming from a tilting bundle as above, one can construct the braid group action  on $D^b\Coh(X)$  induced by mutations with respect to subcategories $\langle S_i\rangle, i\in I$. We expect that on one hand this braid group action induces an action on $D^{\Z_2}\Coh(X)$ and hence on  the algebra $Y_\infty$. On the other hand it induces an action on $K_0(D^b\Coh(X))$, and therefore an action on $\mathbb X$. We expect that this action gives rise to an action of the Weyl group of the generalized root system described in \S~\ref{subsec:shiftedYangian_expectation}. 
The set of roots in $K_0(D^b\Coh(X))$ determines a hyperplane arrangement in $\mathbb X$. The braid group  acts by reflections with respect to these hyperplanes. 

In particular, for two adjacent chambers separated by a root hyperplane, the reflection with respect to the root hyperplane sends the heart of the $t$-structure associated to one chamber to that of the other one. These two chambers determine two Borel subalgebras of $Y_\infty$ which differ by the action of a Weyl group element. 

As  an illustration we explicitly calculate below the effect of the action of the affine Weyl group on Iwahori subgroups in the well-known case of $\widehat{\fs\fl_2}$. 
Then we discuss a similar picture in an example of a particular Calabi-Yau $3$-fold in Remark~\ref{rmk:AffineGrassmannian} 

\begin{exa}\label{ex:sl2_flag}
Let $G:=\SL_2(\C)$ and $T$ be its maximal torus of diagonal matrices. 
Denote by $B$ the Borel subgroup consisting of upper triangular matrices, and 
$B^-$ the one consisting of lower triangular matrices.

Let $G((t)):=SL_2((t))$, $U^-((t)):=\{\begin{bmatrix} 1&0\\ * & 1 \end{bmatrix}\}\subset G((t))$ be the corresponding groups over the field of Laurent series ${\bf C}((t))$. 
The Coxeter presentation of the affine Weyl group is $\widehat{W}=\{s_0, s_1\mid s_0^2=s_1^2=1\}$. It is isomorphic to $\Z_2\ltimes\Z$, where $\Z_2$ is generated by $s_1$ and $\Z$ is generated by $s_0s_1$.
In particular, any element can be uniquely written as $w=(s_0s_1)^m$ or $w=s_1(s_0s_1)^m $ for some $m\in\Z$.
In $G((t))$, we have a subgroup, which is the Tits extension of $\widehat{W}$, with generators $n_{s_0}, n_{s_1}$ defined  as
\[
n_{s_0}:=
\begin{bmatrix}
0& -t^{-1}\\
t &0
\end{bmatrix}, 
n_{s_1}:=
\begin{bmatrix}
0& 1\\
-1 &0
\end{bmatrix}.  \,\ \text{Thus, }
(n_{s_0}n_{s_1})^m=
\begin{bmatrix}
t^{-m}&0\\
0& t^m
\end{bmatrix}
\]
For any Weyl group element $w=(s_0s_1)^m$ (resp. $w=s_1(s_0s_1)^m $), we write the corresponding element in $G((t))$ as $w=(n_{s_0}n_{s_1})^m $
(resp. $w=n_{s_1}(n_{s_0}n_{s_1})^m $) slightly abusing the notation. 

We have the following three subgroups of $G((t))$,
 \begin{align*}
\text{the positive Iwahori }&I^+=\{g \in G[[t]] \mid g(0)\in B \},\\
\text{the level 0 Iwahori }&I^0=T[[t]] U^-((t)),\\
\text{the negative Iwahori }&I^-=\{g \in G[[t^{-1}]] \mid g(\infty)\in B^- \}. 
\end{align*}
Following \cite{MRY}, for $c\in \C, k\in \Z$, we define
\begin{align*}
X_{\alpha_1+k\delta}(c)=
\begin{bmatrix}
1& ct^k\\
0& 1
\end{bmatrix},  \,\ 
X_{-\alpha_1+k\delta}(c)=
\begin{bmatrix}
1& 0\\
ct^k & 1
\end{bmatrix}. 
\end{align*}
The following points are in $I^+$: 
\[
X_{-\alpha_1+\delta}(1)=
\begin{bmatrix}
1& 0\\
t &1
\end{bmatrix} ,    
X_{\alpha_1}(1)=\begin{bmatrix}
1& 1\\
0&1
\end{bmatrix} 
\]
and they determine $I^+$ in the following sense. We take $X_{-\alpha_1+\delta}(1)$ and $X_{\alpha_1}(1)$ to be the set of simple roots. 
Then, any positive root will be $X_{\alpha+k\delta}(c)$ if $\alpha=\alpha_1$ and $k\in \Z_{\geq 0}$ or if $\alpha=-\alpha_1$ and $k\in  \Z_{> 0}$. 
We have 
\begin{align*}
I^+&=\{ a=\begin{bmatrix} 
a_1 & a_2\\
a_3& a_4
\end{bmatrix} \mid a_i\in \C[[t]], a_3\in t\C[[t]], det(a)=1\}\\
&=
 (\prod_{k>0, c} X_{-\alpha_1+k\delta}(c) )
 T([[t]])
  (\prod_{k\geq 0, c} X_{\alpha_1+k\delta}(c) ). 
 \end{align*}
Similarly we can determine $I^-$ and $I^0$ in terms of their simple roots. 
For any  group element $g\in G((t))$, the conjugations $gI^+g^{-1}$,  $gI^-g^{-1}$,  $gI^0g^{-1}$ are also Iwahori subgroups. 
Then, we have the following 3 sets, each consisting a collection of Iwahori subgroups of  $G((t))$,\cite{MRY}
\begin{align*}
G((t))/I^+ 
&\text{, the thin or positive level affine flag variety,} \\
G((t))/I^0
&\text{, the semi-infinite or level 0 affine flag variety,}\\
G((t))/I^-
&\text{, the thick or negative level affine flag variety.}
\end{align*}
With the structure of Weyl group above, 
we have the following orbit decompositions
\[
G((t))=\coprod_{x\in \widehat{W}} I^+xI^+, \,\
G((t))=\coprod_{y\in \widehat{W}} I^0yI^+, \,\ 
G((t))=\coprod_{z\in \widehat{W}} I^-zI^+, 
\]
where in each decomposition the set of orbits are labelled by $\widehat{W}$.
In particular, two Iwahori subgroups from two different components are not conjugate to each other. 
A direct computation determines the conjugations of $I^+$ in terms of simple roots
\begin{align*}
&w  X_{-\alpha_1+\delta}(1) w^{-1}=\begin{cases}
X_{-\alpha_1+(2m+1)\delta}(1) &\text{if $w=(n_{s_0}n_{s_1})^m$; }\\
X_{\alpha_1+(2m+1)\delta}(-1) & \text{if  $w=n_{s_1}(n_{s_0}n_{s_1})^m $. }
\end{cases}
\\
&w  X_{\alpha_1}(1) w^{-1}=\begin{cases}
X_{\alpha_1-2m\delta}(1) &\text{if $w=(n_{s_0}n_{s_1})^m$; }\\
X_{-\alpha_1-(2m)\delta}(-1) & \text{if   $w=n_{s_1}(n_{s_0}n_{s_1})^m $. }
\end{cases}
\end{align*}
\Omit{
\begin{align*}
(n_{s_0}n_{s_1})^m  x_{-\alpha_1+\delta}(1) (n_{s_0}n_{s_1})^{-m} 
&=
(n_{s_0}n_{s_1})^m \begin{bmatrix}
1& 0\\
t &1
\end{bmatrix} (n_{s_0}n_{s_1})^{-m}\\
&=
\begin{bmatrix}
1& 0\\
t^{2m+1} \cdot 1 &1
\end{bmatrix}=x_{-\alpha_1+(2m+1)\delta}(1)
\end{align*}

\begin{align*}
(n_{s_0}n_{s_1})^m  x_{\alpha_1}(1) (n_{s_0}n_{s_1})^{-m} 
&=
(n_{s_0}n_{s_1})^m \begin{bmatrix}
1& 1\\
0&1
\end{bmatrix}  (n_{s_0}n_{s_1})^{-m}\\
&=
\begin{bmatrix}
1& t^{-2m}1\\
0&1
\end{bmatrix} =x_{\alpha_1-2m\delta}(1)
\end{align*}

\begin{align*}
n_{s_1}(n_{s_0}n_{s_1})^m  x_{-\alpha_1+\delta}(1) (n_{s_0}n_{s_1})^{-m} n_{s_1}^{-1}
&=\begin{bmatrix}
1& -t^{2m+1}1 \\
0 &1
\end{bmatrix}=x_{\alpha_1+(2m+1)\delta}(-1)
\end{align*}
\begin{align*}
n_{s_1}(n_{s_0}n_{s_1})^m  x_{\alpha_1}(1) (n_{s_0}n_{s_1})^{-m} n_{s_1}^{-1}
&=\begin{bmatrix}
1& 0\\
 -t^{-2m}1 &1
\end{bmatrix}=x_{-\alpha_1-(2m)\delta}(-1)
\end{align*}
}
Furthermore, $w  X_{-\alpha_1+\delta}(1) w^{-1}$ and $w  X_{\alpha_1}(1) w^{-1}$ determine another Iwahori subgroup $w(I^+)$ which is conjugate to $I^+$. A similar calculation can be done for $I^-$ and $I^0$.
\end{exa}

Finally we make a remark about the comparison of this Weyl group action with the braid group action on the derived category. We follow the notations above. 
\begin{rmk}\label{rmk:AffineGrassmannian}
In the example  $X=\C\times T^*\PP^1$ the root system from \S~\ref{subsec:shiftedYangian_expectation} is that of $\widehat{\fs\fl_2}$. 

We now describe a bijection between $I^+$-orbits on the thin affine flag variety $G((t))/I^+$ with the t-structures  $\mathscr{A}_m^+$ parametrized by the chambers on the PT side of the imaginary root hyperplane. This bijection is so that the simple roots of $I\in G((t))/I^+$ coincide with the simple objects in $\mathscr{A}_m^+$. To describe this bijection, we use the description of the orbits from the decomposition $G((t))=\coprod_{x\in \widehat{W}} I^+xI^+$ and the Tits extension of Weyl group elements as in Example~\ref{ex:sl2_flag}.
For simplicity, we can choose $m>0$ when $w=(n_{s_0}n_{s_1})^m$, and $m<0$, when  $w=n_{s_1}(n_{s_0}n_{s_1})^m $. (The other choice corresponds to the flop of $X$.) Therefore, the above computation shows that the simple objects 
\[
\{ch[S_0]=[C]-m[\pt], ch[S_1]=-[C]+(m+1)[\pt]\}\] in $\mathscr{A}_m^+$ match up with the simple roots 
\[
\{X_{\alpha_1-m \delta}(1),  X_{-\alpha_1+(m+1)\delta}(1)\}.
\] This bijection is equivariant with respect to the affine Weyl group actions.

Similarly, the $I^+$-orbits in $G((t))/I^-$ are in natural bijection with the $t$-structures  $\mathscr{A}_m^-$, so that the simple roots of $I\in G((t))/I^-$ coincide with the simple objects in $\mathscr{A}_m^-$.

The positive level Iwahori  and the negative ones are not related by the action of the Weyl group. Instead they differ by a homological shift $[1]$, and hence are isomorphic as abstract algebras. The  level zero Iwahori subgroup and its conjugations (as elements in $G((t))/I^0$)  do not occur here in the context of the COHA of perverse coherent sheaves.

In general, however, we do not expect two  Iwahori  subalgebras to be isomorphic in the case when they come from two  t-structures that do not differ by the action of a braid group element. A similar idea of exploring COHA's associated to two different t-structures can be found in \cite{Toda} in an attempt to ``categorify the wall-crossing".  We also mention the idea of the construction of the ``derived COHA" which contains all ``root COHA's" corresponding to  rays in ${\bf R}^2$ having the common vertex at the point $(0,0)$,  which was proposed by  Kontsevich and Soibelman in 2012 (unpublished).
\end{rmk}

\appendix

\section{Residue pushforward formula in critical cohomology}
\label{sec:appendix}
\subsection{Convolution operators on fixed points}\label{app:residue}
Let $X$ be a complex smooth algebraic $H$-variety, where $H$ is a complex affine algebraic group. We assume  that $F:=X^H$ consists of isolated fixed points $F=\{F_i \}$ and $|F|< \infty$. 
Let $i: F\inj X$ be the inclusion. 
Denote by $\varphi_f$ the vanishing cycle complex associated
to the function $f$ on $X$. 
Note that, by assumption, we have
\[
H_{c, H}^*(F, \varphi_{f\circ i}\C)^{\vee}=H_{c, H}^*(F, \C)^{\vee}=\oplus_{i} H_{c, H}^*(F_i, \C)^{\vee}
\]
The classes $i_*(1_i)$ form a basis of $H_{c, H}^*(X, \varphi_{f}\C)^{\vee}\otimes_{R_H} K_H$. 
Let us consider the map on localizations $i_*:H_{c, H}^*(F, \varphi_{f\circ i}\C)^{\vee}\otimes_{R_H} K_H \to H_{c, H}^*(X, \varphi_{f}\C)^{\vee}\otimes_{R_H} K_H$ induced by
\[
i_*: H_{c, H}^*(F, \varphi_{f\circ i}\C)^{\vee} \to H_{c, H}^*(X, \varphi_{f}\C)^{\vee}.
\]
Let $N_F X$ be the normal bundle of $F$ inside $X$, and $e(N_F X)$ be its Euler characteristic. We also have
\[
i^*i_*\alpha=e(N_{F}X), \text{and} \,\  i_*^{-1}\beta=\frac{i^*\beta}{e(N_{F}X)}
\]
for $\alpha\in H_{c, H}^*(F, \varphi_{f\circ i})^{\vee}$, and $\beta\in H_{c, H}^*(X, \varphi_{f})^{\vee}$ (see, e.g., \cite[Proposition~2.16(1)]{D}). 

Consider the correspondence
\[
\xymatrix{
&F_W\ar@{}[d]|-*[@]{\subset}^{i_W} \ar[ld]_{q_F}\ar[rd]^{p_F}&\\
F_X \ar@{}[d]|-*[@]{\subset}_{i_X} &W \ar[ld]_{q}\ar[rd]^{p}& F_Y \ar@{}[d]|-*[@]{\subset}^{i_Y}\\
X&&Y
}
\]
We have
\begin{align}\label{eqn:convolution_localization}
i_{X*}^{-1}q_* p^* i_{Y*}(1_{F_Y})
=&q_{F*} i_{W*}^{-1}p^* i_{Y*}(1_{F_Y})\\ \notag
=& \frac{i_W^*}{e(N_{F_W} W)}p^* i_{Y*}(1_{F_Y})\\ \notag
=& \frac{(p_F)^{*} i_{Y}^*}{e(N_{F_W} W)} i_{Y*}(1_{F_Y})\\ \notag
=& \frac{ (p_F)^* e(N_{F_Y}Y)}{e(N_{F_W} W)} 
\end{align}

\subsection{Jeffery-Kirwan residue formula in critical cohomology}

Here we discuss the Jeffery-Kirwan residue localization formula in the setting of critical cohomology. We follow the proof of Guillemin and Kalkman \cite{GK}, also
taking into account the description of the vanishing cycle complex of Kapranov \cite{Kap}. 
\subsubsection{De Rham model of equivariant critical cohomology}
First, we recall the de Rham model of critical cohomology following \cite{Kap}, taking into account the group action as in \cite[\S~2.10]{Ginz}.
Let $Y$ be smooth complex algebraic variety  endowed with an algebraic $\C^*$-action as well as with a $\C^*$-invariant regular function $f$. In order to use the $C^\infty$-de Rham complexes, we consider sheaves in the  analytic topology. Recall that a {\it fine sheaf} is a sheaf with partition of unity, and that a fine sheaf is acyclic under direct image functor, and hence for any sheaf $\mathcal{F}$ the higher derived direct image can be calculated by applying the direct image functor to a resolution of $\mathcal{F}$ by fine sheaves.
Recall that Poincar\'e Lemma implies that a resolution of the constant sheaf $\underline{\C}_Y$ by fine sheaves is given by the de Rham complex $(\Omega^*_Y, d)$. 
Here we use subscript notation for the sheaves, and $\Omega^*(Y)$ for the space of global sections. The space of compactly supported sections is denoted by $\Omega^*_c(Y)$. Let $\varphi_f(\underline{\C}_Y)$ denote the sheaf of vanishing cycles. A complex representing $\varphi_f(\underline{\C}_Y)$ was given in \cite{Kap}. 
Let $\tilde d:=d+df\wedge$ be the ``twisted" differential on $\Omega_Y$. The complex $(\Omega_Y, \tilde d)$ is quasi-isomorphic to $\varphi_f(\underline{\C}_Y)$.
Let $\pi_Y:Y\to\pt$ be the structure map. Recall that $Y$ is smooth and hence the Verdier dualizing sheaf $D_Y$ is $\underline{\C}_Y$ homologically shifted to the dimension of $Y$.
The cohomology of $Y$ valued in the vanishing cycle complex $H^*_{c}(Y,\varphi_f\C)^\vee:=\pi_{Y*}\varphi_f D_Y$ is then calculated as the cohomology of the complex obtained by taking global sections of the de Rham complex representing  $\varphi_f \underline{\C}_Y[\dim Y]$.

Let $\pi:X\to Y$ be a proper flat map of smooth complex algebraic varieties, and $f$ a regular function on $Y$. We have a map $\pi_*: H^*_{c}(X,\varphi_{f\circ \pi}\C)^\vee\to H^*_{c}(Y,\varphi_f\C)^\vee$ induced by applying $\pi_*$ to the map of complexes on $Y$
\[\pi_*(\Omega_X, \tilde d)\to (\Omega^*_Y, \tilde d),\]
where the individual map of sheaves $\pi_*\Omega^*_X\to \Omega^*_Y$ is given by integration along fibers of $\pi$.

The algebraic $\C^*$-action induces an action of the maximal torus $S^1\subseteq \C^*$ when $Y$ is considered as a smooth manifold.  
In the calculations below we follow \cite[\S~2.10]{Ginz} taking into account the $S^1$-action via the Cartan complex. Note that the Borel construction used in \cite{D} differs from the Cartan model by a completion. 
The cohomology $H^*_{c,S^1}(Y,\varphi_f\C)^\vee$
can be calculated via the de Rham model $\Omega^*_{S^1}(Y,\varphi_f):=\Omega^*(Y)^{S^1}\otimes \C[x]$ with the differential $\tilde{d}=d+i(v)\otimes x+\wedge df$. Here $d$ is the de Rham differential on forms, $v$ is the vector field induced by the $S^1$-action, and $d_{S^1}:=d+i(v)\otimes x$ is the differential calculating the usual  Cartan model of equivariant cohomology. 

\subsubsection{An explicit formula of the Kirwan map}

Let $\overline{X}$ be a smooth complex algebraic variety with a $\C^*$-action of dimension $n$. 
In order to consider de Rham complex of the GIT quotient $\overline{X}/\!/_\xi \C^*$ for a character $\xi$ of $\C^*$, we use Kirwan's theorem to identify $\overline{X}/\!/_\xi \C^*$ with the Hamiltonian reduction. 
More precisely, the GIT quotient $\overline{X}/\!/_\xi \C^*$ is diffeomorphic to $\mu^{-1}(\xi)/S^1$ where $\mu$ is the moment map of the compact Lie group $S^1\subseteq \C^*$. 
Let $X=\mu^{-1}([\xi,\infty))$, which is a smooth manifold with boundary and an $S^1$-action that preserves the boundary.

Abusing the notation,  for any smooth manifold $N$ with boundary, endowed with an $S^1$-action as well as  with a smooth function $f$ which preserves the boundary and $f$, we denote the complex $\Omega^*(N)^{S^1}\otimes \C[x]$ with the differential $\tilde{d}=d+i(v)\otimes x+\wedge df$  by $\Omega^*_{S^1}(N,\varphi_f)$.

In the setup above,  let $X$ be 
 a smooth manifold with boundary endowed with a smooth $S^1$-action. Assume the $S^1$-action on $\partial X$ is locally free, so that $\partial X/S^1\cong M$ is again a smooth algebraic variety. Let $f$ be an $S^1$-invariant  regular function on $\overline{X}$. Then, the restriction of $f$ to $\partial X$ is an $S^1$-invariant  smooth function, which is again obtained from pulling back a regular function on $M$, which by an abuse of notation we still denote by $f$.

We define the Kirwan map $\kappa: H^*_{S^1,c}(X,\varphi_f)^\vee\to H^*_{c}(M,\varphi_f)^\vee$ to be the following map at the level of de Rham complexes. We have the restriction of forms $i^*: \Omega^*_{S^1}(X,\varphi_f)\to\Omega^*_{S^1}(\partial X,\varphi_f)$, as well as the pullback $\pi^*:\Omega^*(M,\varphi_f)\to\Omega^*_{S^1}(\partial X,\varphi_f)$, where $\pi: \partial X\to M$ is the quotient map.  Then, for any closed form $\alpha\in \Omega^*_{S^1}(X,\varphi_f)$,  $\kappa(\alpha)$  is the class of any form $\gamma\in \Omega^*(M,\varphi_f)$ with the property that $\pi^*(\gamma)$ has the same class as  $i^*\alpha$. That is, $i^*\alpha-\pi^*(\gamma)$ is a boundary cycle. The calculation of Guillemin and Kalkman \cite{GK} gives an explicit formula for $\gamma$, which we recall here. 

Let $\theta$ be a $S^1$-invariant one form such that $i(v)\theta=1$, which is well-defined on the complement of $X^{S^1}$.
Consider the formal expression $\nu_0=\frac{\theta}{x+d\theta+\theta\wedge df}$.  From definition, we have $\tilde d (\theta)=x+d\theta+\theta\wedge df$, in particular $\tilde d(x+d\theta+\theta\wedge df)=0$. Therefore, $\tilde d(\nu_0)=\frac{\tilde d (\theta)}{x+d\theta+\theta\wedge df}=1$.
Let $\alpha\in \Omega^*_{S^1}(X,\varphi_f)_c$ such that $\tilde d (\alpha)=0$. Then, $\tilde{d}(\alpha\wedge \nu_0)=\alpha$. 
We consider $\nu_0=\frac{\theta}{x}\sum_{n\geq 0}(\frac{-d\theta-\theta\wedge df}{x})^n$ as a Laurent power series in $x$ with coefficients in $\Omega(X)^{S^1}$.

Now let $\{X_i\}_{i=1,\dots, N}$ be the connected components of the fixed point set, and let $U_i$ be pairwise disjoint tubular neighbourhoods of $X_i$ in $X$ such that $U_i\cap \partial X=\emptyset$.  Let $i^*_k: X_k\to X$ be the embedding.
Assume $\deg\alpha=\dim \partial X-1$. Then, by Stokes theorem, we have \[0=\int_X\alpha=\int_X\tilde d \nu_0\wedge\alpha=\sum_{k=1}^N\int_{\partial U_k}\frac{\theta \alpha}{x+d\theta+\theta\wedge df}+\int_{\partial X}\frac{\theta \alpha}{x+d\theta+\theta\wedge df}.\]
We confine ourselves to the case when $X^{S^1}=\sqcup X_k$ is a finite set.
We now follow the calculation of \cite[\S~10.8.6]{GS} to show that $\int_{\partial U_k}\frac{\theta \alpha}{x+d\theta+\theta\wedge df}$ converges to $\frac{i^*_k\alpha}{e(\nu_k)}$ where $\nu_k$ is the tangent space of $X_k$ and $e(\nu_k)$ is the equivariant Euler class, which can be calculated as the product of the weights with multiplicities on $\nu_k$. Indeed, as in {\it loc. cit.}, we have $\alpha=f(x)+d_{S^1}\beta$ where $f(x)$ is the restriction of $\alpha$ at $X_k$ and $\beta\in\Omega_{S^1}(U_k)$, hence $\nu_0\wedge\alpha=f(x)\nu_0+\beta+d_{S^1}(\nu_0\wedge\beta)$. Integrating on $\partial U_k$, the second term vanishes;  the third term is of higher order in terms of the radius of $U_k$; the first term becomes $f(x)\int_{\partial U_k}\nu_0$, where the calculation in local coordinate in {\it loc. cit. } shows that $\int_{\partial U_k}\nu_0=\frac{1}{x^n e(\nu_k)}$. 
To summarize, we obtain 
\begin{equation}\label{eqn:localization}
\int_{\partial X}\frac{\theta \alpha}{x+d\theta+\theta\wedge df}=\sum_{k=1}^N\frac{ i^*_k\alpha}{x^ne(\nu_k)}
\end{equation}
under the assumption that $\deg\alpha=\dim \partial X-1$. By degree consideration (see \cite[page 21]{Z}), the only powers of $x$ that can occur in the Laurent series of $\nu_0$ are the
positive ones and $x^{-1}$. So we have $\alpha\wedge\nu_0=\nu+\beta x^{-1}$, where $\nu$ has the positive powers of $x$ and  $\beta=\Res_{x=0}(\alpha\wedge\nu_0)$.
\footnote{There is a typo in \cite{Z}, where the residue should be over $x=0$. }
In particular, $\alpha=\tilde d( \nu)+i(v)\beta$ with $i(v)\beta=\pi^*\gamma$ for some $\gamma\in \Omega^*(M,\varphi_f)$.
In other words, 
\begin{equation}\label{eqn:residue}
\kappa(\alpha)=\Res_{x=0}\pi_*\frac{\theta \alpha}{x+d\theta+\theta\wedge df}.
\end{equation}
 
\subsubsection{Conclusion}
Assume that $X$ is endowed with an action of $S\times S^1$, where $S$ is a compact Lie group which acts in a Hamiltonian way. Let $\{E_m\}$  be a system of $S$-spaces the limit of which is the Borel construction of $ES$. Following the consideration in \cite[\S~5.2]{Z}, for any  $\alpha \in H^*_{c,S\times S^1}(X,\varphi_f\C)^\vee$ of a certain degree $i$, we choose a representative in $\Omega^{i}_{S^1}(E_m,\varphi)$ with $i=\dim X_m-1$. In the argument above, replacing $X$ by $X_m:=E_m\times_SX$, taking residue from both sides of \eqref{eqn:localization} and taking into account \eqref{eqn:residue},
we obtain the residue pushforward formula in $S$-equivariant cohomology
\begin{equation}\label{eqn:JefferyKirwan}
\int_X \kappa(\alpha)=\sum_{k=1}^N\Res_{x=0}\frac{i^*_k\alpha}{e(\nu_k)}
\end{equation} 
in $H^*_{c,S}(X,\varphi_f\C)^\vee\to H^*_{c,S}(\pt)^\vee$.

{\bf The authors:}

M.~Rap\v{c}\'{a}k- 
{Center for Theoretical Physics, University of California, Berkeley, CA 	94720 USA},
{miroslav.rapcak@gmail.com}

Y.~Soibelman-
{Department of Mathematics, Kansas State University, Manhattan, KS 66506, USA},
{soibel@math.ksu.edu}

Y.~Yang- {School of Mathematics and Statistics, The University of Melbourne, 813 Swanston Street, Parkville VIC 3010, Australia},
{yaping.yang1@unimelb.edu.au}

 G.~Zhao-
{School of Mathematics and Statistics, The University of Melbourne, 813 Swanston Street, Parkville VIC 3010, Australia},
{gufangz@unimelb.edu.au}

\end{document}